\documentclass[a4paper,11pt,fleqn]{article}
\usepackage{amsmath,amssymb,amsthm,graphicx,subfigure,float,caption,epstopdf,tabularx,color, bm,amsfonts,epic}
\usepackage[top=1in, bottom=1in, left=1.25in, right=1.25in]{geometry}
\usepackage{appendix}
\allowdisplaybreaks
\newtheorem{thm}{Theorem}[section]
\newtheorem{lem}{Lemma}[section]

\newtheorem{rmk}{Remark}[section]
\newtheorem*{prf}{Proof}
\numberwithin{equation}{section}

\newcommand*{\m}{\phantom{m}}

\usepackage{indentfirst}
\graphicspath{{Fig}}

\begin{document}
\title{A Novel Sixth Order Energy-Conserved Method for
Three-Dimensional Time-Domain Maxwell's Equations}
\author{Chaolong Jiang$^1$, \quad Wenjun Cai$^1$, \quad Yushun Wang$^{1,}$\footnote{Correspondence author. Email:
wangyushun@njnu.edu.cn.},  \quad Haochen Li$^2$\\
{\small $^1$ Jiangsu Provincial Key Laboratory for NSLSCS,}\\
{\small School of Mathematical Sciences,  Nanjing Normal University,}\\
{\small  Nanjing 210023, China}\\
{\small $^2$ LMAM, CAPT and School of Mathematical Sciences, }\\
{\small Peking University, Beijing 100871, China}\\
}
\date{}
\maketitle

\begin{abstract}
In this paper, a novel sixth order energy-conserved method is proposed for solving
the three-dimensional time-domain Maxwell's equations.
The new scheme preserves five discrete energy conservation laws,
 three momentum conservation laws, symplectic conservation law
as well as two divergence-free properties and is proved to be unconditionally stable, non-dissipative.
An optimal error estimate is established based on
the energy method, which shows that the proposed method
is of sixth order accuracy in time and spectral accuracy in space in discrete $L^{2}$-norm.
The constant in the error estimate is proved to be only $O(T)$.
Furthermore, the numerical dispersion relation is analyzed in detail and a fast solver is presented to solve the resulting discrete linear equations efficiently. Numerical results are addressed to verify our theoretical analysis.
\\[2ex]
\end{abstract}

\section{Introduction}
\label{sec1}
The Maxwell's equations describe the propagation and scattering of electromagnetic
waves and have a wide variety of applications in science and engineering, including
microwave circuits, radio-frequency, antennas, aircraft radar,
 integrated optical circuits, wireless engineering, etc.
Various applications stimulate the investigation on the construction
of efficient numerical methods for the Maxwell's equations.
A well-known numerical method in computational electromagnetic is the finite-difference time domain
(FDTD) method, which was first introduced by Yee \cite{KSY66} and further developed and analyzed in Refs. \cite{MS94,AT96}.
However, the Yee-based FDTD method is only conditionally stable
so that it may require very small temporal step-size
and suffer impractical computational cost for long time computation.

In recent years, due to 
the superior properties
in long time numerical computation over traditional numerical methods, structure-preserving methods have been proved to
be very powerful in numerical simulations (e.g., see Refs. \cite{FQ10,ELW06,LBR04,MPS98} and references therein).
As one of the most important components of the structure-preserving methods,
symplectic schemes, which can preserve the sympletic conservation law of the Maxwell's equations,
have gained remarkable success in the numerical analysis of
the Maxwell's equations (e.g., see Refs. \cite{CWS13,HLSY01,KHZ10,SCL16,SQS07,ST11,ZSC11} and references therein).
In addition to the symplectic conservation law,
the Maxwell's equations also admit two
divergence-free fields, five
energy conservation laws and three momentum conservation laws,
which are very important invariants for
 long time propagation of the electromagnetic waves \cite{AT96}.
Thus, devising the numerical schemes, which can inherit these
original physical features as much as possible, attracts a lot of interest.
In Ref. \cite{CLL08}, Chen et al. proposed an energy-conserved splitting scheme for two dimensional (2D) Maxwell's equations in an
isotropic, lossless and sourceless medium. Further analysis in three dimensional
(3D) case was investigated in Ref. \cite{CLL10}.
Other energy-conserved splitting schemes can be found in Refs. \cite{HJK14,KHZ10,LY13}.
However, most of existing energy-conserved splitting schemes have only second order accuracy in both time and space at most.
The second order schemes are usually effective for
geometries of moderate electrical size, but, for computing large scale problems,
for problems requiring long-time integration, or
for problems of wave propagations over longer distances, it is indispensable for the development of higher order schemes
 which produce smaller dispersion or phase errors for a given mesh resolution \cite{LY13}.

In fact, various spatial high order methods have been proposed for the Maxwell's
 equations in the literature, such as the high order FDTD methods (e.g., Refs. \cite{Shang99,STDM15,AT96,YP01,ZW04}), the Finite element method (see the monograph \cite{MP03}) and the Fourier pseudo-spectral method \cite{LQH97}.
 However, most of these methods have only second order accuracy in time to the best of our knowledge.
In Ref. \cite{CHWG15}, with the aid of the splitting techniques, an unconditionally stable energy-conserved scheme was presented.
This scheme has fourth order accuracy in time and spectral accuracy in space, respectively.
 Higher order methods can be constructed by using the splitting techniques, 
however, this approach can be tedious in practice computation and error estimate.
In addition, the scheme cannot preserve the two discrete divergence-free fields. 
The resulting large divergence-free
errors may make the predicted solutions less stable in long time simulations of
the electromagnetic wave propagation. Avoidance of this numerical
instability is critical in the development of an effective
numerical method for the Maxwell's equations \cite{CLS04,MOS00}.
Subsequently, another energy-conserved scheme with fourth order accuracy in time and spectral accuracy in space, respectively, was designed in Ref. \cite{CWG16}.
The scheme can preserve the two discrete divergence-free fields of the Maxwell's equations.
However, there has been no reference considering a sixth order energy-conserved method for the 3D Maxwell's equations.

Thus, in this paper, our main contributions are as follows:
\begin{itemize}
\item[1.] We first propose a sixth order energy-conserved scheme for the three-dimensional time-domain Maxwell's equations, which satisfies the following properties.
\begin{itemize}
\item[-] The proposed scheme is unconditionally stable and non-dissipative (numerical dispersion analysis).
\item[-] The proposed scheme is symmetric and can preserve five energy, three momentum conservation laws, symplectic conservation law as well as two discrete divergence-free fields.
\end{itemize}
\item[2.] An optimal error estimate of the proposed scheme is established, which shows that the proposed scheme
is unconditionally convergent with order of $O(\tau^{6}+N^{-r})$ in discrete $L^{2}$-norm,
where $\tau$ is the time step size and $N$ is the collocation points used in the spectral method.
In particular, the constant in the error estimate is only $O(T)$.
\item[3.] The numerical dispersion relation including the phase velocity and the group velocity of the scheme is
analyzed in detail, which gives two novel results:
\begin{itemize}
\item[-] The extra solution branches of the Fourier pseudo-spectral scheme only occur with large wave numbers.
\item[-] The grid-anisotropy of the Fourier pseudo-spectral scheme is direction-independent.
\end{itemize}
\item[4.] A fast solver is developed for solving the resulting linear equations efficiently, which can be extended directly to other literature (e.g., Refs. \cite{CWG16,ZSC11}).
\end{itemize}

The outline of this paper is
organized as follows. In Section 2, the Hamiltonian structures of the Maxwell's equations
 as well as the conservation laws in continuous sense is introduced. A sixth order energy-conserved scheme for the Maxwell's equations is proposed in Section 3.
The discrete conservation laws including the five energy conservation laws, the three momentum conservation laws, the symplecticity
and the two divergence-free fields are rigorously proved in Section 4.
The convergence is analyzed in Section 5.
Numerical dispersion relation is investigated in Section 6.
In Section 7, numerical experiments are
presented. We draw some conclusions in Section 8. Finally, a fast solver is presented in Appendix.

\section{Hamiltonian structures and conservation laws}
In this section, we will introduce the Hamiltonian structures of the Maxwell's equations
 as well as the conservation laws in continuous sense.

The 3D time-domain Maxwell's equations in an
isotropic, lossless and sourceless medium can be written as
\begin{align}\label{ME}
 \left \{
 \aligned
 &\frac{\partial {\bf E}}{\partial t}=\frac{1}{\epsilon}\nabla\times {\bm H}, \\
 &\frac{\partial {\bm H}}{\partial t}=-\frac{1}{\mu}\nabla\times {\bm E}, \\
 &\nabla\cdot(\epsilon{\bm E})=0, \\
 &\nabla\cdot(\mu{\bm H})=0,  
 \endaligned
 \right.
 \end{align}
 where ${\bm E}=(E_{x},E_{y},E_{z})^{T}$ is the electric field intensity,
 ${\bm H}=(H_{x},H_{y},H_{z})^{T}$ is the magnetic field intensity,
constant scalars $\mu$ and $\epsilon$ are the magnetic permeability
and the electric permittivity, respectively.
Here, we assume the periodic boundary conditions (PBCs) on the boundary
of the cuboid domain $\Omega=[x_{L},x_{R}]\times[y_{L},y_{R}]\times[z_{L},z_{R}]$.
The initial conditions are supposed to be
\begin{align}
{\bm E}(x,y,z,0)={\bm E}_{0}(x,y,z),\ {\bm H}(x,y,z,0)={\bm H}_{0}(x,y,z),\ (x,y,z) \in\Omega.
\end{align}

In general, the perfect electric conducting (PEC) boundary conditions are used on $\Omega\times[0,T]$.
However, with the aid of the B\'erenger's PML technique \cite{Berenger94,Berenger96}, the
periodic boundary conditions are valid in practical computations. More importantly, when the periodic boundary conditions are considered,
the spatial derivatives of the Maxwell's equations can be efficiently computed by using the fast Fourier transform (FFT).
For more details, please refer to Ref. \cite{LQH97}.

 The Maxwell's equations (\ref{ME}) are a bi-Hamiltonian system,
 which implies that there are two infinite-dimensional Hamiltonian formulations.

The first one is given by
\begin{align}\label{MEH1}
\frac{\partial}{\partial t}\left(\begin{array}{c}
              {\bm H} \\
              {\bm E}  \\
             \end{array}
\right)=\left(\begin{array}{cc}
              0     &  -{\bm I}_{3\times 3}  \\
              {\bm I}_{3\times 3}    &  0 \\
             \end{array}
\right)\left(\begin{array}{c}
              \epsilon^{-1}\nabla\times {\bm H} \\
              \mu^{-1}\nabla\times {\bm E} \\
             \end{array}
\right)=\mathcal{J}_{1}\left(\begin{array}{cc}
              \frac{\delta \mathcal{H}_{1}}{\delta {\bm H}} \\
              \frac{\delta \mathcal{H}_{1}}{\delta {\bm E}} \\
             \end{array}
\right),
\end{align}
with the helicity Hamiltonian functional \cite{AA83}
\begin{align}\label{MEHF1}
\mathcal{H}_{1}=\int_{\Omega} \left(\frac{1}{2\epsilon}{\bm H}^{T}
(\nabla\times {\bm H})+\frac{1}{2\mu}{\bm E}^{T}(\nabla\times {\bm E})\right)dxdydz.
\end{align}
The other one reads
\begin{align}\label{MEH2}
\frac{\partial}{\partial t}\left(\begin{array}{c}
              {\bm H} \\
              {\bm E}  \\
             \end{array}
\right)=\left(\begin{array}{cc}
              0     &  -(\epsilon \mu)^{-1}\nabla \times  \\
              (\epsilon \mu)^{-1}\nabla \times    &  0 \\
             \end{array}
\right)\left(\begin{array}{c}
              \mu{\bm H}\\
             \epsilon{\bm E} \\
             \end{array}
\right)=\mathcal{J}_{2}\left(\begin{array}{cc}
              \frac{\delta \mathcal{H}_{2}}{\delta {\bm H}} \\
              \frac{\delta \mathcal{H}_{2}}{\delta {\bm E}} \\
             \end{array}
\right),
\end{align}
where the quadratic Hamiltonian functional yields \cite{MW82}
\begin{align}\label{MEHF2}
\mathcal{H}_{2}=\int_{\Omega} \left(\frac{\mu}{2}{\bm H}^{T}
 {\bm H}+\frac{\epsilon}{2}{\bm E}^{T}{\bm E}\right)dxdydz,
\end{align}
which is the electromagnetic energy in the Poynting theorem in classical electromagnetism \cite{JK98}.

With the Hamiltonian formulations (\ref{MEH1}) and (\ref{MEH2}), the Maxwell's equations (\ref{ME})
satisfy the symplectic, the helicity and quadratic conservation laws.
\begin{lem} \cite{ST11} Let ${\bm H}$ and ${\bm E}$ be the solutions of the
Maxwell's equations (\ref{ME}) and satisfy PBCs, then it holds that
\begin{align}
\frac{d}{dt}\int_{\Omega} \omega(x,y,z,t)dxdydz=0, \ \omega=dE_{x}\wedge dH_{x}+dE_{y}\wedge dH_{y}+dE_{z}\wedge d H_{z}.
\end{align}
\end{lem}
\begin{lem} \cite{ST11} Let ${\bm H}$ and ${\bm E}$ be the solutions of the
Maxwell's equations (\ref{ME}) and satisfy PBCs, then the helicity and quadratic conservation laws hold, i.e.,
\begin{align}
&\frac{d}{dt}\int_{\Omega} \left(\frac{1}{2\epsilon}{\bm H}^{T}(\nabla\times {\bm H})
+\frac{1}{2\mu}{\bm E}^{T}(\nabla\times {\bm E})\right)dxdydz=0,\\ &\frac{d}{dt}\int_{\Omega}
 \left(\frac{\mu}{2}{\bm H}^{T} {\bm H}+\frac{\epsilon}{2}{\bm E}^{T} {\bm E}\right)dxdydz=0.
\end{align}
\end{lem}
Additionally, the Maxwell's equations (\ref{ME}) also admit others conservation laws.
\begin{lem} \cite{ST11} Let ${\bm H}$ and ${\bm E}$ be the solutions of the
Maxwell's equations (\ref{ME}) and satisfy PBCs, then there exit the momentum conservation laws
\begin{align}\label{MCL}
&\frac{d}{dt}\int_{\Omega}{\bm H}^{T} \partial_{w}{\bm E}dxdydz=0,
\end{align}
where $w=x,y,$ or $z$, and hereafter.
\end{lem}
\begin{lem} \cite{GZ13}  Let ${\bm H}$ and ${\bm E}$ be the solutions of the
Maxwell's equations (\ref{ME}) and satisfy PBCs, then there admit the energy conservation laws
\begin{align}\label{CECL-III}
&\frac{d}{dt}\int_{\Omega}\left(\frac{\epsilon}{2}\partial_{t}{\bm E}^{T}\partial_{t}{\bm E}+\frac{\mu}{2} \partial_{t}{\bm H}^{T}\partial_{t}{\bm H}\right)dxdydz
=0,\\
&\frac{d}{dt}\int_{\Omega}\left(\epsilon\partial_{w}{\bm E}^{T}\partial_{w}{\bm E}+\mu \partial_{w}{\bm H}^{T}\partial_{w}{\bm H}\right)dxdydz
=0,\\
&\frac{d}{dt}\int_{\Omega}\left(\epsilon\partial_{tw}{\bm E}^{T}\partial_{tw}{\bm E}+\mu \partial_{tw}{\bm H}^{T}\partial_{tw}{\bm H}\right)dxdydz
=0.
\end{align}
\end{lem}
\section{A sixth order  energy-conserved scheme}
In this section, our main goal is to propose a sixth order energy-conserved method for the Maxwell's equations in time.
Let $\Omega_{h}=\{(x_{j},y_{k},z_{m})|x_{j}
=x_{L}+(j-1)h_{x},\ y_{k}=y_{L}+(k-1)h_{y},\  z_{m}=z_{L}+(m-1)h_{z};\  j=1,\cdots,N_{x},\ k=1,\cdots,N_{y},
m=1,\cdots,N_{z}\}$ be a partition of $\Omega$ with the mesh size $h_{w}=\frac{L_{w}}{N_{w}}$, where $N_{w}$ is an even integers and $L_{w}=w_{R}-w_{L}$. Denote $h=\text{max}\{h_x,h_y,h_z\}$. Let $\Omega_{\tau}=\{t_{n}|t_{n}=n\tau; 0\leqslant n \leqslant M\}$ be a uniform partition of $[0,T]$ with the time step
  $\tau=\frac{T}{M}$ and $\Omega_{h\tau}=\Omega_{h}\times\Omega_{\tau}$.

Let $\{U_{j,k,m}^{n}|j=1,\cdots,N_{x},\ k=1,\cdots,N_{y},
m=1,\cdots,N_{z},\ 0\leqslant n \leqslant M\}$ be a mesh functions defined on $\Omega_{h\tau}$. Some notations are
introduced below.
\begin{align*}
\hat{\delta}_{t}U_{j,k,m}^{n-1/2}=\frac{U_{j,k,m}^{n}-U_{j,k,m}^{n-1}}{\tau},\  U_{j,k,m}^{n+\frac{1}{2}}=\frac{U_{j,k,m}^{n+1}+U_{j,k,m}^{n}}{2}.
\end{align*}
Moreover, for any grid functions $U_{j,k,m},V_{j,k,m},(x_j,y_k,z_m)\in\Omega_h$, we define the inner product and the norm by
\begin{align*}
\langle{\bm U},{\bm V}\rangle_{h}=h_{x}h_{y}h_{z}\sum_{j=1}^{N_{x}}\sum_{k=1}^{N_{y}}\sum_{m=1}^{N_{z}}U_{j,k,m}
\bar{V}_{j,k,m},\ ||{\bm U}||_{h}^{2}=\langle{\bm U},{\bm U}\rangle_{h},\ ||{\bm U}||_{h,\infty}=\max\limits_{j,k,m
} |U_{j,k,m}|,
\end{align*}
where $\bar{V}_{j,k,m}$ denotes the conjugate of $V_{j,k,m}$.
For vectors ${\bm U}=[({\bm U}_{x})^{T},({\bm U}_{y})^{T},({\bm U}_{z})^{T}]^{T}$
and ${\bm V}=[({\bm V}_{x})^{T},({\bm V}_{y})^{T},({\bm V}_{z})^{T}]^{T}$, the corresponding inner product and norm are
\begin{align*}
&\langle{\bm U},{\bm V}\rangle_{h}=\langle{\bm U}_{x},{\bm V}_{x}\rangle_{h}+\langle{\bm U}_{y},{\bm V}_{y}\rangle_{h}+\langle{\bm U}_{z},
{\bm V}_{z}\rangle_{h},\\
&||{\bm U}||_{h}^{2}=\langle{\bm U}_{x},{\bm U}_{x}\rangle_{h}+\langle{\bm U}_{y},{\bm U}_{y}\rangle_{h}+\langle{\bm U}_{z},{\bm U}_{z}\rangle_{h}.
\end{align*}
\subsection{Fourier pseudo-spectral approximation in space}
As achieving high order accuracy in time, we should also treat the space with an appropriate discretization.
The Fourier pseudo-spectral method is a very good candidate because of the high order accuracy and the fast Fourier
 transform (FFT) algorithm.  
 Moreover, it is shown that the Fourier pseudo-spectral method exhibits
obvious superiority over the conventional finite difference method in simulating electromagnetic
waves \cite{LQH97}.

Now, let us introduce the following Fourier pseudo-spectral discretization for the space variables.
We define
\begin{align*}
&S_{N}^{'''}=\text{span}\{g_{j}(x)g_{k}(y)g_{m}(z),\ j=1,\cdots,N_{x},\ k=1,\cdots,N_{y},\ m=1,\cdots,N_{z}\},
\end{align*}
as the interpolation space, where $g_j(x)$, $g_k(y)$ and $g_m(z)$ are trigonometric polynomials of degree $N_{x}/2$, $N_{y}/2$ and $N_{z}/2$,
given respectively by
\begin{align*}
  &g_{j}(x)=\frac{1}{N_{x}}\sum_{l=-N_{x}/2}^{N_{x}/2}\frac{1}{a_{l}}e^{\text{i}l\mu_{x} (x-x_{j})},\ g_{k}(y)=\frac{1}{N_{y}}\sum_{s=-N_{y}/2}^{N_{y}/2}\frac{1}{b_{s}}e^{\text{i}s\mu_{y} (y-y_{k})},\nonumber \\
  &g_{m}(z)=\frac{1}{N_{z}}\sum_{q=-N_{z}/2}^{N_{z}/2}\frac{1}{c_{q}}e^{\text{i}q\mu_{z} (z-z_{m})},
\end{align*}
with $a_{l}=\left \{
 \aligned
 &1,\ |l|<\frac{N_x}{2},\\
 &2,\ |l|=\frac{N_x}{2},
 \endaligned
 \right.$, $b_{p}=\left \{
 \aligned
 &1,\ |s|<\frac{N_y}{2},\\
 &2,\ |s|=\frac{N_y}{2},
 \endaligned
 \right.$,\ $c_{q}=\left \{
 \aligned
 &1,\ |q|<\frac{N_z}{2},\\
 &2,\ |q|=\frac{N_z}{2},
 \endaligned
 \right.$ and $\mu_{w}=\frac{2\pi}{L_{w}}$.

We define the interpolation operator $I_{N}: C(\Omega)\to S_{N}^{'''}$ as follows
\begin{align}\label{IC}
I_{N}U(x,y,z,t)=\sum_{j=1}^{N_{x}}\sum_{k=1}^{N_{y}}\sum_{m=1}^{N_{z}}U_{j,k,m}(t)g_{j}(x)g_{k}(y)g_{m}(z),
\end{align}
where $U_{j,k,m}(t)=U(x_{j},y_{k},z_{m},t)$ and its vector form is denoted by
\begin{align*}
&{\bm U}=(U_{1,1,1},\cdots,U_{N_{x},1,1},U_{1,2,1},\cdots, U_{N_{x},2,1},\cdots,U_{1,N_y,N_{z}},\cdots,U_{N_x,N_{y},N_{z}})^{T}.
\end{align*}
Making partial differential with respect to $x$, $y$ and $z$, respectively, and evaluating the resulting expression at collocation points ($x_{j},y_{k},z_{m}$), we can obtain
\begin{align*}
\frac{\partial^{p} I_{N}U(x_{j},y_{k},z_{m},t)}{\partial x^{p}}
=&\sum_{j^{'}=1}^{N_{x}}\sum_{k^{'}=1}^{N_{y}}\sum_{m^{'}=1}^{N_{z}}U_{j^{'},k^{'},m^{'}}(t)\frac{d^{p}g_{j^{'}}(x_{j})}{dx^{p}}g_{k^{'}}(y_{k})g_{m^{'}}(z_{m})\\
=&[({\bm I}_{N_{z}}\otimes {\bm I}_{N_{y}}\otimes{\bm D}_{p}^{x}){\bm U}]_{N_{x}N_{y}(m-1)+N_{x}(k-1)+j},\\
\frac{\partial^{p} I_{N}U(x_{j},y_{k},z_{m},t)}{\partial y^{p}}
=&\sum_{j^{'}=1}^{N_{x}}\sum_{k^{'}=1}^{N_{y}}\sum_{m^{'}=1}^{N_{z}}U_{j^{'},k^{'},m^{'}}(t)g_{j^{'}}(x_{j})\frac{d^{p}g_{k^{'}}(y_{k})}{dy^{p}}g_{m^{'}}(z_{m})\\
=&[({\bm I}_{N_{z}}\otimes {\bm D}_{p}^{y}\otimes {\bm I}_{N_{x}}){\bm U}]_{N_{x}N_{y}(m-1)+N_{x}(k-1)+j},\\
 \frac{\partial^{p} I_{N}U(x_{j},y_{k},z_{m},t)}{\partial z^{p}}
=&\sum_{j^{'}=1}^{N_{x}}\sum_{k^{'}=1}^{N_{y}}\sum_{m^{'}=1}^{N_{z}}U_{j^{'},k^{'},m^{'}}(t)g_{j^{'}}(x_{j})g_{k^{'}}(y_{k})\frac{d^{p}g_{m^{'}}(z_{m})}{dz^{p}}\\
=&[({\bm D}_{p}^{z}\otimes {\bm I}_{N_{y}}\otimes {\bm I}_{N_{x}}){\bm U}]_{N_{x}N_{y}(m-1)+N_{x}(k-1)+j},
\end{align*}
where $\otimes$ is the Kronecker product, ${\bm I}_{N_{w}}$ is the identity matrix of dimension $N_{w}\times N_{w}$ and ${\bm D}_{p}^{w}\in \mathbb{R}^{N_{w}\times N_{w}}$ is the
spectral differential matrix whose
entries are
\begin{align}
({\bm D}_{p}^{w})_{j,l}=\frac{d^{p}g_{l}(w_{j})}{dw^{p}}.
\end{align}
with $j,l=1,\cdots,N_w$. In fact, by careful calculations, one can obtain explicitly
\begin{align*}
&({\bm D}^{x}_{1})_{j,l}=
 \left \{
 \aligned
 &\frac{1}{2}\mu_{x} (-1)^{j+l} \cot(r_{x}),\ &j\neq l,\\
 &0,\quad \quad \quad \quad \quad \quad \quad \ \ \ ~&j=l,
 \endaligned
 \right.
\end{align*}
\begin{align*}
 &({\bm D}^{x}_{2})_{j,l}=
 \left \{
 \aligned
 &\frac{1}{2}\mu_{x}^{2} (-1)^{j+l+1}\csc^{2}(r_{x}),\ &j\neq l,\\
 &-\mu_{w}^{2}\frac{N_{w}^{2}+2}{12},\quad \quad \quad \quad ~ &j=l,
 \endaligned
 \right.
 \end{align*}
\begin{align*}
 &({\bm D}^{x}_{3})_{j,l}=
 \left \{
 \aligned
 &\frac{3\mu_{x}^{3}}{4} (-1)^{j+l}\cos(r_{x})\csc^{3}(r_{x})+\frac{\mu_{x}^{3}N_{x}^{2}}{8}(-1)^{j+l+1}\cot(r_{x}),\ &j\neq l,\\
 &0,\quad \quad \quad \quad \quad \quad \quad~ ~~~~~~~~~~~~~~~~~~~~~~~~~~~~~~~~~~~~~~~~~~~~~~~~&j=l,
 \endaligned
 \right.
 \end{align*}

 \begin{align*}
 &({\bm D}^{x}_{4})_{j,l}=
 \left \{
 \aligned
 &\mu_{w}^{4}(-1)^{j+l}\csc^{2}(r_{x})\Big(\frac{N_{x}^2}{4}-\frac{1}{2}-\frac{3}{2}\cot^{2}(r_{x})\Big),\ &j\neq l,\\
 &\mu_{w}^{4}(\frac{N_{x}^{4}}{80}+\frac{N_{x}^{2}}{12}-\frac{1}{30}),\quad \quad \quad \quad \quad ~~~~~~~~~~~~~ &j=l,
 \endaligned
 \right.
 \end{align*}
and
 \begin{align*}
& ({\bm D}^{x}_{5})_{j,l}=
 \left \{
 \aligned
 &\frac{\mu_{x}^{5}}{32}(-1)^{j+l}\cot(r_{x})\Big[N_{x}^{4}+20\csc^{2}(r_{x})\Big(4+6\cot^{2}(r_{x})-N_{x}^{2}\Big)\Big],\ &j\neq l,\\
 &0,\quad \quad \quad ~~~~~~~~~~~~~~~~~~~~~~~~~~~~~~~~~~~~~~~~~~~~~~~~~~~~~~~~~~~~~~~~~~~ &j=l,
 \endaligned
 \right.
 \end{align*}
where $r_{x}=\mu_{x} \frac{x_{j}-x_{l}}{2}$ and $j, l=1,\cdots, N_{x}$.

Furthermore, the following relationship holds \cite{CQ01}
\begin{align}\label{CQ}
({\bm D}_{p}^{x})_{j,l}=(({\bm D}_{1}^{x})^{p})_{j,l}+(-1)^{j+l}\frac{\mu_{x}^{p}}{2N_{x}}
\left[\left(\text{i}\frac{N_{x}}{2}\right)^{p}+\left(-\text{i}\frac{N_{x}}{2}\right)^{p}\right],
\end{align}
which implies that ${\bm D}_{p}^{x}=({\bm D}_{1}^{x})^{p}$, if $p$ is an odd integer.
Here, we only give the case of $w=x$. The cases of $w=y,z$ are analogous.
\begin{lem}\label{lem2.5_1} For a positive integer $p$, ${\bm D}_{p}^{w}$ has following properties
\begin{itemize}
\item[1.] ${\bm D}_{2p}^{w}$ is symmetric and the elements of ${\bm D}_{2p}^{w}$ satisfy
\begin{align*}
(D_{2p}^{w})_{k,k+m}=(D_{2p}^{w})_{k,k-m},\ (D_{2p}^{w})_{k,m+N_w}=(D_{2p}^{w})_{k,m},
\end{align*}
where $k,m=1,\cdots, N_{w}$.
\item[2.]  ${\bm D}_{2p-1}^{w}$ is skew-symmetric and the elements of ${\bm D}_{2p-1}^{w}$ satisfy
\begin{align*}
(D_{2p-1}^{w})_{k,k+N_{w}/2}=0,\ (D_{2p-1}^{w})_{k,k+m}=-(D_{2p-1}^{w})_{k,k-m},\ (D_{2p-1}^{w})_{k,m+N_w}=(D_{2p-1}^{w})_{k,m},
\end{align*}
where $k,m=1,\cdots, N_{w}$.
\end{itemize}
\end{lem}
\begin{prf} \rm Here, we only present the proof for the case of $w=x$. The proofs for the cases of
$w=y,z$ are analogous. By careful calculation, we can obtain
\begin{align}
\frac{d^{p} g_{j}(x)}{dx^{p}}=
 \left \{
 \aligned
 &\frac{2(\text{i}\mu_{x})^{p}}{N_{x}}\sum_{l=1}^{N_{x}/2}\frac{l^{p}}{a_{l}}\cos(l\mu_{x}(x-x_{j})),\ \ \ p\ \text{is an even integer},\\
 &\frac{2\text{i}(\text{i}\mu_{x})^{p}}{N_{x}}\sum_{l=1}^{N_{x}/2}\frac{l^{p}}{a_{l}}\sin(l\mu_{x}(x-x_{j})),\ \ \ p\ \text{is an odd integer},
 \endaligned
 \right.
\end{align}
which implies that
\begin{align}\label{SX_3.5}
({\bm D}_{p}^{x})_{j,k}
 =\left \{
 \aligned
 &\frac{2(\text{i}\mu_{x})^{p}}{N_{x}}\sum_{l=1}^{N_{x}/2}\frac{l^{p}}{a_{l}}\cos\Big((j-k)\frac{2\pi l}{N_x}\Big),\ \ \ p\ \text{is an even integer},\\
 &\frac{2\text{i}(\text{i}\mu_{x})^{p}}{N_{x}}\sum_{l=1}^{N_{x}/2}\frac{l^{p}}{a_{l}}\sin\Big((j-k)\frac{2\pi l}{N_x}\Big),\ \ \ p\ \text{is an odd integer}.
 \endaligned
 \right.
\end{align}
By simple calculation, we finish the proof.
\qed
\end{prf}
\begin{lem}\label{lem2.5_2} For a positive integer $p$,
$({\bm D}_{1}^{w})^{p}$ has the properties
\begin{itemize}
\item[1.] $({\bm D}_{1}^{w})^{2p}$ is symmetric, and the elements of $({\bm D}_{1}^{w})^{2p}$ satisfy
\begin{align*}
&[({\bm D}_{1}^{w})^{2p}]_{k,k+m}=[({\bm D}_{1}^{w})^{2p}]_{k,k-m},\ \ [(D_{1}^{w})^{2p}]_{k,m+N_w}=[(D_{1}^{w})^{2p}]_{k,m},
\end{align*}
where $k,m=1,\cdots, N_{w}$.
\item[2.]  $({\bm D}_{1}^{w})^{2p-1}$ is skew-symmetric and the elements of $({\bm D}_{1}^{w})^{2p-1}$ satisfy
\begin{align*}
&[({\bm D}_{1}^{w})^{2p-1}]_{k,k+N_{w}/2}=0,\ [({\bm D}_{1}^{w})^{2p-1}]_{k,k+m}=-[({\bm D}_{1}^{w})^{2p-1}]_{k,k-m},\\
&[({\bm D}_{1}^{w})^{2p-1}]_{k,m+N_w}=[({\bm D}_{1}^{w})^{2p-1}]_{k,m},
\end{align*}
where $k,m=1,\cdots, N_{w}$.
\end{itemize}
\end{lem}
\begin{prf}
By noting the relationship (\ref{CQ}) and Eq. \eqref{SX_3.5}, we finish the proof.\qed
\end{prf}

\begin{lem}\label{lem3.3}\cite{GCW14} For the matrix ${\bm D}_{p}^{w},\ p=1,2,\cdots$, we have
\begin{align*}
& {\bm D}_{p}^{w}=
 \left \{
 \aligned
 &\mathcal{F}_{N_{w}}^{-1}{\bm \Lambda}^{p}_{w}\mathcal{F}_{N_{w}},\ \ \ p\ \text{is an odd integer,}\\
 &\mathcal{F}_{N_{w}}^{-1}\tilde{\bm \Lambda}^{p}_{w}\mathcal{F}_{N_{w}},\ \ \ p\ \text{is an even integer,}
 \endaligned
 \right.
 \end{align*}
where ${\bm \Lambda}_{w}$ and $\tilde{\bm \Lambda}_{w}$ are the diagonal matrices whose (non-zero) entries are the scaled wave-numbers
\begin{align*}
&{\bm \Lambda}_{w}=\text{\rm i}\mu_{w}\text{\rm diag}\big(0,1,\cdots,\frac{N_{w}}{2}-1,0,-\frac{N_{w}}{2}+1,\cdots,-2,-1\big),\\
&\tilde{\bm \Lambda}_{w}=\text{\rm i}\mu_{w}\text{\rm diag}\big(0,1,\cdots,\frac{N_{w}}{2}-1,\frac{N_{w}}{2},-\frac{N_{w}}{2}+1,\cdots,-2,-1\big),
\end{align*}
and $\mathcal{F}_{N_{w}}$ is the matrix of DFT coefficients with entries given by
$(\mathcal{F}_{N_{w}})_{j,k}=\omega_{N_{w}}^{-j,k},$
$\omega_{N_{w}}=e^{\text{i}\frac{2\pi}{N_{w}}},\ (\mathcal{F}_{N_{w}}^{-1})_{j,k}=\frac{1}{N_{w}}\omega_{N_{w}}^{j,k}$.
\end{lem}


Next, applying the Fourier pseudo-spectral methods to the Hamiltonian formulation (\ref{MEH1}), we can obtain
\begin{align}\label{MESH1}
\frac{d}{d t}\left(\begin{array}{c}
              {\bm H} \\
              {\bm E}  \\
             \end{array}
\right)=\left(\begin{array}{cc}
              0     &  -{\bm I}  \\
              {\bm I}    &  0 \\
             \end{array}
\right)\left(\begin{array}{c}
              \nabla_{{\bm H}}{\mathcal{{\bar H}}_{1}} \\
              \nabla_{{\bm E}}{\mathcal{{\bar H}}_{1}} \\
             \end{array}
\right)=\left(\begin{array}{cc}
              0     &  -{\bm I}  \\
              {\bm I}    &  0 \\
             \end{array}
\right)\left(\begin{array}{c}
              \epsilon^{-1}{\bm D}{\bm H} \\
              \mu^{-1}{\bm D}{\bm E} \\
             \end{array}
\right)=\hat{\mathcal{J}}_{1}\nabla \mathcal{{\bar H}}_{1},
\end{align}
where ${\bm I}\in \mathbb{R}^{3s_1\times 3s_1},\ s_1=N_{x}\times N_{y}\times N_{z}$,
\begin{align}\label{SMEH1}
\mathcal{{\bar H}}_{1}=\frac{1}{2\epsilon}{\bm H}^{T}({\bm D}{\bm H})+\frac{1}{2\mu}{\bm E}^{T}({\bm D}{\bm E}),
\end{align}
  and
\begin{align}
{\bm D}&=\left(\begin{array}{ccc}
              0     &  -{\bm D}_{1}^{z}\otimes {\bm I}_{N_{y}}\otimes {\bm I}_{N_x} & {\bm I}_{N_x}\otimes {\bm D}_{1}^{y}\otimes {\bm I}_{N_x} \\
              {\bm D}_{1}^{z}\otimes {\bm I}_{N_{y}}\otimes {\bm I}_{N_x}    &  0 & -{\bm I}_{N_{z}}\otimes {\bm I}_{N_{y}}\otimes{\bm D}_{1}^{x} \\
              -{\bm I}_{N_x}\otimes {\bm D}_{1}^{y}\otimes {\bm I}_{N_x} &{\bm I}_{N_{z}}\otimes {\bm I}_{N_{y}}\otimes{\bm D}_{1}^{x} & 0\\
             \end{array}
\right)\nonumber\\
&:=\left(\begin{array}{ccc}
              0     &  -{\bm D}_{3} &{\bm D}_{2}  \\
              {\bm D}_{3}    &  0 & -{\bm D}_{1} \\
              -{\bm D}_{2} &{\bm D}_{1} & 0\\
             \end{array}
\right),
\end{align}
is symmetric structure matrix corresponding to the discretization of the operator $\nabla\times$.
\begin{rmk}\label{rmk2.1}
The matrices ${\bm D}_{1}$, ${\bm D}_{2}$, and ${\bm D}_{3}$  admit the following properties
\begin{itemize}
\item[1.] Commutative law of multiplication
\begin{align*}
{\bm D}_{s_{2}}{\bm D}_{s_{3}}={\bm D}_{s_{3}}{\bm D}_{s_{2}},
\end{align*}
where $\ s_2,s_3=1,2,3,$ and $s_2\not=s_3$.
\item[2.] Skew symmetry
\begin{align*}
{\bm D}_p^{T}=-{\bm D}_p,\ p=1,2,3.
\end{align*}
\end{itemize}
 Moreover, one should note that 
the components of ${\bm H}$ and ${\bm E}$ in Eqs. \eqref{MESH1}-\eqref{SMEH1} are the values of grid functions and different from those in (\ref{ME}).
Without being confused, the notations ${\bm H}$ and ${\bm E}$ are still be adopted in the subsequent sections.

\end{rmk}
\subsection{Sixth order AVF method approximation in time}
In this subsection, the sixth order average vector field (AVF) method \cite{LWQ14} will be
employed to discretize finite dimensional Hamiltonian system (\ref{MESH1}) in time.
 For ordinary differential equations (ODEs)
\begin{align}\label{ODE}
\frac{\text{d}{\bm y}}{\text{d}t}=f({\bm y}),\ {\bm y}(0)={\bm y}_{0}\in \mathbb{R}^{2d},
\end{align}
the AVF method is defined by
\begin{align}\label{2AVF}
\frac{{\bm y}^{n+1}-{\bm y}^{n}}{\tau}=\int_{0}^{1} f((1-\xi){\bm y}^{n}+\xi
{\bm y}^{n+1})d\xi,
\end{align}
where $\tau$ is the time step.
The average vector field (AVF) method \eqref{2AVF}, which was first derived in Ref. \cite{MQR99}, and then identified as a B-series
method in Ref. \cite{QM08}, is affine-covariant, of order 2 and self-adjoint \cite{CMOQ10}.
The remarkable advantage of the AVF method (\ref{2AVF}) is that it can preserve exactly the energy integral for any Hamiltonian
system with a constant structure matrix.
With the help of the concrete
formulas of the substitution law \cite{CHV05,FQ10} for the trees of order 5, the second order  AVF method (\ref{2AVF}) was extended to sixth order  for
the Hamiltonian system with a constant structure matrix  \cite{LWQ14}.
The sixth order AVF method can be rewritten as a compact form \cite{LWQ14}
\begin{equation}\label{6AVF}
\frac{{\bm y}^{n+1}-{\bm y}^{n}}{\tau}=~ \tilde{\bm S}(\frac{{\bm y}^{n+1}+{\bm y}^{n}}{2})\int_{0}^{1}\nabla H((1-\xi){\bm y}^{n}+\xi {\bm y}^{n+1})d\xi,
\end{equation}
where
\begin{align*}
\tilde{\bm S}({\bm y})=~ & \Big[{\bm I}_{2d\times 2d}-\frac{\tau^{2}}{12}{\bm S}\mathcal{\bm H}{\bm S}\mathcal{\bm H}+\frac{\tau^{4}}{720}\Big(6S\mathcal{\bm H}{\bm S}\mathcal{\bm H}{\bm S}\mathcal{\bm H}{\bm S}\mathcal{\bm H}-{\bm S}\mathcal{\bm T}{\bm S}\mathcal{\bm T}
+{\bm S}\mathcal{\bm H}{\bm S}\mathcal{\bm H}{\bm S}\mathcal{\bm T}\\
&-{\bm S}\mathcal{\bm T}{\bm S}\mathcal{\bm H}{\bm S}\mathcal{\bm H}
   -\frac{3}{2}{\bm S}\mathcal{\bm H}{\bm S}\mathcal{\bm L} -\frac{3}{2}{\bm S}\mathcal{\bm L}{\bm S}\mathcal{\bm H}+3{\bm S}\mathcal{\bm R}{\bm S}\mathcal{\bm H}+3{\bm S}\mathcal{\bm H}{\bm S}\mathcal{\bm R}\Big)\Big]{\bm S},
\end{align*}
${\bm S}$ is the constant structure matrix, $H({\bm y}): \mathbb{R}^{2d} \rightarrow \mathbb{R}$ is Hamiltonian function, and
 the symmetric matrices $\mathcal{\bm H}({\bm y})$, $\mathcal{\bm T}({\bm y})$, $\mathcal{\bm L}({\bm y})$ and $\mathcal{\bm R}({\bm y})$ are given by the Einstein summation convention
\begin{align*}
& \mathcal{\bm H}_{ij}:=\frac{\partial^{2}H}{\partial y_{i}\partial y_{j}},\quad
\mathcal{\bm T}_{ij}:=\frac{\partial^{3}H}{\partial y_{i}\partial y_{j}\partial y_{k}}S^{kl}\frac{\partial H}{\partial y_{l}},\\
& \mathcal{\bm L}_{ij}:=\frac{\partial^{4}H}{\partial y_{i}\partial y_{j}\partial y_{k}\partial y_{l}}
S^{km}\frac{\partial H}{\partial y_{m}}S^{ln}\frac{\partial H}{\partial y_{n}},
\quad \mathcal{\bm R}_{ij}:=\frac{\partial^{3}H}{\partial y_{i}\partial y_{j}\partial y_{k}}S^{kl}
\frac{\partial^{2} H}{\partial y_{l}\partial y_{m}}S^{mn}\frac{\partial H}{\partial y_{n}}.
\end{align*}
Some others works on the AVF method can be found in Refs. \cite{CGM12,DO11}.

Applying the sixth order AVF method (\ref{6AVF}) to Eq. (\ref{MESH1}), we can obtain
\begin{align}\label{MEHS2}
\frac{1}{\tau}
\Big(
  \begin{array}{c}
    {\bm H}^{n+1}-{\bm H}^{n} \\
    {\bm E}^{n+1}-{\bm E}^{n} \\
  \end{array}
\Big)=&\Big[\Big(
                \begin{array}{cc}
                  {\bm I} & 0 \\
                  0 & {\bm I} \\
                \end{array}
              \Big)-\frac{\tau^2}{12}\Big(
                                         \begin{array}{cc}
                                           -c^2 {\bm D}^2 & 0 \\
                                           0 & -c^2 {\bm D}^2 \\
                                         \end{array}
                                       \Big)
                                       -
  \frac{\tau^4}{120}\Big(
                                         \begin{array}{cc}
                                           -c^4 {\bm D}^4 & 0 \\
                                           0 & -c^4{\bm D}^4 \\
                                         \end{array}
                                       \Big)
  \Big]\nonumber\\
&\Big(
          \begin{array}{cc}
            0 & -{\bm I} \\
            {\bm I} & 0 \\
          \end{array}
        \Big)\Big(
                 \begin{array}{c}
                   \varepsilon^{-1}{\bm D}\int_0^1 ((1-\xi){\bm H}^n+\xi{\bm H}^{n+1})d\xi \\
                   \mu^{-1}{\bm D}\int_0^1 ((1-\xi){\bm E}^n+\xi{\bm E}^{n+1})d\xi \\
                 \end{array}
               \Big),
\end{align}
where $c=1/\sqrt{\epsilon\mu}$. The integration in Eq. (\ref{MEHS2}) can be calculated exactly to give

\begin{align}\label{SO-ECS1}
&\frac{1}{\tau}\Big(
  \begin{array}{c}
    {\bm H}^{n+1}-{\bm H}^{n} \\
    {\bm E}^{n+1}-{\bm E}^{n} \\
  \end{array}
\Big)=\Big[\Big(
                \begin{array}{cc}
                  0 & -{\bm D} \\
                  {\bm D}   &0  \\
                \end{array}
              \Big)+\frac{c^2\tau^2}{12}\Big(
                                         \begin{array}{cc}
                                          0& -{\bm D}^3 \\
                                           {\bm D}^3 & 0 \\
                                         \end{array}
                                       \Big)
                                       +\frac{c^4\tau^4}{120}\Big(
                                         \begin{array}{cc}
                                           0  & -{\bm D}^5 \\
                                           {\bm D}^5 & 0 \\
                                         \end{array}
                                       \Big)
  \Big]\nonumber\\
 &~~~~~~~~~~~~~~~~~~~~~~~~~~~~~~ \Big(
                 \begin{array}{c}
                   \varepsilon^{-1}({\bm H}^{n+1}+{\bm H}^{n})/2 \\
                   \mu^{-1}({\bm E}^{n+1}+{\bm E}^{n})/2 \\
                 \end{array}
               \Big),
\end{align}
which comprises our sixth order energy-conserved scheme for the 3D Maxwell's equations. In fact,
Eqs. (\ref{SO-ECS1}) can be expressed as an equivalent form
\begin{align}\label{SO-ECS2}
\Big(\begin{array}{cc}
              \frac{2\mu }{\tau}{\bm I}     &{\hat {\bm D}}  \\
              -{\hat {\bm D}}    &  \frac{2\epsilon }{\tau}{\bm I} \\
             \end{array}
\Big)\Big(\begin{array}{c}
              {\bm H}^{n+1} \\
              {\bm E}^{n+1}  \\
             \end{array}
\Big)=\Big(\begin{array}{cc}
              \frac{2\mu }{\tau}{\bm I}& -{\hat {\bm D}}  \\
              {\hat {\bm D}}    &  \frac{2\epsilon }{\tau}{\bm I} \\
             \end{array}
\Big)\Big(\begin{array}{c}
  {\bm H}^{n} \\
  {\bm E}^{n} \\
\end{array}
\Big),
\end{align}
where ${\hat {\bm D}}={\bm D}+\frac{c^{2}\tau^{2}}{12}{\bm D}^{3}+\frac{c^{4}\tau^{4}}{120}{\bm D}^{5}$.
\begin{rmk} If exchanging
$n+1\longleftrightarrow n$ and $\tau\longleftrightarrow -\tau$, then we find that the scheme (\ref{SO-ECS1}) unaltered.
Thus, according to Ref. \cite{ELW06}, the scheme (\ref{SO-ECS1}) is symmetric.
Moreover, the proposed scheme is uniquely solvable (see Appendix).
\end{rmk}

\section{Discrete conservation laws and divergence preservation}
In this section, we will rigorously prove the scheme (\ref{SO-ECS1}) satisfies the discrete version of the five energy conservation laws, the three momentum conservation laws, the symplecticity as well as the two discrete divergence-free fields. But before, we introduce a new semi-norm for any grid functions $U_{j,k,m},(x_j,y_k,z_m)\in\Omega_h$,
\begin{align*}
||{\bm U}||_{{\bm D}_{1}^{x},h}=\langle{\bm D}_{1}{\bm U},{\bm D}_{1}{\bm U}\rangle_{h}^{\frac{1}{2}},\ ||{\bm U}||_{{\bm D}_{1}^{y},h}
=\langle{\bm D}_{2}{\bm U},{\bm D}_{2}{\bm U}\rangle_{h}^{\frac{1}{2}},\
||{\bm U}||_{{\bm D}_{1}^{z},h}=\langle{\bm D}_{3}{\bm U},{\bm D}_{3}{\bm U}\rangle_{h}^{\frac{1}{2}},
\end{align*}
 and
\begin{align*}
||{\bm U}||_{{\bm D}_{1}^{w},h}^{2}=||{\bm U}_{x}||_{{\bm D}_{1}^{w},h}^{2}+||{\bm U}_{y}||_{{\bm D}_{1}^{w},h}^{2}+||{\bm U}_{z}||_{{\bm D}_{1}^{w},h}^{2},
\end{align*}
for a vector ${\bm U}=[({\bm U}_{x})^{T},({\bm U}_{y})^{T},({\bm U}_{z})^{T}]^{T}$.
\begin{thm} The solutions ${\bm H}^{n}$ and ${\bm E}^{n}$ of the proposed
scheme (\ref{SO-ECS1}) satisfy the discrete sympletic conservation law
\begin{align}\label{SCL1}
\omega^{n+1}=\omega^{n},\ \omega^{n}=d{\bm E}^{n}\wedge d{\bm H}^{n},\ n=0,\cdots,M-1,
\end{align}
where
\begin{align*}
&d{\bm E}^{n}\wedge d{\bm H}^{n}=d{\bm E}_{x}^{n}\wedge d{\bm H}_{x}^{n}+d{\bm E}_{y}^{n}\wedge d{\bm H}_{y}^{n}+d{\bm E}_{z}^{n}\wedge d{\bm H}_{z}^{n},\\
&d{\bm E}_{w}^{n}\wedge d{\bm H}_{w}^{n}=\sum_{j,k,m}d{\bm E}_{w_{j,k,m}}^{n}\wedge d{\bm H}_{w_{j,k,m}}^{n}.
\end{align*}
\end{thm}
\begin{prf}\rm Eq. (\ref{SO-ECS1}) can be rewritten as
\begin{align}\label{SCL2}
\left(\begin{array}{cc}
              0     &I  \\
              -I    &  0 \\
             \end{array}
\right)\left(\begin{array}{c}
              \frac{{\bm H}^{n+1}-{\bm H}^{n}}{\tau} \\
              \frac{{\bm E}^{n+1}-{\bm E}^{n}}{\tau}   \\
             \end{array}
\right)=\left(\begin{array}{cc}
             \epsilon^{-1}{\hat {\bm D}}     & 0 \\
              0    &  \mu^{-1}{\hat {\bm D}} \\
             \end{array}
\right)\left(\begin{array}{c}
  \frac{{\bm H}^{n+1}+{\bm H}^{n}}{2} \\
  \frac{{\bm E}^{n+1}+{\bm E}^{n}}{2} \\
\end{array}
\right).
\end{align}
The variational equations associated with Eq. (\ref{SCL2}) read
\begin{align}\label{SX_3.3}
\left(\begin{array}{cc}
              0     &I  \\
              -I    &  0 \\
             \end{array}
\right)\left(\begin{array}{c}
              \frac{d{\bm H}^{n+1}-d{\bm H}^{n}}{\tau} \\
              \frac{d{\bm E}^{n+1}-d{\bm E}^{n}}{\tau}   \\
             \end{array}
\right)=\left(\begin{array}{cc}
             \epsilon^{-1}{\hat {\bm D}}     & 0 \\
              0    &  \mu^{-1}{\hat {\bm D}} \\
             \end{array}
\right)\left(\begin{array}{c}
  \frac{d{\bm H}^{n+1}+d{\bm H}^{n}}{2} \\
  \frac{d{\bm E}^{n+1}+d{\bm E}^{n}}{2} \\
\end{array}
\right).
\end{align}
Taking the wedge product with $\left(\begin{array}{c}
  \frac{d{\bm H}^{n+1}+d{\bm H}^{n}}{2} \\
  \frac{d{\bm E}^{n+1}+d{\bm E}^{n}}{2} \\
\end{array}
\right)$ on both sides of Eq. \eqref{SX_3.3} and noting the fact
\begin{align}
\left(\begin{array}{c}
  \frac{d{\bm H}^{n+1}+d{\bm H}^{n}}{2} \\
  \frac{d{\bm E}^{n+1}+d{\bm E}^{n}}{2} \\
\end{array}
\right)\wedge\left(\begin{array}{cc}
             \epsilon^{-1}{\hat {\bm D}}     & 0 \\
              0    &  \mu^{-1}{\hat {\bm D}} \\
             \end{array}
\right)\left(\begin{array}{c}
  \frac{d{\bm H}^{n+1}+d{\bm H}^{n}}{2} \\
  \frac{d{\bm E}^{n+1}+d{\bm E}^{n}}{2} \\
\end{array}
\right)=0,
\end{align}
we obtain the discrete sympletic conservation law (\ref{SCL1}).
\qed
\end{prf}
\begin{rmk}
Since the system of the Maxwell's equations (\ref{ME}) is linear, we can get the  sympletic conservation law for the method (\ref{SO-ECS1}).
Usually, energy-conserved methods do not imply the preservation of the symplectic conservation law
for the given system.
\end{rmk}

Since the AVF method can preserve Hamiltonian energy automatically, we can derive the following energy conservation laws, immediately.
\begin{thm}  The solutions ${\bm H}^{n}$ and ${\bm E}^{n}$ of the proposed
scheme (\ref{SO-ECS1}) satisfy the discrete energy conservation laws
\begin{align}\label{ECL-I1}
\mathcal{E}_{1}^{n+1}=\mathcal{E}_{1}^{n},\ \mathcal{E}_{1}^{n}=\frac{1}{2\epsilon}\langle{\bm H}^{n},{\bm D}{\bm H}^{n}\rangle_{h}+\frac{1}{2\mu}\langle{\bm E}^{n},{\bm D}{\bm E}^{n}\rangle_{h},\ n=0,\cdots,M-1,
\end{align}
and
\begin{align}\label{ECL-II1}
\mathcal{E}_{2}^{n+1}=\mathcal{E}_{2}^{n},\ \mathcal{E}_{2}^{n}=\frac{\mu}{2}||{\bm H}^{n}||_{h}^{2}+\frac{\epsilon}{2}||{\bm E}^{n}||_{h}^{2},\ n=0,\cdots,M-1.
\end{align}
\end{thm}

In addition, the scheme (\ref{SO-ECS1}) also satisfies the following energy conservation laws.
\begin{thm} The solutions ${\bm H}^{n}$ and ${\bm E}^{n}$ of the
scheme (\ref{SO-ECS1}) satisfy the discrete energy conservation law
\begin{align}\label{ECL-III1}
\mathcal{E}_{3}^{n+1}=\mathcal{E}_{3}^{n},\ \mathcal{E}_{3}^{n}=\frac{\mu}{2}||\hat{\delta}_{t}{\bm H}^{n-1/2}||_{h}^{2}
+\frac{\epsilon}{2}||\hat{\delta}_{t}{\bm E}^{n-1/2}||_{h}^{2},\ n=1,\cdots,M-1.
\end{align}
\end{thm}
\begin{prf}\rm  We can deduce from (\ref{SO-ECS1})
\begin{align}\label{ECL-III2}
\frac{1}{\tau}\Bigg(\begin{array}{c}
              \mu{\hat {\delta}}_{t}{\bm H}^{n+1/2}-\mu{\hat {\delta}}_{t}{\bm H}^{n-1/2}\\
              \epsilon{\hat {\delta}}_{t}{\bm E}^{n+1/2}-\epsilon{\hat {\delta}}_{t}{\bm E}^{n-1/2} \\
             \end{array}
\Bigg)=\Bigg(\begin{array}{cc}
              0     &-{\hat {\bm D}}  \\
              {\hat {\bm D}}    &  0 \\
             \end{array}
\Bigg)\Bigg(\begin{array}{c}
  \frac{1}{2}({\hat {\delta}}_{t}{\bm H}^{n+1/2}+{\hat {\delta}}_{t}{\bm H}^{n-1/2}) \\
  \frac{1}{2} ({\hat {\delta}}_{t}{\bm E}^{n+1/2}+{\hat {\delta}}_{t}{\bm E}^{n-1/2}) \\
\end{array}
\Bigg).
\end{align}
Making the inner product with $\frac{1}{2}\Big(\begin{array}{c}
              {\hat {\delta}}_{t}{\bm H}^{n+1/2}+{\hat {\delta}}_{t}{\bm H}^{n-1/2}\\
              {\hat {\delta}}_{t}{\bm E}^{n+1/2}+{\hat {\delta}}_{t}{\bm E}^{n-1/2} \\
             \end{array}
\Big)$ on both sides of Eq. (\ref{ECL-III2}) and by virtue of
 the skew-symmetric property of the matrix $\Big(\begin{array}{cc}
              0     &-{\hat {\bm D}}  \\
              {\hat {\bm D}}    &  0 \\
             \end{array}
\Big)$, we finish the proof.
\qed
\end{prf}
\begin{thm}\label{SX_thm3.4} The solutions ${\bm H}^{n}$ and ${\bm E}^{n}$ of the proposed
scheme (\ref{SO-ECS1}) also possess the discrete energy conservation laws
\begin{align}\label{ECL-IV}
\mathcal{E}_{4_{w}}^{n+1}=\mathcal{E}_{4_{w}}^{n},
\ \mathcal{E}_{4_{w}}^{n}=\mu||{\bm H}^{n}||_{{\bm D}_{1}^{w},h}^{2}+\epsilon||{\bm E}^{n}||_{{\bm D}_{1}^{w},h}^{2},\ n=0,\cdots,M-1,
\end{align}
and
\begin{align}\label{ECL-V}
\mathcal{E}_{5_{w}}^{n+1}=\mathcal{E}_{5_{w}}^{n},
\ \mathcal{E}_{5_{w}}^{n}=\mu||{\hat \delta_{t}}{\bm H}^{n-1/2}||_{{\bm D}_{1}^{w},h}^{2}
+\epsilon||{\hat \delta_{t}}{\bm E}^{n-1/2}||_{{\bm D}_{1}^{w},h}^{2},\ n=1,\cdots,M-1.
\end{align}
\end{thm}
\begin{prf} \rm We define block diagonal matrices
${\bm M}_{x}=\Bigg(\begin{array}{ccc}
              {\bm D}_{1}     &\ &\   \\
              \    & {\bm D}_{1}&\  \\
                \    & \ & {\bm D}_{1} \\
             \end{array}
\Bigg),$\
${\bm M}_{y}=\Bigg(\begin{array}{ccc}
              {\bm D}_{2}     &\ &\   \\
              \    & {\bm D}_{2}&\  \\
                \    & \ & {\bm D}_{2} \\
             \end{array}
\Bigg),$
and\
${\bm M}_{z}=\Bigg(\begin{array}{ccc}
              {\bm D}_{3}     &\ &\   \\
              \    & {\bm D}_{3}&\  \\
                \    & \ & {\bm D}_{3} \\
             \end{array}
\Bigg).$ According to Remark \ref{rmk2.1}, it is easy to see that 
${\bm M}_{w}{\hat {\bm D}}={\hat {\bm D}}{\bm M}_{w}$. Then,
left-multiplying (\ref{SO-ECS1}) with block diagonal matrix
${\bm G}_{w}=\Big(\begin{array}{cc}
              {\bm M}_{w}    &\   \\
              \    & {\bm M}_{w} \\
             \end{array}
\Big),$\
we have
\begin{align}\label{ECL-IV1}
\frac{1}{\tau}\left(\begin{array}{c}
              \mu ({\bm M}_{w}{\bm H}^{n+1})-\mu({\bm M}_{w}{\bm H}^{n})\\
              \epsilon ({\bm M}_{w}{\bm E}^{n+1})-\epsilon ({\bm M}_{w}{\bm E}^{n}) \\
             \end{array}
\right)=\left(\begin{array}{cc}
              0     &-{\hat {\bm D}}  \\
              {\hat {\bm D}}    &  0 \\
             \end{array}
\right)\left(\begin{array}{c}
  \frac{{\bm M}_{w}{\bm H}^{n+1}+{\bm M}_{w}{\bm H}^{n}}{2} \\
  \frac{{\bm M}_{w}{\bm E}^{n+1}+{\bm M}_{w}{\bm E}^{n}}{2} \\
\end{array}
\right).
\end{align}
Computing the inner product with $\left(\begin{array}{c}
              {\bm M}_{w}{\bm H}^{n+1}+{\bm M}_{w}{\bm H}^{n}\\
              {\bm M}_{w}{\bm E}^{n+1}+{\bm M}_{w}{\bm E}^{n} \\
             \end{array}
\right)$ on both sides of Eq. (\ref{ECL-IV1}) and we can derive that
\begin{align}\label{ECL-IV2}
\mu ||{\bm M}_{w}{\bm H}^{n+1}||^{2}+\epsilon ||{\bm M}_{w}{\bm E}^{n+1}||_{h}^{2}
=\mu ||{\bm M}_{w}{\bm H}^{n}||^{2}+\epsilon ||{\bm M}_{w}{\bm E}^{n}||_{h}^{2}.
\end{align}
Note the notation introduced for grid vector function and we can get the discrete energy conservation law (\ref{ECL-IV}).
\qed
\end{prf}

Left-multiplying (\ref{ECL-III2}) with the block diagonal matrix ${\bm G}_{w}$, and taking the inner product with
\begin{align*}
\Big(\begin{array}{c}
              {\hat {\delta}}_{t}{\bm M}_{w}{\bm H}^{n+1/2}+{\hat {\delta}}_{t}{\bm M}_{w}{\bm H}^{n-1/2}\\
              {\hat {\delta}}_{t}{\bm M}_{w}{\bm E}^{n+1/2}+{\hat {\delta}}_{t}{\bm M}_{w}{\bm E}^{n-1/2} \\
             \end{array}
\Big), \
\end{align*}
we can obtain the second discrete energy conservation law of Theorem \ref{SX_thm3.4}.

Next, we can show that the resulting numerical scheme (\ref{SO-ECS1}) preserves
the corresponding discrete momentum conservation laws.
\begin{thm} The solutions ${\bm H}^{n}$ and ${\bm E}^{n}$ of the proposed
scheme (\ref{SO-ECS1}) capture momentum conservation laws
\begin{align}\label{MCL-I}
\mathcal{M}_{w}^{n+1}=\mathcal{M}_{w}^{n},\ \mathcal{M}_{w}^{n}=\langle{\bm H}^{n},{\bm M}_{w}{\bm E}^{n}\rangle_{h},\ n=0,\cdots,M-1,
\end{align}
where ${\bm M}_{w}$ is defined as above.
\end{thm}
\begin{prf}\rm  Eq. (\ref{ECL-IV1}) can be rewritten as
\begin{align}\label{MCL-I1}
\frac{1}{\tau}\Big(\begin{array}{c}
               {\bm M}_{w}{\bm H}^{n+1}-{\bm M}_{w}{\bm H}^{n}\\
              -({\bm M}_{w}{\bm E}^{n+1}-{\bm M}_{w}{\bm E}^{n}) \\
             \end{array}
\Big)=\Big(\begin{array}{c}
  \mu^{-1}{\hat {\bm D}}{\bm M}_{w}{\bm E}^{n+1/2}\\
  \epsilon^{-1}{\hat {\bm D}}{\bm M}_{w}{\bm H}^{n+1/2}\\
\end{array}
\Big).
\end{align}
Then taking the inner product with $\Big(\begin{array}{c}
              {\bm E}^{n+1/2}\\
             {\bm H}^{n+1/2}\\
             \end{array}
\Big)$ on both sides of Eq. (\ref{MCL-I1}) and we have
\begin{align}
\frac{1}{\tau}\Big(\langle{\bm H}^{n+1},{\bm M}_{w}{\bm E}^{n+1})-({\bm H}^{n},{\bm M}_{w}{\bm E}^{n}\rangle_{h}\Big)
&=-\mu^{-1}\langle{\bm E}^{n+1/2},{\hat {\bm D}}{\bm M}_{w}{\bm E}^{n+1/2}\rangle_{h}\nonumber\\
&-\epsilon^{-1}\langle{\bm H}^{n+1/2},{\hat {\bm D}}{\bm M}_{w}{\bm H}^{n+1/2}\rangle_{h}.
\end{align}
By noting $({\hat {\bm D}}{\bm M}_{w})^{T}=-{\hat {\bm D}}{\bm M}_{w}$, we can obtain (\ref{MCL-I}).
\qed
\end{prf}

Finally,
we show that the proposed scheme (\ref{SO-ECS1}) can preserve the two discrete divergence-free fields of the Maxwell's equations.
\begin{thm}  For the scheme (\ref{SO-ECS1}), the following discrete divergence-free fields hold
\begin{align}\label{DP}
{\bar \nabla}\cdot(\epsilon {\bm E}^{n+1})={\bar \nabla}\cdot(\epsilon {\bm E}^{n})
\  \text{and}\  {\bar \nabla}\cdot(\mu {\bm H}^{n+1})={\bar \nabla}\cdot(\mu {\bm H}^{n}),\ n=0,\cdots,M-1,
\end{align}
where ${\bar \nabla}=\left(\begin{array}{c}
              {\bm {\bm D}}_{1}\\
             {\bm {\bm D}}_{2}\\
             {\bm {\bm D}}_{3}\\
             \end{array}
\right)$, corresponding to the discretization of the operator $\nabla$ and
\begin{align*}
&{\bar \nabla}\cdot(\mu {\bm H})={\bm D}_{1}(\mu{\bm H}_{x})+{\bm D}_{2}(\mu{\bm H}_{y})
+{\bm D}_{3}(\mu{\bm H}_{z}),\\
&{\bar \nabla}\cdot(\epsilon {\bm E})={\bm D}_{1}(\epsilon{\bm E}_{x})
+{\bm D}_{2}(\epsilon{\bm E}_{y})+{\bm D}_{3}(\epsilon{\bm E}_{z}).
\end{align*}
\end{thm}
\begin{prf} \rm Eq. (\ref{SO-ECS1}) is equivalent to
\begin{align}\label{DP1}
&\frac{\mu{\bm H}^{n+1}-\mu{\bm H}^{n}}{\tau}=-{\hat {\bm D}}\frac{{\bm E}^{n+1}+{\bm E}^{n}}{2},\\\label{DP2}
&\frac{\epsilon{\bm E}^{n+1}-\epsilon{\bm E}^{n}}{\tau}={\hat {\bm D}}\frac{{\bm H}^{n+1}+{\bm H}^{n}}{2}.
\end{align}
Left-multiplying (\ref{DP1}) and (\ref{DP2}) with ${\bar \nabla}\cdot$, we then can have
\begin{align}
&\frac{{\bar \nabla}\cdot(\mu{\bm H}^{n+1})-{\bar \nabla}\cdot(\mu{\bm H}^{n})}{\tau}=
-{\bar \nabla}\cdot\left({\hat {\bm D}}\frac{{\bm E}^{n+1}+{\bm E}^{n}}{2}\right),\\
&\frac{{\bar \nabla}\cdot(\epsilon{\bm E}^{n+1})-{\bar \nabla}\cdot(\epsilon{\bm E}^{n})}{\tau}
={\bar \nabla}\cdot\left({\hat {\bm D}}\frac{{\bm H}^{n+1}+{\bm H}^{n}}{2}\right),
\end{align}
whilst the right terms are equal to zero by careful calculation. This completes the proof.
\qed
\end{prf}
\begin{rmk} Eqs. \eqref{DP} are not, strictly speaking, divergence-free, as ${\bar \nabla}\cdot(\epsilon {\bm E}^{0})\not={\bm 0}$ and ${\bar \nabla}\cdot(\epsilon {\bm H}^{0})\not={\bm 0}$. However, ${\bar \nabla}$ is a higher approximation to $\nabla$ in space, so, by choosing appropriate collocation points, $||{\bar \nabla}\cdot(\epsilon {\bm E}^{0})||_{h,\infty}$ and $||{\bar \nabla}\cdot(\epsilon {\bm H}^{0})||_{h,\infty}$ can reach machine precision in practical computation (see Ref. \cite{CWG16} and Section 7). Thus, we also refer to Eqs. \eqref{DP} as divergence free.
\end{rmk}
\section{Convergence analysis}
In this section, we will establish error estimate of the proposed scheme (\ref{SO-ECS1}).
For simplicity,
we let $\Omega=[0,2\pi]^{3}$, $L^{2}(\Omega)$ with the inner product $\langle \cdot , \cdot \rangle$ and
the norm $\Arrowvert\cdot\Arrowvert$ defined previously.
For any positive integer $r$,
the semi-norm and the norm of $H^{r}(\Omega)$ are denoted by $\arrowvert\cdot\arrowvert_{r}$
and $\Arrowvert\cdot\Arrowvert_{r}$, respectively.
$\Arrowvert\cdot\Arrowvert_{0}$ is denoted by $\Arrowvert\cdot\Arrowvert$ for simplicity.
Let $C_{p}^{\infty}$ be the set of infinitely differentiable functions with period $2\pi$ defined
on $\Omega$ for all variables. $H_{p}^{r}(\Omega)$ is the closure of $C_{p}^{\infty}$ in $H^{r}(\Omega)$.

let $N_x=N_y=N_z=N$, the interpolation space $S_{N}^{'''}$ can be written as
\begin{align*}
&S_{N}^{'''}=\Big\{ u| u=\sum_{\arrowvert j\arrowvert,\arrowvert k\arrowvert,\arrowvert m\arrowvert\leq\frac{N}{2}}
\frac{{\hat u}_{j,k,m}}{c_{j}c_{k}c_{m}}e^{\text{i}(jx+ky+mz)}: {\hat u}_{\frac{N}{2},k,m}={\hat u}_{-\frac{N}{2},k,m},\\
&~~~~~~~~~~~~~~~~~~~{\hat u}_{j,\frac{N}{2},m}={\hat u}_{j,-\frac{N}{2},m},\
{\hat u}_{j,k,\frac{N}{2}}={\hat u}_{j,k,-\frac{N}{2}}\Big\},
\end{align*}
where $c_{l}=1,\ |l|<\frac{N}{2},\ c_{-\frac{N}{2}}=c_{\frac{N}{2}}=2$.

We denote
\begin{align}
S_{N}=\Big\{u| u=\sum_{\arrowvert j\arrowvert,\arrowvert k\arrowvert,\arrowvert m\arrowvert\leq\frac{N}{2}}
{\hat u}_{j,k,m}e^{\text{i}(jx+ky+mz)}\Big\}.
\end{align}
It is remarked that $S_{N}^{'''}\subseteq S_{N}$ and $S_{N-2}\subseteq S^{'''}_{N}$.
We denote $P_{N}:[L^{2}(\Omega)]^{3}\to [S_{N}]^{3}$ as the orthogonal projection operator
and recall the interpolation operator $I_{N}:[C(\Omega)]^{3}\to [S_{N}^{'''}]^{3}$.
By noting $P_{N}\partial_{w}u=\partial_{w} P_{N}u,\ w=x,y,z,$
we can see that $\nabla \times$ and $P_{N}$ satisfy the commutative law.
\begin{lem}\label{lem4.1}\cite{CWG16} For ${\bm u}\in [S_{N}^{'''}]^{3}$,
$\Arrowvert {\bm u}\Arrowvert\leqslant \Arrowvert {\bm u}\Arrowvert_{h}\leqslant 2\sqrt{2}\Arrowvert {\bm u}\Arrowvert$.
\end{lem}
\begin{lem}\label{lem4.2} \cite{CQ82} If $0\leqslant \alpha \leqslant r$ and ${\bm u}\in [H_{p}^{r}(\Omega)]^{3}$, then
\begin{align}
\Arrowvert P_{N}{\bm u}-{\bm u}\Arrowvert_{\alpha}\leqslant CN^{\alpha-r}\arrowvert {\bm u}\arrowvert_{r},
\end{align}
 and in addition if $r>\frac{3}{2}$ then
\begin{align}
\Arrowvert I_{N}{\bm u}-{\bm u}\Arrowvert_{\alpha}\leqslant CN^{\alpha-r}\arrowvert {\bm u}\arrowvert_{r}.
\end{align}
\end{lem}
\begin{lem}\label{lem4.3} If ${\bm u}^{*}=P_{N-2}{\bm u},\ {\bm u}\in[H_{p}^{r}(\Omega)]^{3},\ N>2,\ r>\frac{3}{2}$,
 then $\Arrowvert {\bm u}^{*}-{\bm u}\Arrowvert_{N}\leqslant CN^{-r} \arrowvert {\bm u}\arrowvert_{r}$.
 \end{lem}
\begin{prf}\rm By virtue of Lemmas \ref{lem4.1} and \ref{lem4.2}, we have
\begin{align*}
&\Arrowvert {{\bm u}-{\bm u}^{*}}\Arrowvert_{h}=\Arrowvert I_{N}({\bm u}-{\bm u}^{*})\Arrowvert_{h}\\
&~~~~~~~~~~~~~~\leqslant2\sqrt{2}\Arrowvert {I_{N}({\bm u}-{\bm u}^{*})}\Arrowvert
=2\sqrt{2}\Arrowvert {I_{N}{\bm u}-{\bm u}^{*}}\Arrowvert\\
&~~~~~~~~~~~~~~\leqslant2\sqrt{2}(\Arrowvert I_{N}{\bm u}-{\bm u}\Arrowvert+\Arrowvert {\bm u}-{\bm u}^{*}\Arrowvert
)\leqslant CN^{-r}\arrowvert {\bm u}
\arrowvert_{r}.
\end{align*}
\end{prf}
\begin{thm}\label{SX_thm4.1} Suppose that the exact periodic solution components ${\bm H}$ and ${\bm E}$ are
smooth enough: ${\bm H},{\bm E}\in C^{7}\left(0,T;[H^{r}_{p}(\Omega)]^{3}\right),r>\frac{10}{2}$.
The initial conditions are ${\bm H}_{0},\ {\bm E}_{0}\in\left([H^{r}_{p}(\Omega)]^{3}\right),r>\frac{10}{2}$.
For $n\geqslant 0$, let ${\bm H}^{n}$ and ${\bm E}^{n}$ be the solutions of the scheme (\ref{SO-ECS1}).
Then, for any fixed $T>0$, there exists a positive constant $C$ independent of
$\tau,h_{x},h_{y},$ and $h_{z}$, which may be vary in different cases, such that
\begin{align}\label{CAT}
\max\limits_{0 \leq n\leq
M}\left(\epsilon \Arrowvert {\bm H}(t_{n})-{\bm H}^{n}\Arrowvert_{h}^{2}+\mu\Arrowvert{\bm E}(t_{n})
-{\bm E}^{n}\Arrowvert_{h}^{2}\right)^{\frac{1}{2}}\leqslant CT(\tau^{6}+N^{-r}).
\end{align}
\end{thm}
\begin{prf}\rm
Let
\begin{align*}
{\bm E}^{*}=P_{N-2} {\bm E},\ {\bm H}^{*}=P_{N-2} {\bm H}.
\end{align*}
The projection of Eqs. (\ref{ME}) are written as
\begin{align}\label{PEs}
\mu\partial_{t}{\bm H}^{*}+\nabla\times {\bm E}^{*}=0,\
\epsilon\partial_{t}{\bm E}^{*}-\nabla\times {\bm H}^{*}=0.
\end{align}
Let
\begin{align}\label{SX_4.6}
&{\bm\Phi}^{n+\frac{1}{2}}=
\mu\delta_{t}({\bm H}^{*})^{n}+\nabla\times ({\bm E}^{*})^{n+\frac{1}{2}}
+\frac{c^{2}\tau^{2}}{12}(\nabla\times)^{3}({\bm E}^{*})^{n+\frac{1}{2}}
+\frac{c^{4}\tau^{4}}{120}(\nabla\times)^{5}({\bm E}^{*})^{n+\frac{1}{2}},\\\label{SX_4.7}
&{\bm\Psi}^{n+\frac{1}{2}}=\epsilon\delta_{t}({\bm E}^{*})^{n}-\nabla\times ({\bm H}^{*})^{n+\frac{1}{2}}
-\frac{c^{2}\tau^{2}}{12}(\nabla\times)^{3}({\bm H}^{*})^{n+\frac{1}{2}}-\frac{c^{4}\tau^{4}}{120}(\nabla\times)^{5}({\bm H}^{*})^{n+\frac{1}{2}},
\end{align}
where $n=0,\cdots,M-1$.
Eq. \eqref{SX_4.6} can rewritten as
\begin{align}\label{PEQ}
&{\bm\Phi}^{n+\frac{1}{2}}={\bm\Phi}_{1}^{n+\frac{1}{2}}+{\bm\Phi}_{2}^{n+\frac{1}{2}},
\end{align}
where
\begin{align*}
&{\bm\Phi}_{1}^{n+\frac{1}{2}}=
\mu\delta_{t}({\bm H}^{*}-{\bm H})^{n}+\nabla\times ({\bm E}^{*}-{\bm E})^{n+\frac{1}{2}}
+\frac{c^{2}\tau^{2}}{12}(\nabla\times)^{3}({\bm E}^{*}-{\bm E})^{n+\frac{1}{2}}\\
&~~~~~~~~~+\frac{c^{4}\tau^{4}}{120}(\nabla\times)^{5}({\bm E}^{*}-{\bm E})^{n+\frac{1}{2}},
\end{align*}
and
\begin{align*}
&{\bm\Phi}_{2}^{n+\frac{1}{2}}=\mu\delta_{t}{\bm H}^{n}+\nabla\times {\bm E}^{n+\frac{1}{2}}
+\frac{c^{2}\tau^{2}}{12}(\nabla\times)^{3}{\bm E}^{n+\frac{1}{2}}
+\frac{c^{4}\tau^{4}}{120}(\nabla\times)^{5}{\bm E}^{n+\frac{1}{2}},
\end{align*}
By using lemma \ref{lem4.2} and Taylor's expansion at the node $t_{n+\frac{1}{2}}$, we can obtian
\begin{align}\label{CAT1}
&\Arrowvert {\bm\Phi}_{1}^{n+\frac{1}{2}}\Arrowvert\leqslant C(\tau^{6}+N^{-r}),\  \Arrowvert {\bm\Phi}_{2}^{n+\frac{1}{2}}\Arrowvert
\leqslant C\tau^{6},\ n=0,\cdots,M-1.
\end{align}
Thus, we have
\begin{align}\label{CAT1}
&\Arrowvert {\bm\Phi}^{n+\frac{1}{2}}\Arrowvert\leqslant C(\tau^{6}+N^{-r}),\ n=0,\cdots,M-1.
\end{align}
By noting ${\bm\Phi}^{n+\frac{1}{2}}\in[S_{N}^{'''}]^{3}$, then, we can deduce from Lemma \ref{lem4.1} that
\begin{align}\label{SX_4.12}
&\Arrowvert {\bm\Phi}^{n+\frac{1}{2}}\Arrowvert_{h}\leqslant C(\tau^{6}+N^{-r}),\ n=0,\cdots,M-1.
\end{align}
By the analogous argument to ${\bm\Phi}^{n+\frac{1}{2}}$, we have
\begin{align}\label{SX_4.13}
&\Arrowvert {\bm\Psi}^{n+\frac{1}{2}}\Arrowvert_{h}\leqslant C(\tau^{6}+N^{-r}),\ n=0,\cdots,M-1.
\end{align}
Thus, from Eqs. \eqref{SX_4.12} and \eqref{SX_4.13}, we have
\begin{align}
&(\Arrowvert {\bm\Phi}^{n+\frac{1}{2}}\Arrowvert_{h}^{2}+\Arrowvert {\bm\Psi}^{n+\frac{1}{2}}\Arrowvert_{h}^{2})^{\frac{1}{2}}
\leqslant C(\tau^{6}+N^{-r}),\ n=0,\cdots,M-1,
\end{align}
where the inequation $\sqrt{a^{2}+b^{2}}\leqslant a+b,\ {\forall} a,b\geqslant 0$ is used.

Note that if ${\bm E}^{*},\ {\bm H}^{*}\in S_{N}^{'''}$, we can see that
\begin{align*}
&\partial_{x}\widetilde {{ U}}^{*}(x_{j},y_{k},z_{m})=\partial_{x}I_{N}\widetilde {{ U}}^{*}(x_{j},y_{k},z_{m})=[{\bm D}_{1}\widetilde {{\bm U}}^{*}]_{{N_{x}N_{y}(m-1)+N_{x}(k-1)+j}},\\
&\partial_{y}\widetilde {{ U}}^{*}(x_{j},y_{k},z_{m})=\partial_{y}I_{N}\widetilde {{ U}}^{*}(x_{j},y_{k},z_{m})=[{\bm D}_{2}\widetilde {{\bm U}}^{*}]_{{N_{x}N_{y}(m-1)+N_{x}(k-1)+j}},\\
&\partial_{z}\widetilde {{ U}}^{*}(x_{j},y_{k},z_{m})=\partial_{z}I_{N}\widetilde {{ U}}^{*}(x_{j},y_{k},z_{m})=[{\bm D}_{3}\widetilde {{\bm U}}^{*}]_{{N_{x}N_{y}(m-1)+N_{x}(k-1)+j}},
\end{align*}
where \begin{align*}
&\widetilde{{\bm U}}^{*}=(U_{1,1,1}^{*},\cdots,U_{N_{x},1,1}^{*},U_{1,2,1}^{*},\cdots, U_{N_{x},2,1}^{*},\cdots,U_{1,N_y,N_{z}}^{*},\cdots,U_{N_x,N_{y},N_{z}}^{*})^{T}.
\end{align*}
By noting ${\bm E}^{*},\ {\bm H}^{*}\in S_{N}^{'''}$, we can get
\begin{align}\label{SX_4.15}
&{\bm\Phi}^{n+\frac{1}{2}}=
\mu\delta_{t}({\bm H}^{*})^{n}+ {\bm D}({\bm E}^{*})^{n+\frac{1}{2}}
+\frac{c^{2}\tau^{2}}{12}{\bm D}^{3}({\bm E}^{*})^{n+\frac{1}{2}}
+\frac{c^{4}\tau^{4}}{120}{\bm D}^{5}({\bm E}^{*})^{n+\frac{1}{2}},\\\label{SX_4.16}
&{\bm\Psi}^{n+\frac{1}{2}}=\epsilon\delta_{t}({\bm E}^{*})^{n}
- {\bm D}({\bm H}^{*})^{n+\frac{1}{2}}
-\frac{c^{2}\tau^{2}}{12}{\bm D}^{3}({\bm H}^{*})^{n+\frac{1}{2}}
-\frac{c^{4}\tau^{4}}{120}{\bm D}^{5}({\bm H}^{*})^{n+\frac{1}{2}},
\end{align}
where the components of the above vectors
are the values of grid functions.

The scheme (\ref{SO-ECS1}) can be rewritten as
\begin{align}\label{SX_4.17}
&\mu\delta_{t}{\bm H}^{n}
+{\bm D}{\bm E}^{n+\frac{1}{2}}+\frac{c^{2}\tau^{2}}{12}{\bm D}^{3}{\bm E}^{n+\frac{1}{2}}
+\frac{c^{4}\tau^{4}}{120}{\bm D}^{5}{\bm E}^{n+\frac{1}{2}}=0,\\\label{SX_4.18}
&\epsilon\delta_{t}{\bm E}^{n}-{\bm D}{\bm H}^{n+\frac{1}{2}}
-\frac{c^{2}\tau^{2}}{12}{\bm D}^{3}{\bm H}^{n+\frac{1}{2}}
-\frac{c^{4}\tau^{4}}{120}{\bm D}^{5}{\bm H}^{n+\frac{1}{2}}=0.
\end{align}
Let ${ \mathcal{\bm H}}={\bm H}^{*}-{\bm H}$ and ${\mathcal{ \bm E}}
={\bm E}^{*}-{\bm E}$. Subtracting Eqs. \eqref{SX_4.15}-\eqref{SX_4.16} from Eqs. \eqref{SX_4.17}-\eqref{SX_4.18}, respectively,
we can obtain the error equations as follows
\begin{align}\label{SX_4.19}
&\mu\delta_{t}{\mathcal{H}}^{n}+ {\bm D}{\mathcal{E}}^{n}
+\frac{c^{2}\tau^{2}}{12}{\bm D}^{3}{\mathcal{E}}^{n+\frac{1}{2}}
+\frac{c^{4}\tau^{4}}{120}{\bm D}^{5}{\mathcal{E}}^{n+\frac{1}{2}}={\bm\Phi}^{n+\frac{1}{2}},\\\label{SX_4.20}
&\epsilon\delta_{t}{\mathcal{E}}^{n}- {\bm D}{\mathcal{H}}^{n+\frac{1}{2}}
-\frac{c^{2}\tau^{2}}{12}{\bm D}^{3}{\mathcal{H}}^{n+\frac{1}{2}}
-\frac{c^{4}\tau^{4}}{120}{\bm D}^{5}{\mathcal{H}}^{n+\frac{1}{2}}={\bm\Psi}^{n+\frac{1}{2}}.
\end{align}
Computing the inner product of (\ref{SX_4.19}) and (\ref{SX_4.20}) with ${\mathcal{H}}^{n+\frac{1}{2}}$ and ${\mathcal{ E}}^{n+\frac{1}{2}}$, respectively, we can have
\begin{align}\label{SX_4.21}
&\langle\mu\delta_{t}{\mathcal{ H}}^{n},{\mathcal{H}}^{n+\frac{1}{2}}\rangle_{h}
+\langle{\bm D}{\mathcal{E}}^{n+\frac{1}{2}},{\mathcal{H}}^{n+\frac{1}{2}}\rangle_{h}
+\frac{c^{2}\tau^{2}}{12}\langle{\bm D}^{3}{\mathcal{E}}^{n+\frac{1}{2}},{\mathcal{H}}^{n+\frac{1}{2}}\rangle_{h}\nonumber\\
&~~~~~~~~~~~~~~~~~~~~~~~~~~~
+\frac{c^{4}\tau^{4}}{120}\langle{\bm D}^{5}{\mathcal{E}}^{n+\frac{1}{2}},{\mathcal{H}}^{n+\frac{1}{2}}\rangle_{h}
=\langle {\bm\Phi}^{n+\frac{1}{2}},{\mathcal{H}}^{n+\frac{1}{2}}\rangle_{h},\\\label{SX_4.22}
&\langle\epsilon\delta_{t}{\mathcal{E}}^{n},{\mathcal{E}}^{n+\frac{1}{2}}\rangle_{h}
-\langle{\bm D}{\mathcal{H}}^{n+\frac{1}{2}},{\mathcal{E}}^{n+\frac{1}{2}}\rangle_{h}
-\frac{c^{2}\tau^{2}}{12}\langle{\bm D}^{3}{\mathcal{H}}^{n+\frac{1}{2}},
{\mathcal{E}}^{n+\frac{1}{2}}\rangle_{h}\nonumber\\
&~~~~~~~~~~~~~~~~~~~~~~~~~~~
-\frac{c^{4}\tau^{4}}{120}\langle{\bm D}^{5}{\mathcal{H}}^{n+\frac{1}{2}},
{\mathcal{E}}^{n+\frac{1}{2}}\rangle_{h}=\langle {\bm \Psi}^{n+\frac{1}{2}},{\mathcal{ E}}^{n+\frac{1}{2}}\rangle_{h}.
\end{align}
Adding (\ref{SX_4.21}) to (\ref{SX_4.22}) and making use of the complete square formulation, then, we can gain the following energy identity
\begin{align}\label{EI}
&\Arrowvert \sqrt{\mu}{\mathcal{H}}^{n+1}-\frac{\tau}{2\sqrt{\mu}}{\bm \Phi}^{n+\frac{1}{2}}\Arrowvert_{h}^{2}
+\Arrowvert \sqrt{\epsilon}{\mathcal{ E}}^{n+1}-\frac{\tau}{2\sqrt{\epsilon}}{\bm \Psi}^{n+\frac{1}{2}}\Arrowvert_{h}^{2}\nonumber\\
&~~~~~~~~~~~~~~~~~~~~~~~~~=\Arrowvert \sqrt{\mu}{\mathcal{ H}}^{n}
+\frac{\tau}{2\sqrt{\mu}}{\bm \Phi}^{n+\frac{1}{2}}\Arrowvert_{N}^{2}+\Arrowvert \sqrt{\epsilon}{\mathcal{ E}}^{n}
+\frac{\tau}{2\sqrt{\epsilon}}{\bm \Psi}^{n+\frac{1}{2}}\Arrowvert_{h}^{2}.
\end{align}
With the use of the triangle inequality of the norm and (\ref{EI}), we can get
\begin{align*}
&\Big(\Arrowvert \sqrt{\mu}{\mathcal{H}}^{n+1}\Arrowvert_{h}^{2}+\Arrowvert \sqrt{\epsilon}{\mathcal{E}}^{n+1}
\Arrowvert_{h}^{2}\Big)^{\frac{1}{2}}\\
&\leqslant\Big(\Arrowvert \sqrt{\mu}{\mathcal{H}}^{n+1}-\frac{\tau}{2\sqrt{\mu}}{\bm \Phi}^{n+\frac{1}{2}}\Arrowvert_{h}^{2}
+\Arrowvert \sqrt{\epsilon}{\mathcal{ E}}^{n+1}-\frac{\tau}{2\sqrt{\epsilon}}{\bm \Psi}^{n+\frac{1}{2}}\Arrowvert_{h}^{2}\Big)^{\frac{1}{2}}
+\frac{\tau}{2}\Big(\Arrowvert \frac{1}{\sqrt{\mu}}{\bm \Phi}^{n+\frac{1}{2}}\Arrowvert^{2}_{h}
+\Arrowvert \frac{1}{\sqrt{\epsilon}}{\bm \Psi}^{n+\frac{1}{2}}\Arrowvert^{2}_{h}\Big)^{\frac{1}{2}}\\
&=\Big(\Arrowvert \sqrt{\mu}{\mathcal{H}}^{n}+\frac{\tau}{2\sqrt{\mu}}{\bm \Phi}^{n+\frac{1}{2}}\Arrowvert_{h}^{2}
+\Arrowvert \sqrt{\epsilon}{\mathcal{ E}}^{n}
+\frac{\tau}{2\sqrt{\epsilon}}{\bm \Psi}^{n+\frac{1}{2}}\Arrowvert_{h}^{2}\Big)^{\frac{1}{2}}
+\frac{\tau}{2}\Big(\Arrowvert \frac{1}{\sqrt{\mu}}{\bm \Phi}^{n+\frac{1}{2}}\Arrowvert^{2}_{h}
+\Arrowvert \frac{1}{\sqrt{\epsilon}}{\bm \Psi}^{n+\frac{1}{2}}\Arrowvert^{2}_{h}\Big)^{\frac{1}{2}}\\
&\leqslant\Big(\Arrowvert \sqrt{\mu}{\mathcal{H}}^{n}\Arrowvert_{h}^{2}
+\Arrowvert \sqrt{\epsilon}{ \mathcal{ E}}^{n}\Arrowvert_{h}^{2}\Big)^{\frac{1}{2}}
+\tau\Big(\Arrowvert \frac{1}{\sqrt{\mu}}{\bm \Phi}^{n+\frac{1}{2}}\Arrowvert^{2}_{h}
+\Arrowvert \frac{1}{\sqrt{\epsilon}}{\bm \Psi}^{n+\frac{1}{2}}\Arrowvert^{2}_{h}\Big)^{\frac{1}{2}}.
\end{align*}
Recursively, applying the above inequality from time level $n-1$ to 0, we have
\begin{align}\label{Eq5}
&\Big(\Arrowvert \sqrt{\mu}{\mathcal{H}}^{n}\Arrowvert_{h}^{2}
+\Arrowvert \sqrt{\epsilon}{\mathcal{ E}}^{n}\Arrowvert_{h}^{2}\Big)^{\frac{1}{2}}
\leqslant\Big(\Arrowvert \sqrt{\mu}{\mathcal{H}}^{0}\Arrowvert_{h}^{2}+\Arrowvert
\sqrt{\epsilon}{\mathcal{ E}}^{0}\Arrowvert_{N}^{2}\Big)^{\frac{1}{2}}\nonumber\\
&~~~~~~~~~~~~~~~~~~~~~~~~~~~~~~~~~~~~~~~+\sum_{l=0}^{n-1}\tau\Big(\Arrowvert
 \frac{1}{\sqrt{\mu}}\Phi^{l+\frac{1}{2}}\Arrowvert^{2}_{h}
+\Arrowvert \frac{1}{\sqrt{\epsilon}}\Psi^{l+\frac{1}{2}}\Arrowvert^{2}_{h}\Big)^{\frac{1}{2}}.
\end{align}

By virtue of Lemma \ref{lem4.3}, and recalling that ${\bm H}^{0}={\bm H}(0), \ {\bm E}^{0}={\bm E}(0)$,
we have
\begin{align*}
&\Arrowvert{\mathcal{H}}^{0}\Arrowvert_{h}=\Arrowvert{\bm H}(0)-P_{N-2}{\bm H}(0)\Arrowvert_{h}\leqslant CN^{-r},\
\end{align*}
and
\begin{align*}
\Arrowvert{\mathcal{E}}^{0}\Arrowvert_{h}=\Arrowvert{\bm E}(0)-P_{N-2}{\bm E}(0)\Arrowvert_{h}\leqslant CN^{-r},
\end{align*}
which further implies that
\begin{align}\label{Eq6}
(\Arrowvert\sqrt{\mu}{\mathcal{H}}^{0}\Arrowvert_{h}^{2}+\Arrowvert\sqrt{\epsilon}{\mathcal{E}}^{0}
\Arrowvert_{h}^{2})^{\frac{1}{2}}\leqslant CN^{-r}.
\end{align}

With noting $n\tau<T$, we can deduce from (\ref{Eq5}) and (\ref{Eq6}) that
\begin{align}\label{ES2}
\Big(\Arrowvert \sqrt{\mu}{\mathcal{H}}^{n}\Arrowvert_{h}^{2}+\Arrowvert
 \sqrt{\epsilon}{\mathcal{ E}}^{n}\Arrowvert_{h}^{2}\Big)^{\frac{1}{2}}\leqslant CT(\tau^{6}+N^{-r}),\ n=0,\cdots,M.
\end{align}
Making use of (\ref{ES2}) and the  inequation, $a+b\leqslant\sqrt{2a^{2}+2b^{2}},\ {\forall}a,b\geqslant 0$, we can obtain
\begin{align}\label{SX_4.27}
\Arrowvert \sqrt{\mu}{\mathcal{H}}^{n}\Arrowvert_{h}+\Arrowvert
\sqrt{\epsilon}{ \mathcal{ E}}^{n}\Arrowvert_{h}\leqslant CT(\tau^{6}+N^{-r}),\ n=0,\cdots,M.
\end{align}
With Lemma \ref{lem4.3} and Eq. \eqref{SX_4.27}, the following error estimate can be established
\begin{align*}
&\Big(\mu\Arrowvert {\bm H}(t_{n})-{\bm H}^{n}\Arrowvert_{h}^{2}
+\epsilon\Arrowvert {\bm E}(t_{n})-{\bm E}^{n}\Arrowvert_{h}^{2}\Big)^{\frac{1}{2}}
\leqslant \sqrt{\mu}\Arrowvert {\bm H}(t_{n})-{\bm H}^{n}\Arrowvert_{h}
+\sqrt{\epsilon}\Arrowvert {\bm E}(t_{n})-{\bm E}^{n}\Arrowvert_{h}\\
&\leqslant \sqrt{\mu}\Arrowvert {\bm H}(t_{n})-P_{N-2}{\bm H}(t_{n})\Arrowvert_{h}
+\sqrt{\mu}\Arrowvert {\mathcal{H}}^{n}\Arrowvert_{h}
+\sqrt{\epsilon}\Arrowvert {\bm E}(t_{n})-P_{N-2}{\bm E}(t_{n})\Arrowvert_{h}
+\sqrt{\epsilon}\Arrowvert {\mathcal{E}}^{n}\Arrowvert_{h}\\
&\leqslant CT(\tau^{6}+N^{-r}), n=0,\cdots,M.
\end{align*}
This ends the proof.
\qed
\end{prf}
\section{Numerical dispersion relation}
In this section,
the numerical dispersion relation of the sixth order energy-conserved scheme (\ref{SO-ECS1}) will be investigated.
Let the elements of ${\bm u}$ satisfy $u_{j}=u_{j+N_{w}}$ with $u_{j}=u_{0}e^{-\text{i}\kappa_{w}jh_{w}}$,
  and denote $[({\bm D}_{1}^{w})^{p}]_{j,k}=(\hat{d}_{p}^{w})_{j,k},\ [{\bm D}_{p}^{w}]_{j,k}=(d_{p}^{w})_{j,k}$, where $j,k=1,\cdots,N_{w}$. With the help of Lemmas \ref{lem2.5_1}-\ref{lem2.5_2} and Eq. (\ref{CQ}), for a fixed positive integer $p$ , we can obtain the following results:\\
When $p$ is an odd integer, we have
\begin{align}\label{DF1}
[({\bm D}_{1}^{w})^{p}{\bm u}]_{j}&=\sum_{k=1}^{N_{w}}(\hat{d}_{p}^{w})_{j,k}u_{k}=\sum_{k=1}^{N_{w}}(\hat{d}_{p}^{w})_{j,j+k}u_{j+k}\nonumber\\
&=\sum_{k=1}^{\frac{N_{w}}{2}-1}(\hat{d}_{p}^{w})_{j,j+k}u_{j+k}+\sum_{k=\frac{N_{w}}{2}+1}^{N_{w}}(\hat{d}_{p}^{w})_{j,j+k}u_{j+k}\nonumber\\
&=\sum_{k=1}^{\frac{N_{w}}{2}-1}(\hat{d}_{p}^{w})_{j,j+k}u_{j+k}+\sum_{k=-\frac{N_{w}}{2}+1}^{0}(\hat{d}_{p}^{w})_{j,j+k}u_{j+k}\nonumber\\
&=\sum_{k=j-\frac{N_{w}}{2}+1}^{j+\frac{N_{w}}{2}-1}(\hat{d}_{p}^{w})_{j,k}u_{k}
=\sum_{k=j-\frac{N_{w}}{2}+1}^{j-1}(\hat{d}_{p}^{w})_{j,k}u_{k}
+\sum_{k=j+1}^{j+\frac{N_{w}}{2}-1}(\hat{d}_{p}^{w})_{j,k}u_{k}\nonumber\\
&=-2\text{i}u_{0}e^{-\text{i}\kappa_{w}jh_{w}}\sum_{k=1}^{\frac{N_{w}}{2}-1}(\hat{d}_{p}^{w})_{j,j+k}
\sin(kh_{w}\kappa_{w})\nonumber\\
&=-2\text{i}u_{0}e^{-\text{i}\kappa_{w}jh_{w}}\sum_{k=1}^{\frac{N_{w}}{2}-1}(d_{p}^{w})_{j,j+k}
\sin(kh_{w}\kappa_{w})=u_{0}e^{-\text{i}\kappa_{w}jh_{w}}\bar{d}_{p}^{w},
\end{align}
and
\begin{align}\label{DF2}
[{\bm D}_{p}^{w}{\bm u}]_{j}=-2\text{i}u_{0}e^{-\text{i}\kappa_{w}jh_{w}}\sum_{k=1}^{\frac{N_{w}}{2}-1}(d_{p}^{w})_{j,j+k}
\sin(kh_{w}\kappa_{w})=u_{0}e^{-\text{i}\kappa_{w}jh_{w}}\tilde{d}_{p}^{w},
\end{align}
where $\bar{d}_{p}^{w}=\tilde{d}_{p}^{w}=-2\text{i}\sum_{k=1}^{\frac{N_{w}}{2}-1}(d_{p}^{w})_{j,j+k}
\sin(kh_{w}\kappa_{w})$ is a pure imaginary number.\\
When $p$ is an even integer, by the analogous argument, we have
\begin{align}\label{DF3}
[({\bm D}_{1}^{w})^{p}{\bm u}]_{j}&=\sum_{k=1}^{N_{w}}(\hat{d}_{p}^{w})_{j,k}u_{k}\nonumber\\
&=u_{0}e^{-\text{i}\kappa_{w}jh_{w}}\Big[2\sum_{k=1}^{\frac{N_{w}}{2}-1}
(\hat{d}_{p}^{w})_{j,j+k}\cos(\kappa_{w}h_{w}k)\nonumber\\
&~~~~+(\hat{d}_{p}^{w})_{j,j}
+(\hat{d}_{p}^{w})_{j,j+\frac{N_{w}}{2}}\cos(\kappa_{w}h_{w}N_{w}/2)\Big]\nonumber\\
&=u_{0}e^{-\text{i}\kappa_{w}jh_{w}}\bar{d}_{p}^{w},
\end{align}
and
\begin{align}\label{DF4}
[{\bm D}_{p}^{w}{\bm u}]_{j}
&=u_{0}e^{-\text{i}\kappa_{w}jh_{w}}\Big[2\sum_{k=1}^{\frac{N_{w}}{2}-1}
(d_{p}^{w})_{j,j+k}\cos(\kappa_{w}h_{w}k)\nonumber\\
&~~~~+(d_{p}^{w})_{j,j}
+(d_{p}^{w})_{j,j+\frac{N_{w}}{2}}\cos(\kappa_{w}h_{w}N_{w}/2)\Big]\nonumber\\
&=u_{0}e^{-\text{i}\kappa_{w}jh_{w}}\tilde{d}_{p}^{w},
\end{align}
where $\bar{d}_{p}^{w}=2\sum_{k=1}^{\frac{N_{w}}{2}-1}
(\hat{d}_{p}^{w})_{j,j+k}\cos(\kappa_{w}h_{w}k)+(\hat{d}_{p}^{w})_{j,j}
+(\hat{d}_{p}^{w})_{j,j+\frac{N_{w}}{2}}\cos(\kappa_{w}h_{w}N_{w}/2)$ and
$\tilde{d}_{p}^{w}=2\sum_{k=1}^{\frac{N_{w}}{2}-1}
(d_{p}^{w})_{j,j+k}\cos(\kappa_{w}h_{w}k)+(d_{p}^{w})_{j,j}+(d_{p}^{w})_{j,j
+\frac{N_{w}}{2}}\cos(\kappa_{w}h_{w}N_{w}/2)$ are real numbers.
\begin{rmk} The numerical dispersion relation
of the proposed scheme (\ref{SO-ECS1}) will be established, when $\bar{d}_{p}^{w},\ p=1,2,3,4$ are obtained from Eqs. (\ref{DF1}) and (\ref{DF3}), respectively.
Further, by virtue of Eqs. (\ref{DF1})-(\ref{DF3}),
the numerical dispersion relation of the Fourier pseudo-spectral
method for generally linear Hamiltonian PDEs can also be well established.
\end{rmk}

Take a plane wave solution of the Maxwell's equations (\ref{ME}) as
\begin{align}\label{DA1}
\left(
\begin{array}{c}
              {\bm H}\\
             {\bm E}\\
\end{array}
\right)=\left(
\begin{array}{c}
              {\bm H}_{0}\\
             {\bm E}_{0}\\
\end{array}
\right)e^{-\text{i}(k_{x}x+k_{y}y+k_{z}z-\omega t)},
\end{align}
where ${\bm H}_{0}=((H_{x})_{0},(H_{y})_{0},(H_{z})_{0})^{T},
{\bm E}_{0}=((E_{x})_{0},(E_{y})_{0},(E_{z})_{0})^{T}$ denote an arbitrary constant vector, $\omega$ is the frequency and
$k_{x}$, $k_{y}$ and $k_{z}$ are the wave number along the $x$-direction, $y$-direction and $z$-direction, respectively.
Substituting (\ref{DA1}) into the Maxwell's equations Eq. (\ref{ME}) leads to the exact dispersion relation
\begin{align}\label{exact}
\omega^{2}=c^{2}(k_{x}^{2}+k_{y}^{2}+k_{z}^{2}).
\end{align}
Associated with the dispersion relation is two important quantities: the phase velocity $v_{p}$ and the group velocity $v_{g}$
\begin{align}
v_{p}=\frac{\omega}{|{\bm \kappa}|}{\hat{\bm \kappa}},\ v_{g}=\nabla_{{\bm \kappa}}\omega({\bm \kappa}),
\end{align}
where the wave number vector ${\bm \kappa}=(\kappa_{x},\kappa_{y},\kappa_{z})$,
and ${\hat{\bm \kappa}}=\frac{{\bm \kappa}}{|{\bm \kappa}|}$ is the unit vector.
 The phase velocity $v_{p}$ describes the speed at which the phase of a wave propagate and the
direction of the normal vector to the propagating wavefront. The group velocity $v_{g}$
describes the speed at which the envelope
of a wave packet propagates and gives the direction of the normal vector to the constant
$\omega$-surface of the dispersion relation \cite{WB99}.
Here, we discuss the three dimensional problem in ($x,y,z$). The exact group velocity
 $v_{g}$ and the wave number vector ${\bm\kappa}$ can be expressed as in spherical coordinates
\begin{align}\label{sc1}
&v_{g}=|v_{g}|(\sin\alpha\cos\beta,\sin\alpha\sin\beta,\cos\alpha),\\\label{sc2}
&{\bm \kappa}=(\kappa_{x},\kappa_{y},\kappa_{z})=|\kappa|(\sin\phi\cos\theta,\sin\phi\sin\theta,\cos\phi),
\end{align}
where $|v_{g}|$ and $\alpha,\beta$ are the magnitude and angles of $v_{g}$, respectively,
$|{\bm \kappa}|$ and $\theta,\phi$ are the magnitude and angles of ${\bm \kappa}$, respectively.

Next, we take the numerical solution of (\ref{ME}) to be
\begin{align}\label{DA2}
U_{w_{j,k,m}}^{n}=(U_{w})_{0}e^{-\text{i}(k_{x}jh_{x}+k_{y}kh_{y}+k_{z}mh_{z}-n\tau\omega)},\  U=H\  \text{or}\  E,
\end{align}
where
\begin{align*}
-\pi\leqslant h_{x}\kappa_{x}\leqslant\pi, -\pi\leqslant h_{y}\kappa_{y}\leqslant\pi,
-\pi\leqslant h_{z}\kappa_{z}\leqslant\pi, -\pi\leqslant \tau\omega\leqslant\pi.
\end{align*}
Let the increasing factor $\lambda=e^{\text{i}\omega \tau}$.
Substituting
\begin{align}\label{DA2}
U_{w_{j,k,m}}^{n}=(U_{w})_{0}\lambda^{n} e^{-\text{i}(k_{x}jh_{x}+k_{y}kh_{y}+k_{z}mh_{z})},\  U=H \ \text{or}\  E,
\end{align}
into the component form of the scheme (\ref{SO-ECS2}), then, we can obtain
\begin{align}\label{MEDA}
\lambda\left(
\begin{array}{cc}
              {\bm I}_{3\times 3}& c_{1}{\bm A}\\
             -c_{2}{\bm A}& {\bm I}_{3\times 3}\\
\end{array}
\right)\left(
\begin{array}{c}
              {\bm H}_{0}\\
             {\bm E}_{0}\\
\end{array}
\right)=\left(
\begin{array}{cc}
              {\bm I}_{3\times 3}& -c_{1}{\bm A}\\
             c_{2}{\bm A}& {\bm I}_{3\times 3}\\
\end{array}
\right)\left(
\begin{array}{c}
              {\bm H}_{0}\\
             {\bm E}_{0}\\
\end{array}
\right),
\end{align}
with $c_{1}=\frac{\tau}{2\mu}, c_{2}=\frac{\tau}{2\epsilon}$, and
\begin{align*}
{\bm A}=\left(
\begin{array}{ccc}
              0& -a_{z}&a_{y}\\
             a_{z}& 0&-a_{x}\\
             -a_{y}& a_{x} &0
\end{array}
\right),
\end{align*}
where
\begin{align}
&a_{x}=\bar{d}_{1}^{x}-\frac{c^{2}\tau^{2}}{12}\Big(\bar{d}_{3}^{x}
+\bar{d}_{1}^{x}(\bar{d}_{2}^{y}+\bar{d}_{2}^{z})\Big)\nonumber\\
&~~~~~~~+\frac{c^{4}\tau^{4}}{120}\Big(\bar{d}_{5}^{x}+\bar{d}_{1}^{x}(\bar{d}_{4}^{y}+2\bar{d}_{2}^{y}\bar{d}_{2}^{z}+\bar{d}_{4}^{z})
+2\bar{d}_{3}^{x}(\bar{d}_{2}^{y}+\bar{d}_{2}^{z})\Big),
\end{align}
 $a_{w}, w=y,z$ can be obtained by exchanging the index $x\longleftrightarrow w$ in $a_{x}$,
and
\begin{align*}
&{\bar d}_{1}^{w}=2\text{i}\sum_{k=1}^{\frac{N_{w}}{2}-1}\Big[\frac{(-1)^{k}\pi}{N_{w}h_{w}}\cot(\frac{k\pi}{N_{w}})\sin(kh_{w}\kappa_{w})\Big],\\
&\bar{d}_{2}^{w}=\Big[N_{w}-\frac{N_{w}^{2}+2}{3}+(-1)^{\frac{N_{w}}{2}}(N_{w}-2)
\cos(\frac{N_{w}}{2}h_{w}\kappa_{w})\Big]\Big(\frac{\pi}{N_{w}h_{w}}\Big)^{2}\nonumber\\
&~~~~~~~~+\sum_{k=1}^{\frac{N_{w}}{2}-1}(-1)^{k}\Big[\frac{2\pi^{2}}{N_{w}h_{w}^{2}}
-\Big(\frac{2\pi}{N_{w}h_{w}}\Big)^{2}\csc^{2}\Big(\frac{k\pi}{N_{w}}\Big)\Big]\cos(kh_{w}\kappa_{w}),\\
&\bar{d}_{3}^{w}=-2\text{i}\sum_{k=1}^{\frac{N_{w}}{2}-1}\frac{(-1)^{k}}{N_{w}}\Big(\frac{\pi}{h_{w}}\Big)^{3}\Big[\cot(\frac{k\pi}{N_{w}})
-\frac{6}{N_{w}^{2}}\cos(\frac{k\pi}{N_{w}})\csc^{3}(\frac{k\pi}{N_{w}})\Big]\sin(kh_{w}\kappa_{w}),\\
&\bar{d}_{4}^{w}=\Big(\frac{\pi}{h_{w}}\Big)^{4}\Big[\frac{1}{5}-\frac{1}{N_{w}}+\frac{4}{3N_{w}^{2}}-\frac{8}{15N_{w}^{4}}+(-1)^{\frac{N_{w}}{2}}
\Big(\frac{4}{N_{w}^{2}}-\frac{8}{N_{w}^{4}}-\frac{1}{N_{w}}\Big)\cos\Big(\frac{N_{w}}{2}h_{w}\kappa_{w}\Big)\Big]\nonumber\\
&~~~~+2\sum_{k=1}^{\frac{N_{w}}{2}-1}(-1)^{k}\Big(\frac{\pi}{h_{w}}\Big)^{4}
\Big[\csc^{2}(\frac{k\pi}{N_{w}})\Big(\frac{4}{N_{w}^{2}}-\frac{8}{N_{w}^{4}}
-\frac{24}{N_{w}^{4}}\cot^{2}(\frac{k\pi}{N_{w}})\Big)-\frac{1}{N_{w}}\Big]\cos(kh_{w}\kappa_{w}),\\
&\bar{d}_{5}^{w}=2\text{i}\sum_{k=1}^{\frac{N_{w}}{2}-1}(-1)^{k}\Big(\frac{\pi}{h_{w}}\Big)^{5}\cot(\frac{k\pi}{N_{w}})
\Big[\frac{1}{N_{w}}+\frac{20}{N_{w}^{5}}\csc^{2}(\frac{k\pi}{N_{w}})\Big(4
+6\cot^{2}(\frac{k\pi}{N_{w}})-N_{w}^{2}\Big)\Big]\\
&~~~~~~~\sin(kh_{w}\kappa_{w}).
\end{align*}
Since Eq. (\ref{MEDA}) holds for nonzero $U_{w_{0}}$, this implies that
\begin{align}
|\lambda {\bm I}_{6\times 6}-\mathcal{A}^{-1}(2{\bm I}_{6\times 6}-\mathcal{A})|=0,
\end{align}
where
\begin{align*}
\mathcal{A}=\left(
\begin{array}{cc}
              {\bm I}_{3\times 3}& c_{1}{\bm A}\\
             -c_{2}{\bm A}& {\bm I}_{3\times 3}\\
\end{array}
\right).
\end{align*}
By careful calculation, we have
\begin{align}\label{MEDA2}
&\lambda_{1,2}=1,\ \lambda_{3,4}=\frac{1+\text{i}\sqrt{c_{1}c_{2}}\sqrt{|a_{x}|^{2}+|a_{y}|^{2}+|a_{z}|^{2}}}
{1-\text{i}\sqrt{c_{1}c_{2}}\sqrt{|a_{x}|^{2}+|a_{y}|^{2}+|a_{z}|^{2}}},\nonumber\\
&\lambda_{5,6}=\frac{1-\text{i}\sqrt{c_{1}c_{2}}\sqrt{|a_{x}|^{2}+|a_{y}|^{2}+|a_{z}|^{2}}}
{1+\text{i}\sqrt{c_{1}c_{2}}\sqrt{|a_{x}|^{2}+|a_{y}|^{2}+|a_{z}|^{2}}}.
\end{align}
It follows from \eqref{MEDA2} that the modulus of the characteristic values are equal to one. Therefore, we have the following theorem.
\begin{thm}
The proposed sixth order energy-conserved scheme is unconditionally stable and non-dissipative.
\end{thm}
Sequentially, inserting $\lambda=e^{\text{i}\omega\tau}$ into the characteristic values (\ref{MEDA2})
results in numerical dispersion relation of the scheme (\ref{SO-ECS1})
\begin{align}\label{NDA}
 \tan^{2}\frac{\omega \tau}{2}=\frac{\tau^{2}c^{2}}{4}\left(|a_{x}|^{2}+|a_{y}|^{2}+|a_{z}|^{2}\right).
\end{align}
It is remarked that the dispersion relation (\ref{NDA}) converges to the theoretical dispersion relation (\ref{exact})
provided that $\tau, h_{x}, h_{y}, h_{z}\rightarrow 0$.

Then, the normalized phase velocity is given by
\begin{align}
&\Big| \frac{v_{p}}{c}\Big|=\Big| \frac{\omega}{c\vert\kappa\vert}\Big|,
\end{align}
It is obvious that the normalized phase velocity of the Maxwell's equations (\ref{ME}) is equal to one.
Supposing $N_{x}=N_{y}=N_{z}=N$, $S=c\frac{\tau}{h}$ and the number of points per wavelength $N_{\lambda}=\frac{\lambda}{h}$ ($\lambda$ is the wavelength), the normalized phase velocity of the scheme (\ref{SO-ECS1}) can be expressed as
\begin{align}
&\vert\frac{v_{p}}{c}\vert=\frac{N_{\lambda}}{\pi S}\tan^{-1}\left(\frac{1}{2}\sqrt{|c\tau a_{x}|^{2}+|c\tau a_{y}|^{2}+|c\tau a_{z}|^{2}}\right).
\end{align}

Finally, we calculate the magnitude of the group velocity $v_{g}$ by
\begin{align}\label{vg1}
|v_{g}|=\sqrt{(v_{g})_{x}^{2}+(v_{g})_{y}^{2}+(v_{g})_{z}^{2}}=\sqrt{\left(\frac{\partial \omega}{\partial \kappa_{x}}\right)^{2}+\left(\frac{\partial \omega}{\partial \kappa_{y}}\right)^{2}+\left(\frac{\partial \omega}{\partial \kappa_{z}}\right)^{2}}.
\end{align}
It's clear that the exact group velocity of the Maxwell's equations (\ref{ME}) is equal to one. From numerical dispersion relation
(\ref{NDA}), we have
\begin{align}\label{v_g1}
&\frac{\partial \omega}{\partial \kappa_{w}}=\mathbb{A}
\frac{\partial(|a_{x}|^{2}+|a_{y}|^{2}+|a_{z}|^{2})}{\partial \kappa_{w}}
=-\mathbb{A}\left(\frac{\partial a_{x}^{2}}{\partial \kappa_{w}}+\frac{\partial a_{y}^{2}}{\partial
 \kappa_{w}}+\frac{\partial a_{z}^{2}}{\partial \kappa_{w}}\right),
\end{align}
where
\begin{align}\label{v_g2}
\mathbb{A}=\frac{c}{2}\frac{1}
{[1+(\frac{c\tau}{2})^2(|a_{x}|^{2}+|a_{y}|^{2}+|a_{z}|^{2})]\sqrt{|a_{x}|^{2}+|a_{y}|^{2}+|a_{z}|^{2}}},
\end{align}
and
\begin{align}\label{v_g3}
&\frac{\partial a_{x}^{2}}{\partial \kappa_{x}}=2a_{x}\Big\{\frac{\partial \bar{d}_{1}^{x}}{\partial \kappa_{x}}
-\frac{c^2\tau^2}{12}\Big[\frac{\partial \bar{d}_{3}^{x}}{\partial \kappa_{x}}
+\frac{\partial \bar{d}_{1}^{x}}{\partial \kappa_{x}}(\bar{d}_{2}^{y}+\bar{d}_{2}^{z})\Big]\nonumber\\
&~~~~~~~~+
\frac{c^4\tau^4}{120}\Big[\frac{\partial \bar{d}_{5}^{x}}{\partial \kappa_{x}}+\frac{\partial \bar{d}_{1}^{x}}{\partial \kappa_{x}}(\bar{d}_{4}^{y}+2\bar{d}_{2}^{y}\bar{d}_{2}^{z}+\bar{d}_{4}^{z}
)+2\frac{\partial \bar{d}_{3}^{x}}{\partial \kappa_{x}}(\bar{d}_{2}^{y}
+\bar{d}_{2}^{z})\Big]\Big\},
\end{align}
\begin{align}\label{v_g4}
&\frac{\partial a_{x}^{2}}{\partial \kappa_{y}}=2a_{x}\Big\{-\frac{c^2\tau^2}{12}\bar{d}_{1}^{x}\frac{\partial \bar{d}_{2}^{y}}{\partial \kappa_{y}}
+\frac{c^4\tau^4}{120}\Big[(2\bar{d}_{1}^{x}\bar{d}_{2}^{z}+2\bar{d}_{3}^{x})\frac{\partial \bar{d}_{2}^{y}}{\partial \kappa_{y}}
+\bar{d}_{1}^{x}\frac{\partial \bar{d}_{4}^{y}}{\partial \kappa_{y}}\Big]\Big\}.
\end{align}
Others can be obtained by the following tricks
\begin{enumerate}
\item $\frac{\partial a_{w}^{2}}{\partial \kappa_{w}}, w=y,z$ are obtained by exchanging
 $x\longleftrightarrow w$ in Eq. (\ref{v_g3}),
\item $\frac{\partial a_{y}^{2}}{\partial \kappa_{x}},\  \frac{\partial a_{x}^{2}}{\partial \kappa_{z}}$ and
$\frac{\partial a_{z}^{2}}{\partial \kappa_{y}}$ are obtained by exchanging $y\longleftrightarrow x$, $y\longleftrightarrow z$ and
$x\longleftrightarrow z$ in Eq. (\ref{v_g4}), respectively.
\item $\frac{\partial a_{z}^{2}}{\partial \kappa_{x}}$ and $\frac{\partial a_{y}^{2}}{\partial \kappa_{z}}$ are obtained by
exchanging $y\longleftrightarrow z$ in
$\frac{\partial a_{y}^{2}}{\partial \kappa_{x}}$ and $x\longleftrightarrow y$ in $\frac{\partial a_{x}^{2}}{\partial \kappa_{z}}$, respectively.
\end{enumerate}
Thus, the numerical group velocity of the scheme (\ref{SO-ECS1})
can be obtained by inserting $\frac{\partial \omega}{\partial \kappa_{w}}, w=x,y,z$ into Eq. (\ref{vg1}).
\section{Numerical experiments}
In this section, we will investigate the numerical behavior of the proposed method (\ref{SO-ECS1}) presented in
Section 3 and analyze the numerical dispersion relation derived in Section 4.
All diagrams presented below refer to the numerical integration of the Maxwell's equations (\ref{ME})
with $\mu=\epsilon=1$.
\subsection{Benchmark test}
In this subsection, we will focus on the rate of convergence, the two divergence-free properties as well as the discrete energy and momentum conservation laws.
Furthermore, some comparisons between our scheme and the composition scheme which is derived from the composition method (see Page. 44 in Ref. \cite{ELW06}) via the AVF scheme in Ref. \cite{CWG16} and the Gauss scheme obtained by using the Gauss methods of order 6 (see Page. 34 in Ref. \cite{ELW06}) for Eq. \eqref{MESH1} will be presented.
In our computations, a fast solver is employed to solve the linear system (\ref{SO-ECS1}) efficiently.
 For more details of the fast solver, please refer to Appendix.
We assume that the domain is $\Omega=[0,2]\times[0,2]\times[0,2]$ with periodic boundary conditions.
The exact solutions of Eqs. (\ref{ME}) are \cite{CLL08}
\begin{align}\label{ES}
&E_{x}=\frac{k_{y}-k_{z}}{\epsilon \sqrt{\mu} w}\cos(w\pi t)\cos(k_{x}\pi x)\sin(k_{y} \pi y)\sin(k_{z}\pi z),\\
&H_{x}=\sin(w\pi t)\sin(k_{x}\pi x)\cos(k_{y} \pi y)\cos(k_{z}\pi z),\\
&E_{y}=\frac{k_{z}-k_{x}}{\epsilon \sqrt{\mu} w}\cos(w\pi t)\sin(k_{x}\pi x)\cos(k_{y} \pi y)\sin(k_{z}\pi z),\\
&H_{y}=\sin(w\pi t)\cos(k_{x}\pi x)\sin(k_{y} \pi y)\cos(k_{z}\pi z),\\
&E_{z}=\frac{k_{x}-k_{y}}{\epsilon \sqrt{\mu} w}\cos(w\pi t)\sin(k_{x}\pi x)\sin(k_{y} \pi y)\cos(k_{z}\pi z),\\
&H_{z}=\sin(w\pi t)\cos(k_{x}\pi x)\cos(k_{y} \pi y)\sin(k_{z}\pi z),
\end{align}
where $k_{x}=1,k_{y}=2,k_{z}=-3$ and $w=(k_{x}^{2}+k_{y}^{2}+k_{z}^{2})/(\epsilon\mu)$.

We set the exact solutions (6.1-6.6) at $t=0$ on the domain $\Omega=[0,2]\times[0,2]\times[0,2]$ as initial
conditions. In order to evaluate the numerical errors, the following formulas are used
\begin{align}
L^{\infty}=\max\Big\{\Arrowvert \mu{\bm H}(t_{n})-\mu{\bm H}^{n}\Arrowvert_{h,\infty},\Arrowvert \epsilon{\bm E}(t_{n})-\epsilon{\bm E}^{n}\Arrowvert_{h,\infty}\Big\},
\end{align}
and
\begin{align}
L^{2}=\left(\mu\Arrowvert{\bm H}(t_{n})-{\bm H}^{n}\Arrowvert_{h}^{2}+\epsilon\Arrowvert{\bm E}(t_{n})-{\bm E}^{n}\Arrowvert_{h}^{2}\right)^{\frac{1}{2}}.
\end{align}
The rate of convergence is defined as
\begin{align}
\text{Rate}=\frac{\ln( error_{1}/error_{2})}{\ln (\tau_{1}/\tau_{2})},
\end{align}
where $\tau_{l}, error_{l}, (l=1,2)$ are step sizes and errors with the step size $\tau_{l}$, respectively.


\begin{table}[H]
\tabcolsep=9pt
\small
\renewcommand\arraystretch{1.3}
\centering
\caption{The temporal errors and rates of convergence of different schemes with different time steps and $N_{x}=N_{y}=N_{z}=16$ at $T=1$.}\label{Tab611}
\begin{tabularx}{1.2\textwidth }{XXXXXXX}\hline
{Scheme}&{$\tau$}  &{$L^{\infty}$} & {Rate}  & {$L^{2}$} & {Rate}&{CPU (s)} \\     
\hline
 {\eqref{SO-ECS1}} &{0.01}  & {2.5295e-08} & {-} & {4.5198e-08}&{-}&{ 0.6}\\
 {}&{0.005} & {3.9585e-10} & {5.9978}  & {7.0731e-10} & {5.9978}&{1.0}  \\   
 {}&{0.0025}  & {6.2084e-12} & {5.9946} & {1.1053e-11}&{5.9998}&{1.7} \\
 {}&{0.00125} & {1.9063e-13} & {-}  & {2.7046e-13} & {-}&{2.9}  \\
 {composition}&{0.01}  & {7.3383e-08} & {-} & {1.3112e-07}&{-}&{1.3}\\
  {}&{0.005} & {1.1481e-09} & {5.9981}  & { 2.0515e-09} & {5.9981}&{2.0}  \\   
 {}&{0.0025}  & {1.8013e-11} & {5.9941} & {3.2060e-11}&{5.9999}&{3.3} \\
 {}&{0.00125} & {5.3180e-13} & {-} & {6.2894e-13} & {-}  &{6.1}\\
 {Gauss}&{0.01}  & { 2.9804e-10} & {-} & { 5.3255e-10}&{-}&{2.0}\\
  {}&{0.005} & {4.6428e-12} & {6.004}  & { 8.3228e-12} & {5.9997}&{2.7}  \\   
 {}&{0.0025}  & {9.3703e-14} & {5.6308} & {1.4772e-13}&{5.8161}&{3.9} \\
 {}&{0.00125} & {8.3932e-14} & {-} & {1.4881e-13} & {-}  &{6.2}\\
\hline
\end{tabularx}
\end{table}

\begin{table}[H]
\tabcolsep=9pt
\small
\renewcommand\arraystretch{1.2}
\centering
\caption{The spatial error of the proposed scheme with different spatial steps and $\tau=10^{-3}$ at $T=1$.}\label{Tab612}
\begin{tabularx}{\textwidth}{XXXXXXX}\hline
{$N_x\times N_y\times N_z$}  &{$L^{\infty}$} & {$L^{2}$}&{CPU (s)} \\     
\hline
 {$8\times 8 \times 8$}  & {8.8374e-13} &{1.0989e-13}&{1.5}\\
 {$16\times 16 \times 16$} & {1.1513e-13} & {1.3643e-13}&{3.8}   \\   
 {$32\times 32 \times 32$}  & {1.7786e-13} & {2.2746e-13}&{33.1} \\
 {$64\times 64 \times 64$}  & {1.9196e-13}  & {2.4856e-13}&{450.1}\\
\hline
\end{tabularx}
\end{table}

Table \ref{Tab611} shows the temporal errors and rates of convergence of different schemes with different time steps
and $N_{x}=N_{y}=N_{z}=16$ at $T=1$.
It is clear that all schemes can reach the theoretical order of 6, which verifies the theoretical analysis. Compared with the composition scheme, the proposed scheme shows
distinct advantage in the errors and CPU time.
 The Gauss scheme provides smaller error than our scheme, but our scheme is the most efficient. In addition, since the Gauss scheme contains triple intermediate variables, the memory requirements of the Gauss scheme are at least triple as large as our scheme. We display the spatial error of the scheme (\ref{SO-ECS1}) in Table \ref{Tab612}. As illustrated in  Table \ref{Tab612},
the spatial error of the proposed scheme is negligible and up to machine precision, which confirms that, for
sufficiently smooth problems, the Fourier pseudo-spectral method is of arbitrary order
of accuracy.  We omit the spatial error of the composition and Gauss schemes because they are also discretized by the Fourier pseudo-spectral method in sapce.

Then, we investigate the errors of the energy and momentum conservation laws over the time interval $t\in[0,100]$ 
with $N_{x}=N_{y}=N_{z}=16$ and $\tau=0.01$ in Fig. \ref{fig_611}. It is clear to see that the errors on the energy and momentum of the Gauss scheme is the smallest, whereas, the ones of the composition scheme is the largest. It's worth noting that energy invariants
$\mathcal{E}_{3}$ and $\mathcal{E}_{5x}$ are larger than others,
we suspect that it is probably due to the fact that they have been divided by the time size $\tau$ at each time level.
The errors of the invariants $\mathcal{E}_{4w}$, $\mathcal{E}_{5w}$ and  $\mathcal{M}_{w},\ w=y,z$ are not presented since they are close to the case of $x$.

Finally, we check the two discrete divergence-free fields over the time interval $t\in[0,100]$ with $N_{x}=N_{y}=N_{z}=16$ and $\tau=0.01$ in Fig. \ref{fig_612}.
Numerical results show that the two discrete divergence-free fields reach machine precision in long time computation,
which conforms the fact that, with appropriate collocation points, our scheme can preserve the two divergence-free fields exactly.

\begin{figure}[H]
\centering\begin{minipage}[t]{60mm}
\includegraphics[width=60mm]{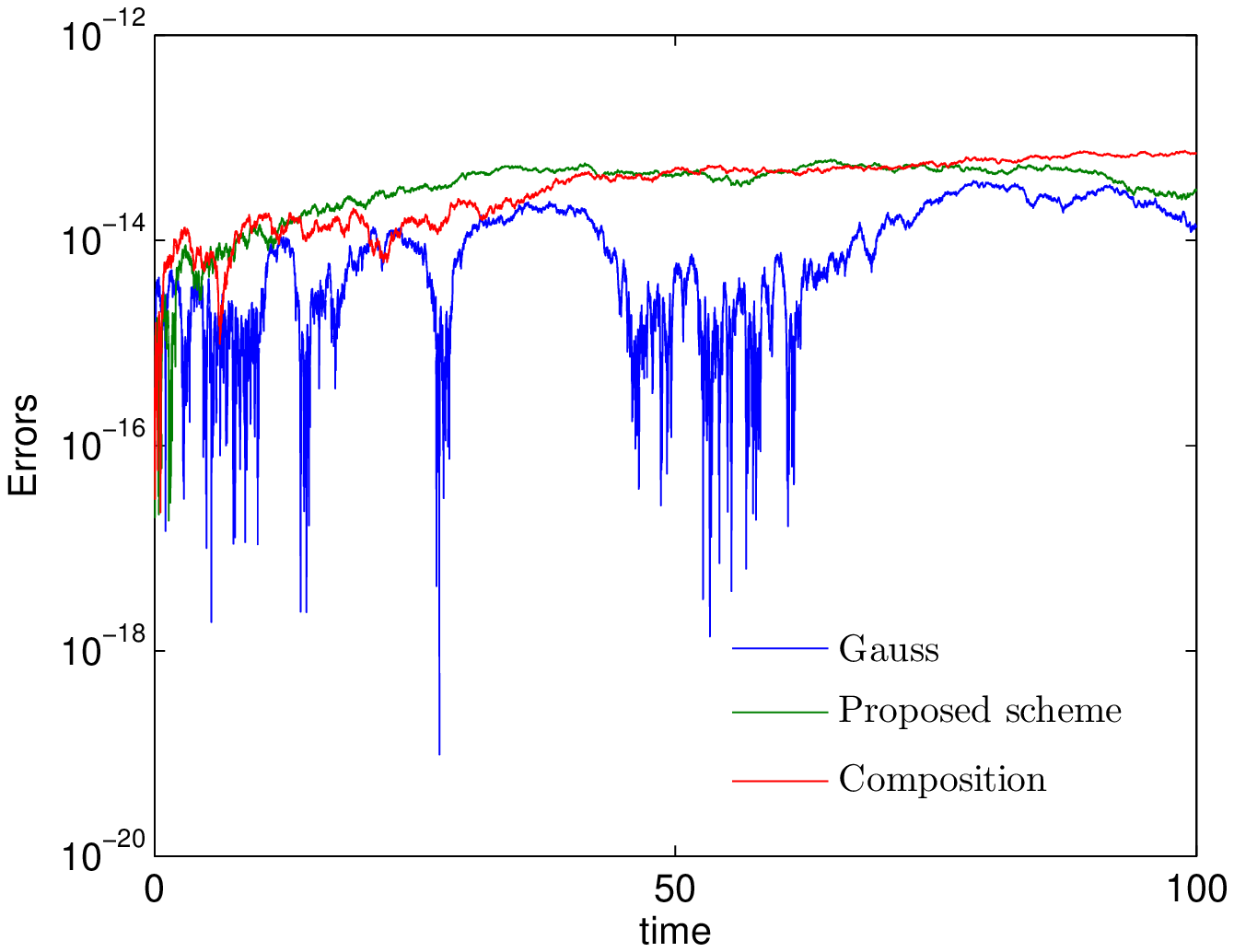}
{\footnotesize  (a) Errors in invariant $\mathcal{E}_{1}$.}
\end{minipage}\ \
\begin{minipage}[t]{60mm}
\includegraphics[width=60mm]{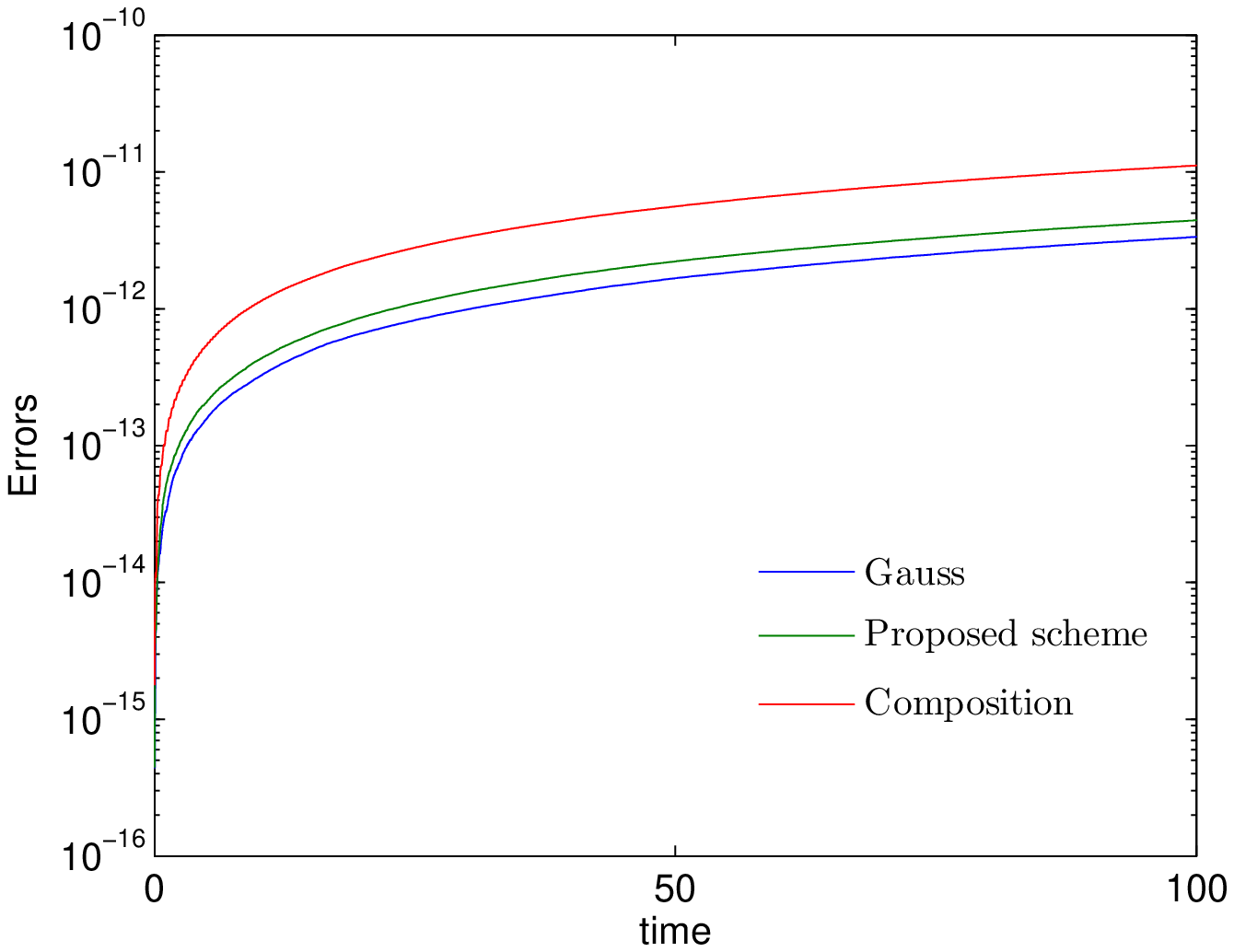}\\
{\footnotesize  (b) Errors in invariant $\mathcal{E}_{2}$.}
\end{minipage}
\end{figure}
\begin{figure}[H]
\centering\begin{minipage}[t]{60mm}
\includegraphics[width=60mm]{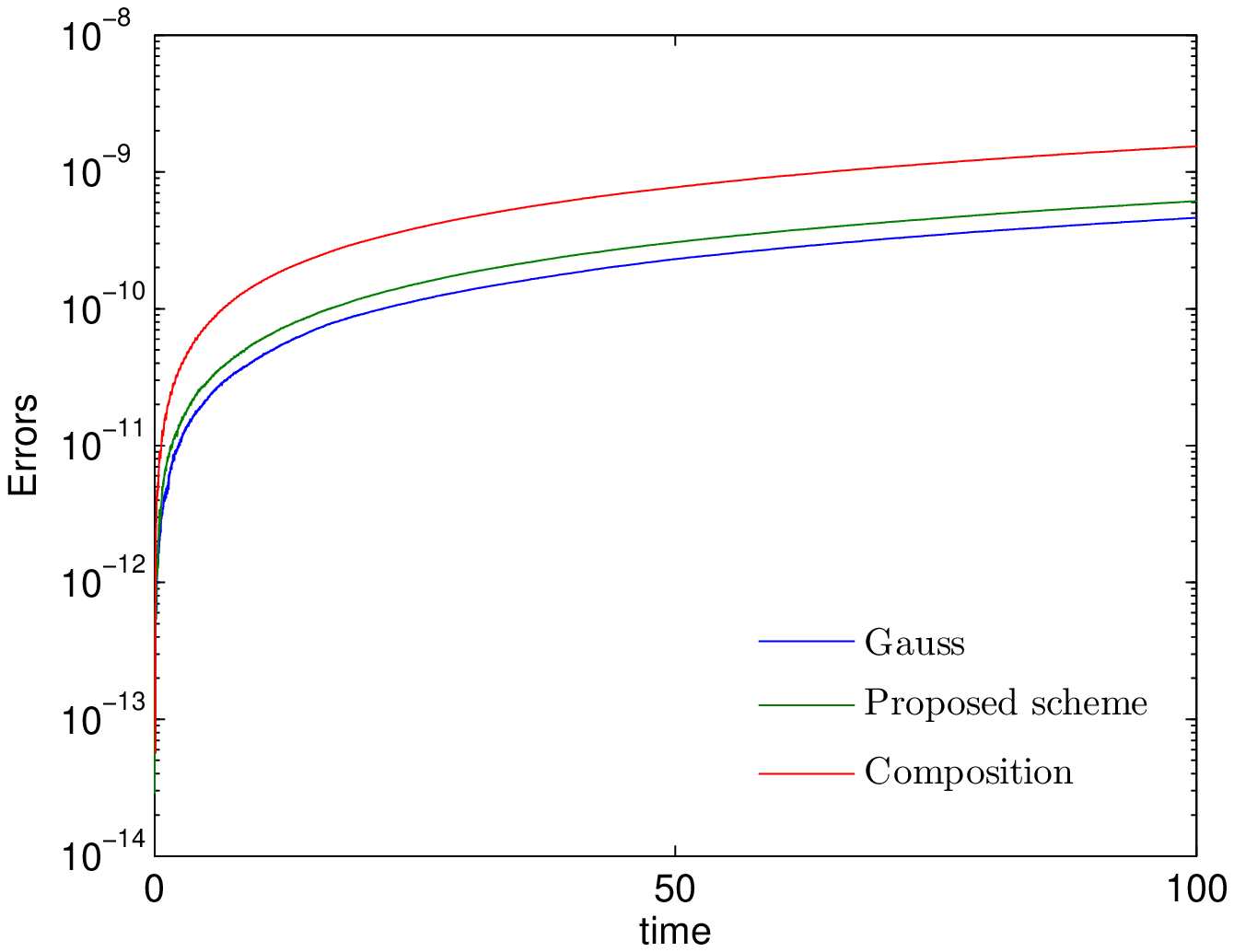}\\
{\footnotesize  (c) Errors in invariant $\mathcal{E}_{3}$.}
\end{minipage}\ \
\begin{minipage}[t]{60mm}
\includegraphics[width=60mm]{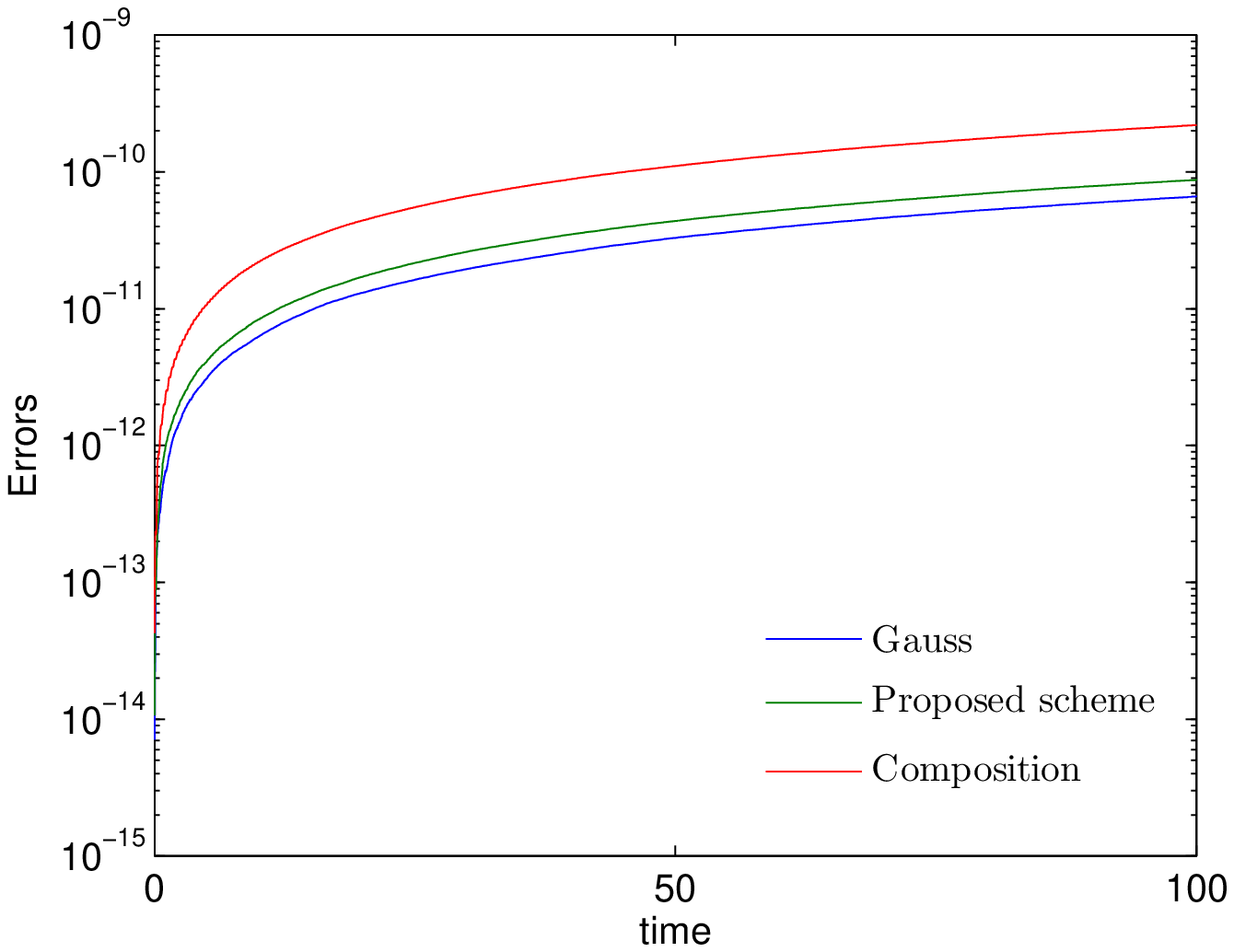}
{\footnotesize  (d) Errors in invariant $\mathcal{E}_{4x}$.}
\end{minipage}
\end{figure}

\begin{figure}[H]
\centering\begin{minipage}[t]{60mm}
\includegraphics[width=60mm]{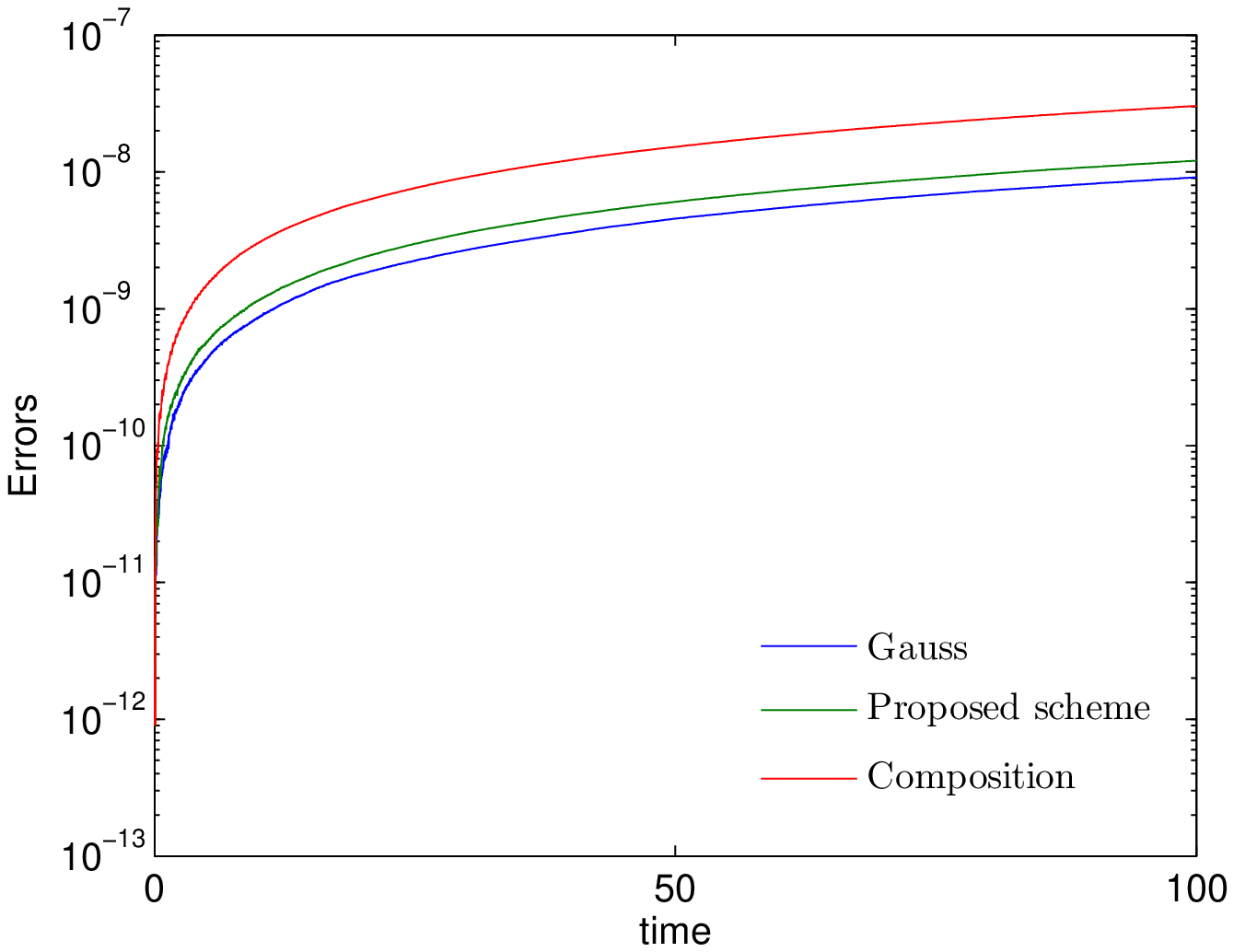}\\
{\footnotesize  (e) Errors in invariantinvariant $\mathcal{E}_{5x}$.}
\end{minipage}
\begin{minipage}[t]{60mm}
\includegraphics[width=60mm]{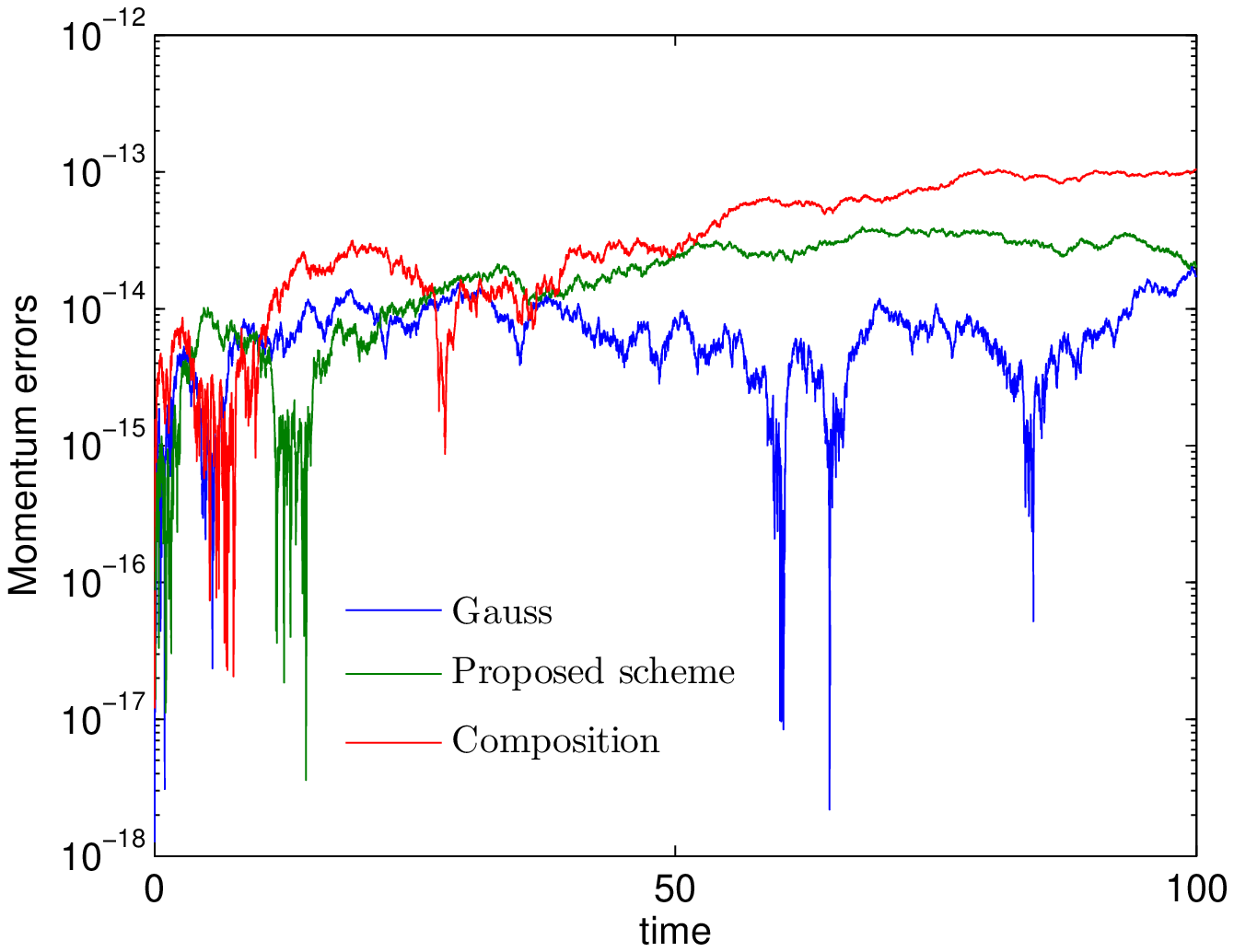}
{\footnotesize (f) Errors in invariant$\mathcal{M}_{x}$.}
\end{minipage}
\caption{The errors in invariants over the time interval $t\in[0,100]$ with $N_x=N_y=N_z=16$ and $\tau=0.01$.}\label{fig_611}
\end{figure}

\begin{figure}[H]
\centering\begin{minipage}[t]{60mm}
\includegraphics[width=60mm]{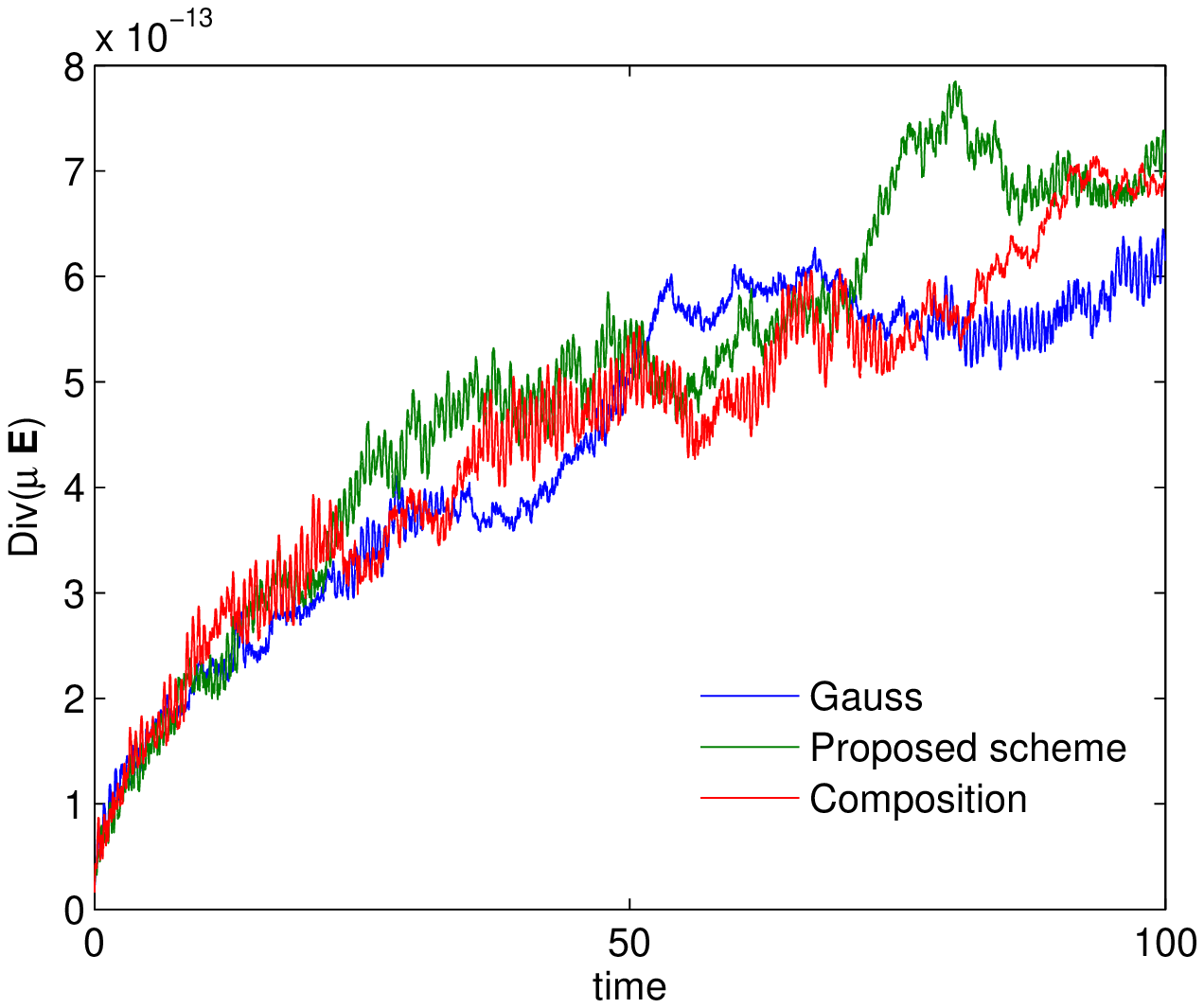}
\end{minipage}\ \
\begin{minipage}[t]{60mm}
\includegraphics[width=60mm]{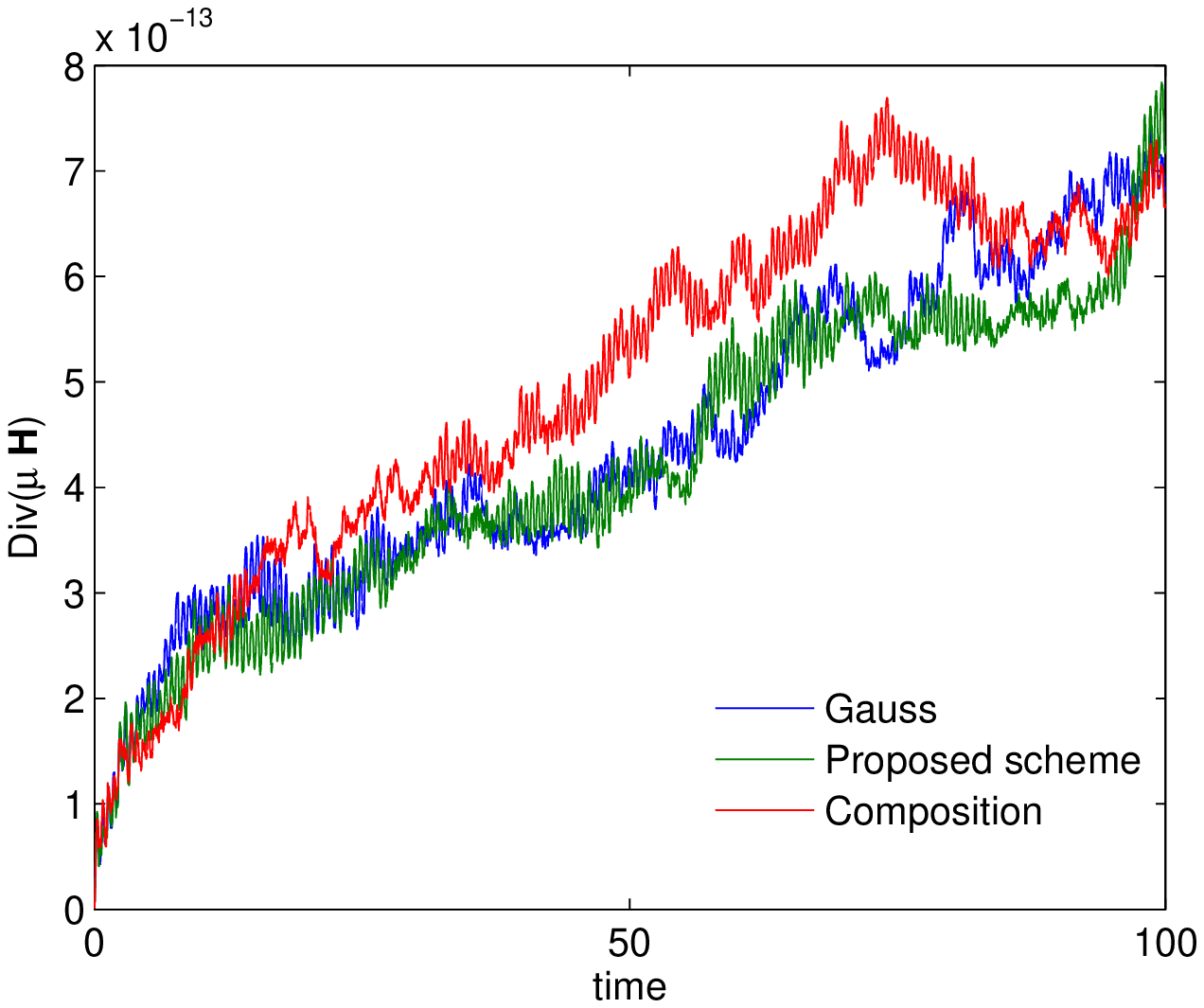}
\end{minipage}
\caption{The two discrete divergence-free fields over the time interval $t\in[0,100]$ with $N_x=N_y=N_z=16$ and $\tau=0.01$.}\label{fig_612}
\end{figure}

\subsection{Numerical dispersion analysis }
In this subsection, we will investigate the numerical dispersion relation
of the proposed method (\ref{SO-ECS1}) including the numerical phase velocity and the numerical group velocity.
In the following discussion, for the sake of simplicity,
only uniform cell is considered, i.e., $N_{x}=N_{y}=N_{z}=N=150$.
\subsubsection{Normalized numerical phase velocity}
\begin{figure}[H]
\centering\begin{minipage}[t]{60mm}
\includegraphics[width=60mm]{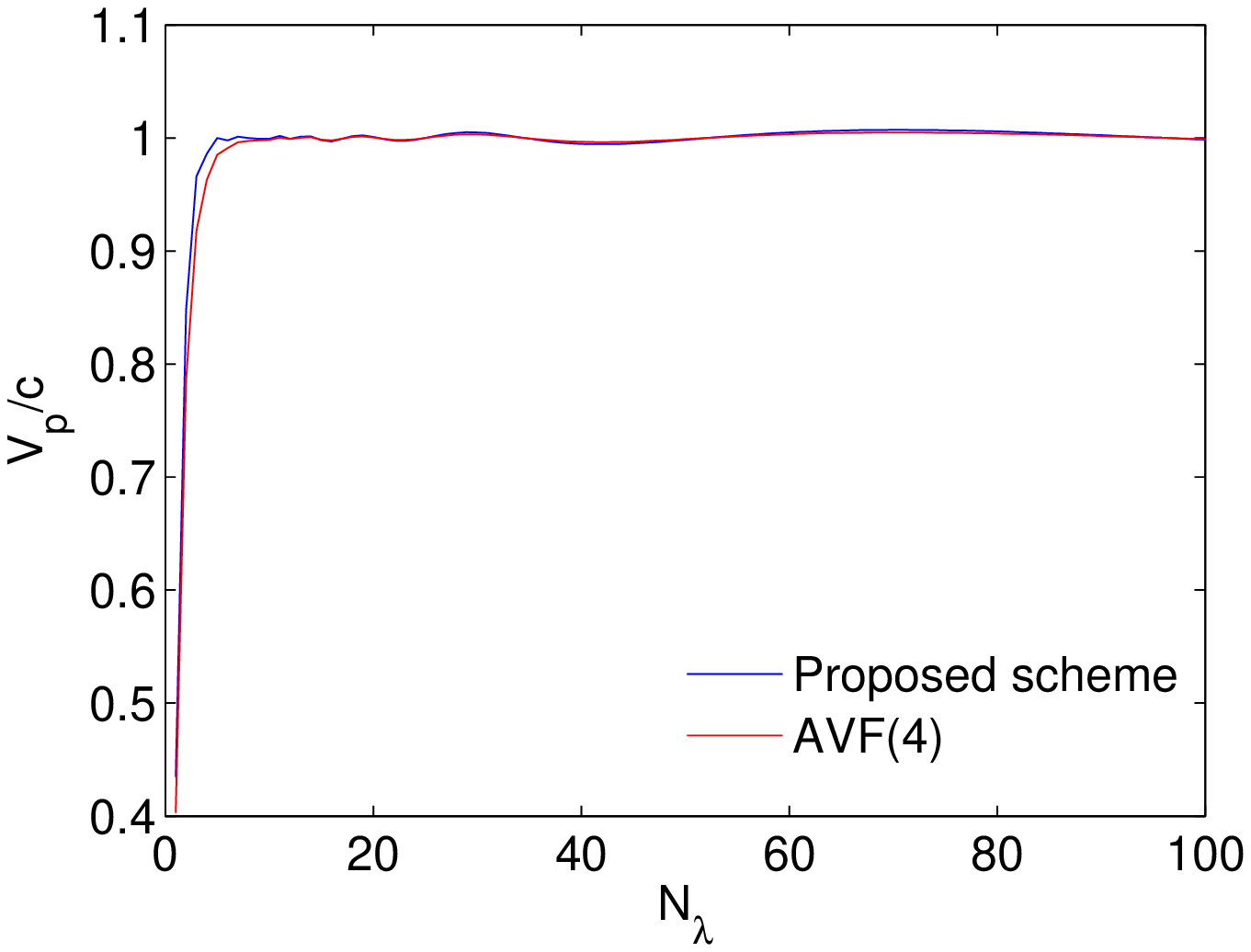}
{\footnotesize (a)}
\end{minipage}\
\begin{minipage}[t]{60mm}
\includegraphics[width=60mm]{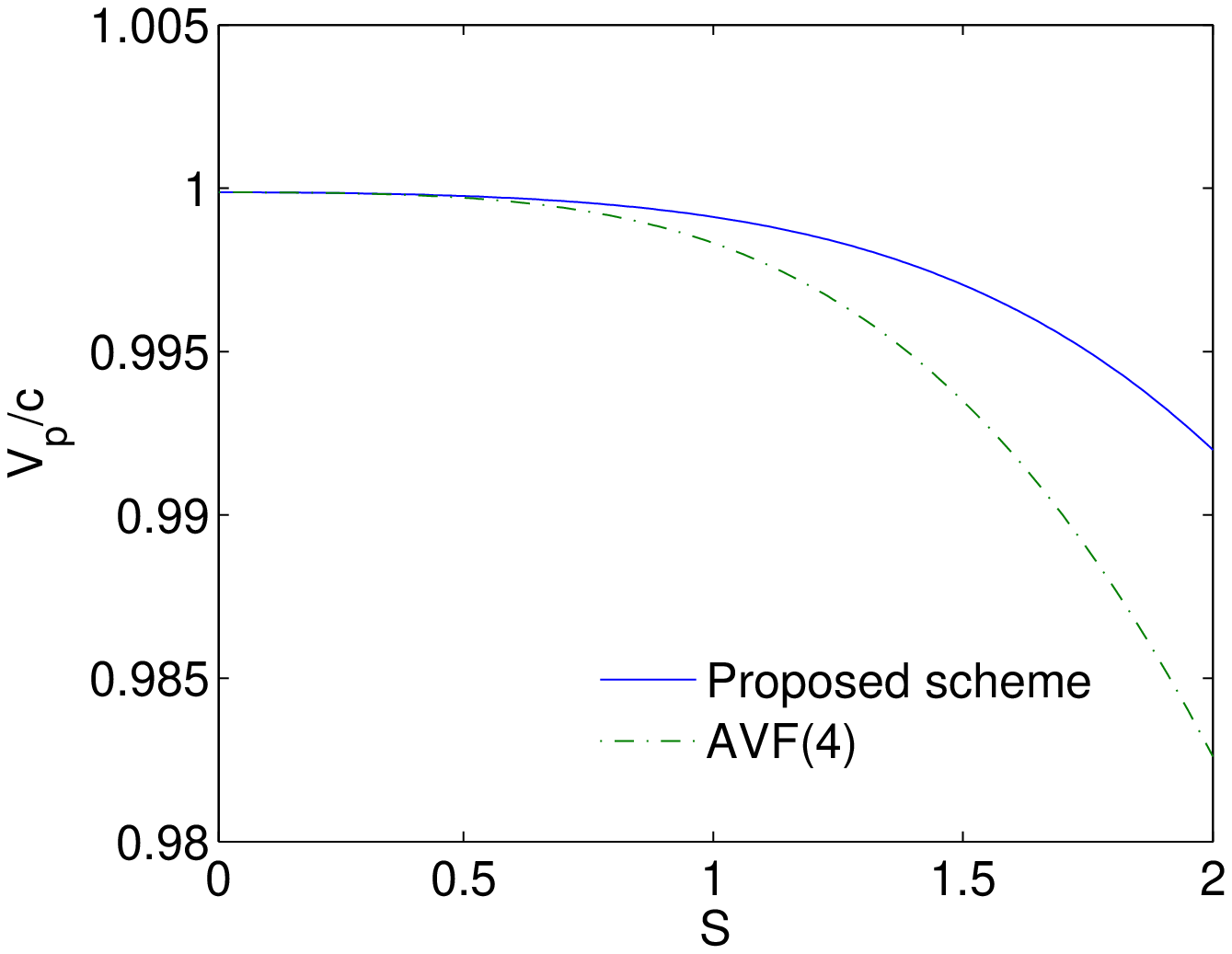}
{\footnotesize (b)}
\end{minipage}\\
\caption{(a) Normalized numerical phase velocities versus the number of points per
wavelength $N_{\lambda}$ with  $\phi=\pi/4$ and $\theta=3\pi/8$.
(b) Normalized numerical phase velocities versus the CFL number $S$ with $N_{\lambda}=5$ with $\phi=\pi/4$ and $\theta=3\pi/8$.} \label{fig_614}
\end{figure}

\begin{figure}[H]
\centering\begin{minipage}[t]{60mm}
\includegraphics[width=60mm]{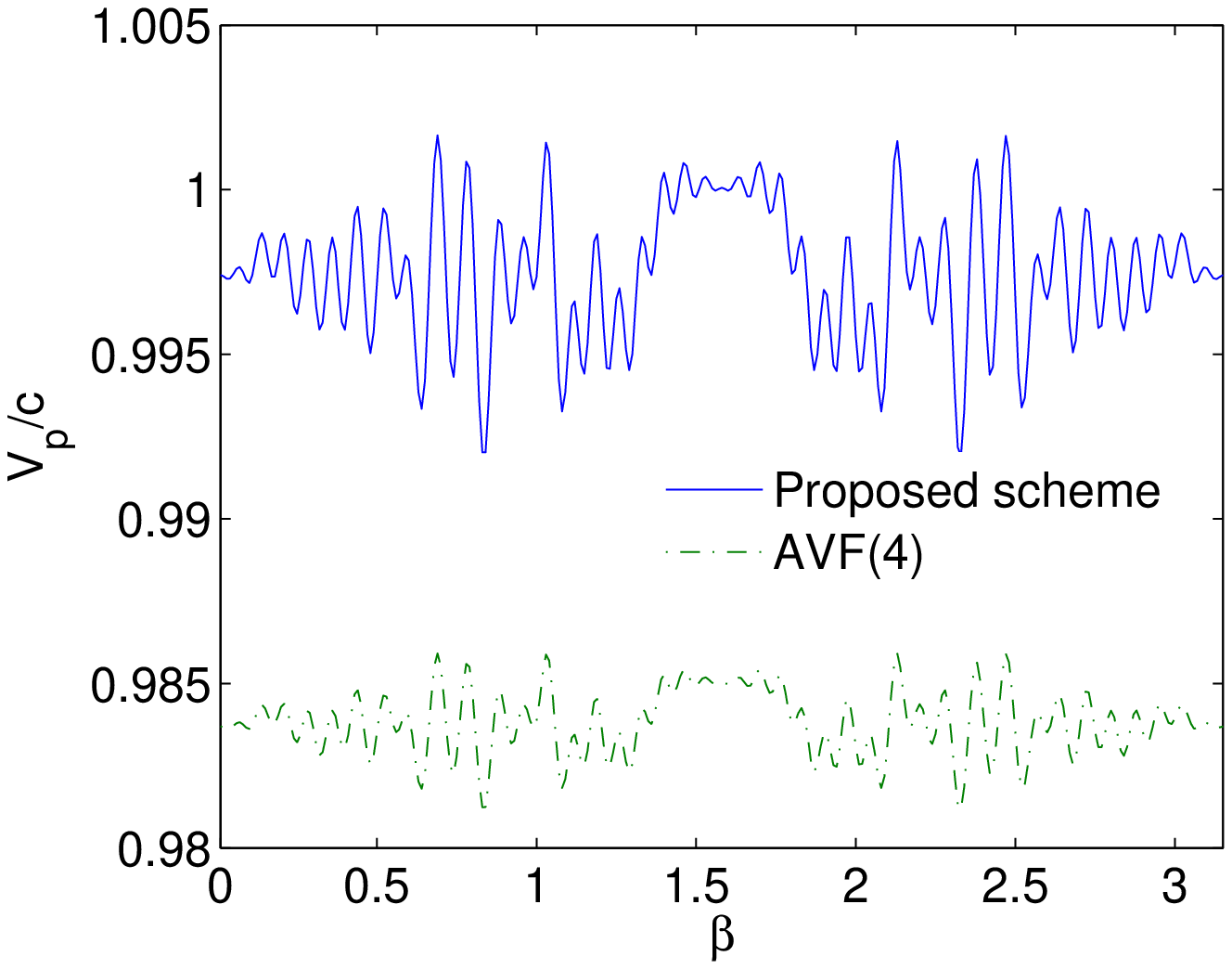}
{\footnotesize (a)}
\end{minipage}\
\begin{minipage}[t]{60mm}
\includegraphics[width=60mm]{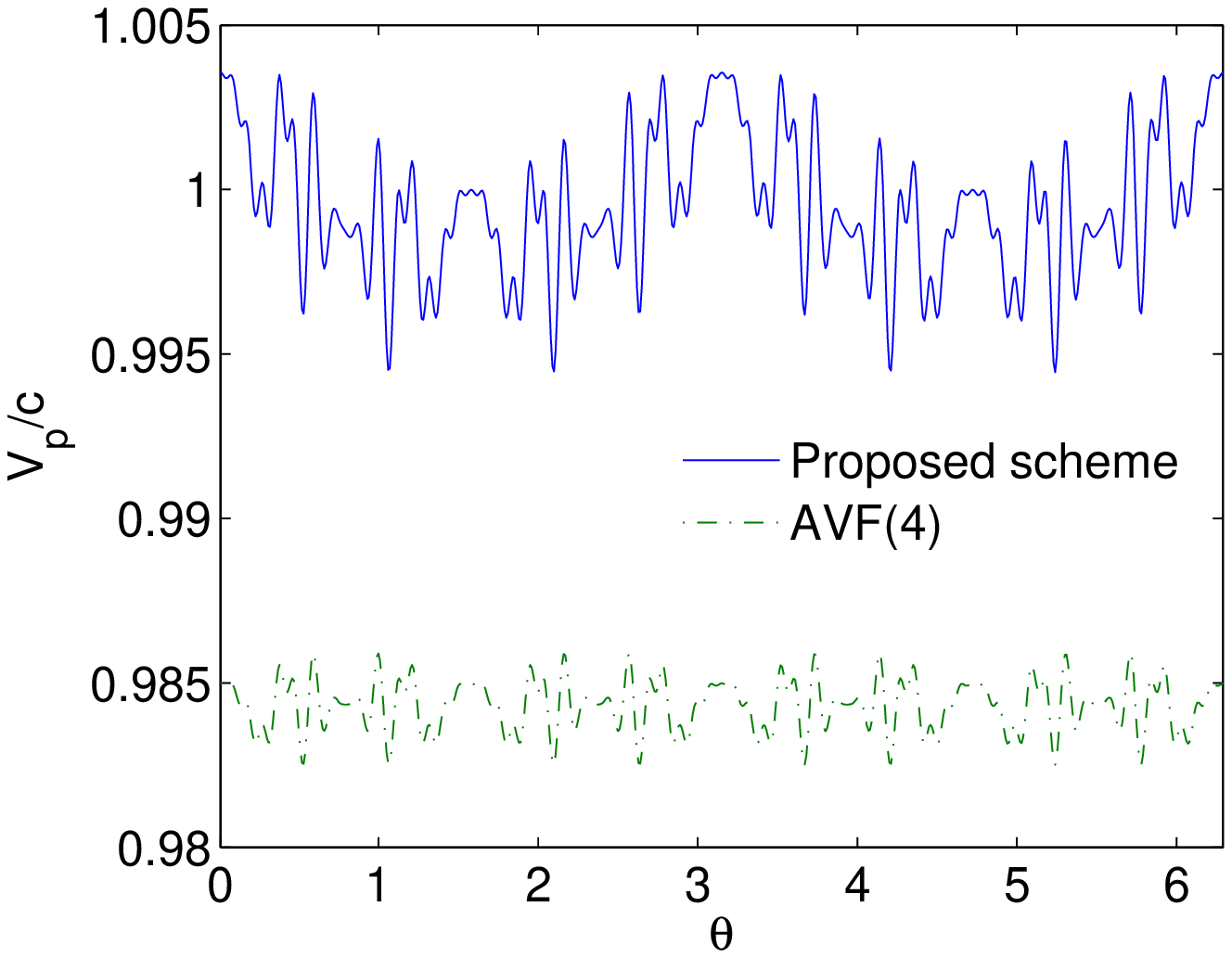}
{\footnotesize (b)}
\end{minipage}\\
\caption{ (a) Normalized numerical phase velocities versus the propagation angles
$\phi$ with $S=1$, $N_{\lambda}=5$ and $\theta=3\pi/8$.
(b) Normalized numerical phase velocities versus the propagation angles with $\theta$ with $S=1$, $N_{\lambda}=5$ and $\phi=\pi/4$.}\label{fig_615}
\end{figure}

First, we will focus on the normalized numerical phase velocity by comparing our scheme (\ref{SO-ECS1}) with the AVF(4) scheme \cite{CWG16}.
Figs. \ref{fig_614} and \ref{fig_615} show the comparisons of normalized numerical phase velocities
versus the Courant-Friedrich-Levy (CFL) number $S$, the number of points per
wavelength $N_{\lambda}$ and the propagation angles $\phi$ and $\theta$.
As illustrated in Figs. \ref{fig_614} and \ref{fig_615}, the proposed scheme (\ref{SO-ECS1}) is
the best one whose normalized numerical phase velocity is the
closest to the analytic solution one, which implies that the
proposed method has lower numerical phase error than the AVF(4) scheme.


\subsubsection{Some analysis on the numerical dispersion relation}
\begin{figure}[H]
\centering\begin{minipage}[t]{60mm}
\includegraphics[width=60mm]{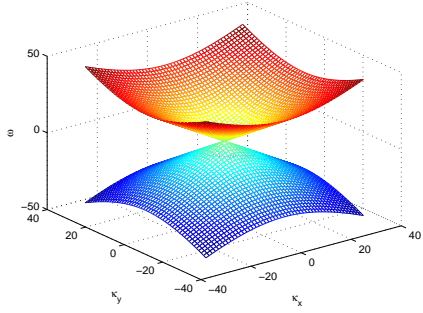}
{\footnotesize (a) Exact dispersion (\ref{exact})}
\end{minipage}\quad\quad
\begin{minipage}[t]{60mm}
\includegraphics[width=60mm]{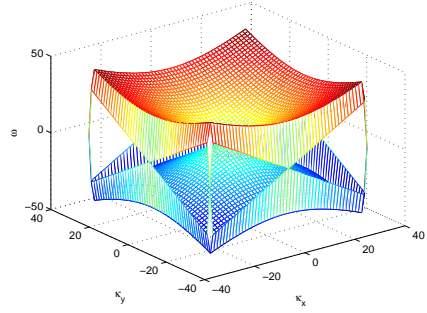}
{\footnotesize (b) Scheme (\ref{SO-ECS1})}
\end{minipage}
\\
\caption{The dispersion relation figures on ($\kappa_{x},\kappa_y$)
with $\tau=0.01,\ h=0.1$ and $N=150$ for the Maxwell's equations.}\label{fig611}
\end{figure}

\begin{figure}[H]
\centering\begin{minipage}[t]{60mm}
\includegraphics[width=60mm]{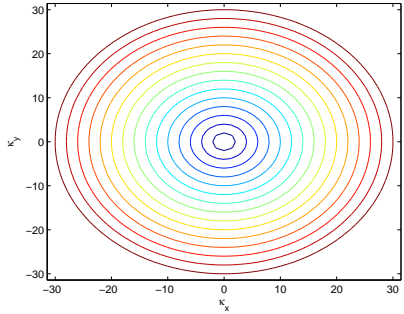}
{\footnotesize (a) Exact dispersion (\ref{exact})}
\end{minipage}\quad\quad
\begin{minipage}[t]{60mm}
\includegraphics[width=60mm]{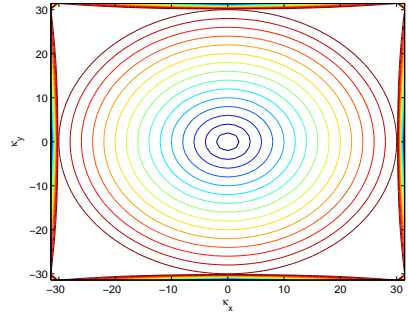}
{\footnotesize (b) Scheme (\ref{SO-ECS1})}
\end{minipage}
\\
\caption{The contour plots on ($\kappa_{x},\kappa_y$)-plane with $\tau=0.01,\ h=0.1$ and $N=150$ for the Maxwell's equations.
 }\label{fig612}
\end{figure}

In this subsection, we make some analysis on the numerical dispersion relation of the proposed scheme (\ref{SO-ECS1}).
In Fig. \ref{fig611}, we first plot the dispersion relation for the frequency $\omega$ as a function of the wave number vector ($\kappa_{x}, \kappa_{y}$).
From Fig. \ref{fig611} (b), one can see that, similar to the conventional FDTD methods \cite{CWS13,ST11}, the proposed method (\ref{SO-ECS1}) gives a
numerical dispersion surface $\omega$ with extra solution
branches which correspond to the existence of the nonphysical parasitic
waves in the numerical solution \cite{AM04}.
But, in contrast to the FDTD methods,
such extra solution branches only occur with large wave numbers.
This fact can be verified in following analysis again.
In order to see the relation clearly, the contour plots of $\omega$ are displayed in Fig \ref{fig612}.
As illustrated in Fig \ref{fig612} (b), we can observe that, except for the large wave numbers $\kappa_{w},w=x,y$,
the numerical contours are circular and the distances and shapes of the numerical contours are
almost consistent with the exact one.
This conforms the fact in Fig. \ref{fig611} (b) and implies
that the numerical propagation speed $|v_{g}|$ is direction-independent.
In addition, we can infer
that the grid-anisotropy \cite{NT82} of the scheme (\ref{SO-ECS1}) causes
that propagation directions $\alpha$ and $\beta$ (see (\ref{sc1})) are independent of the directions.
Here, due to the symmetric property of the
numerical dispersion relation,
we omit the figures of the dispersion relation on the other planes.

In order to verify the above inference, we investigate the magnitude $|v_{g}|$
of the group velocity (\ref{vg1}) in spherical coordinates (\ref{sc2}).
The group velocity, which characterizes the speed of energy transport in wave packets \cite{WB99},
is fundamental to the understanding of linear waves \cite{WB99}.

\begin{figure}[H]
\centering\begin{minipage}[t]{60mm}
\includegraphics[width=60mm]{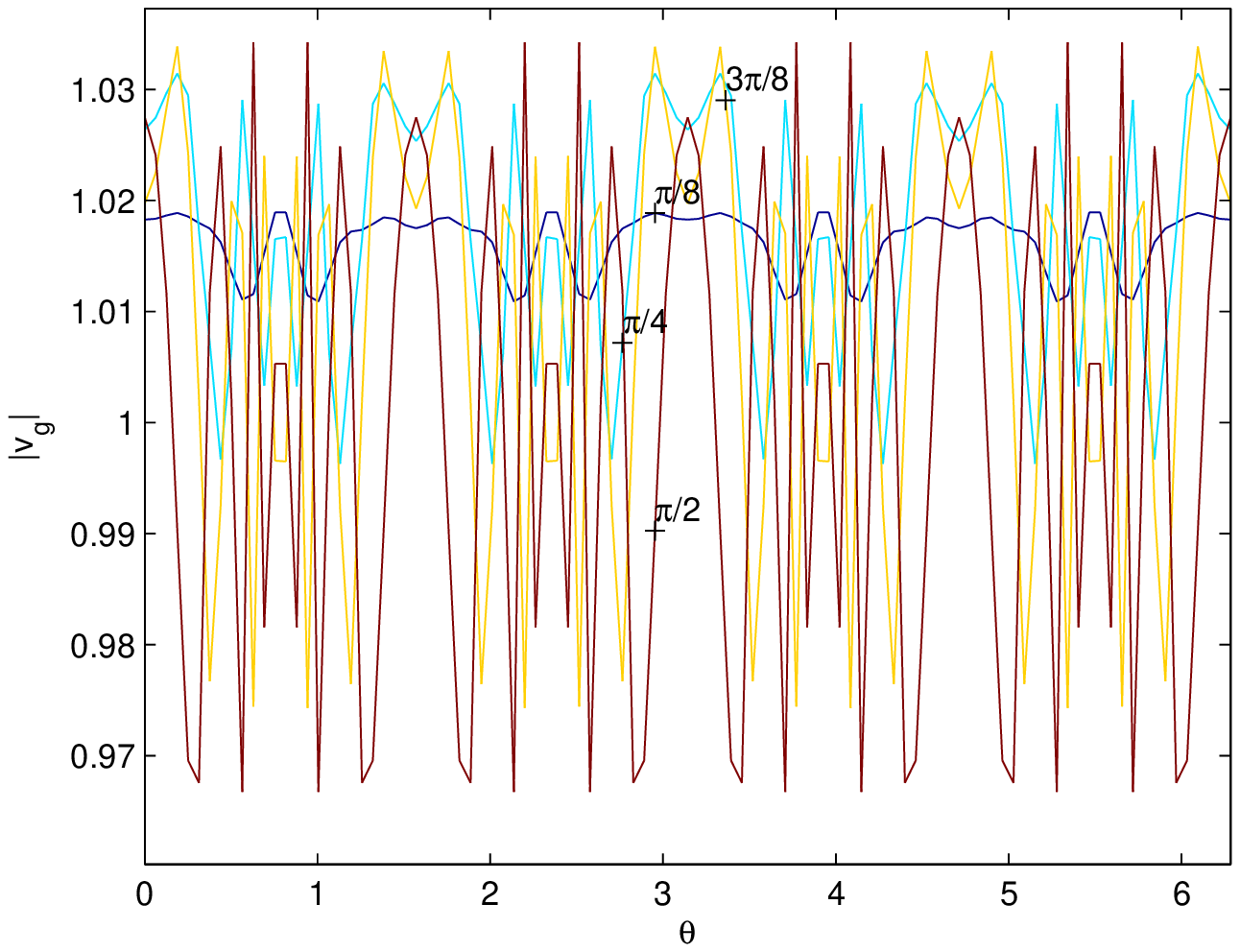}\\
{\footnotesize (a) $|{\bm \kappa}|=2.5\pi$}
\end{minipage}
\begin{minipage}[t]{60mm}
\includegraphics[width=60mm]{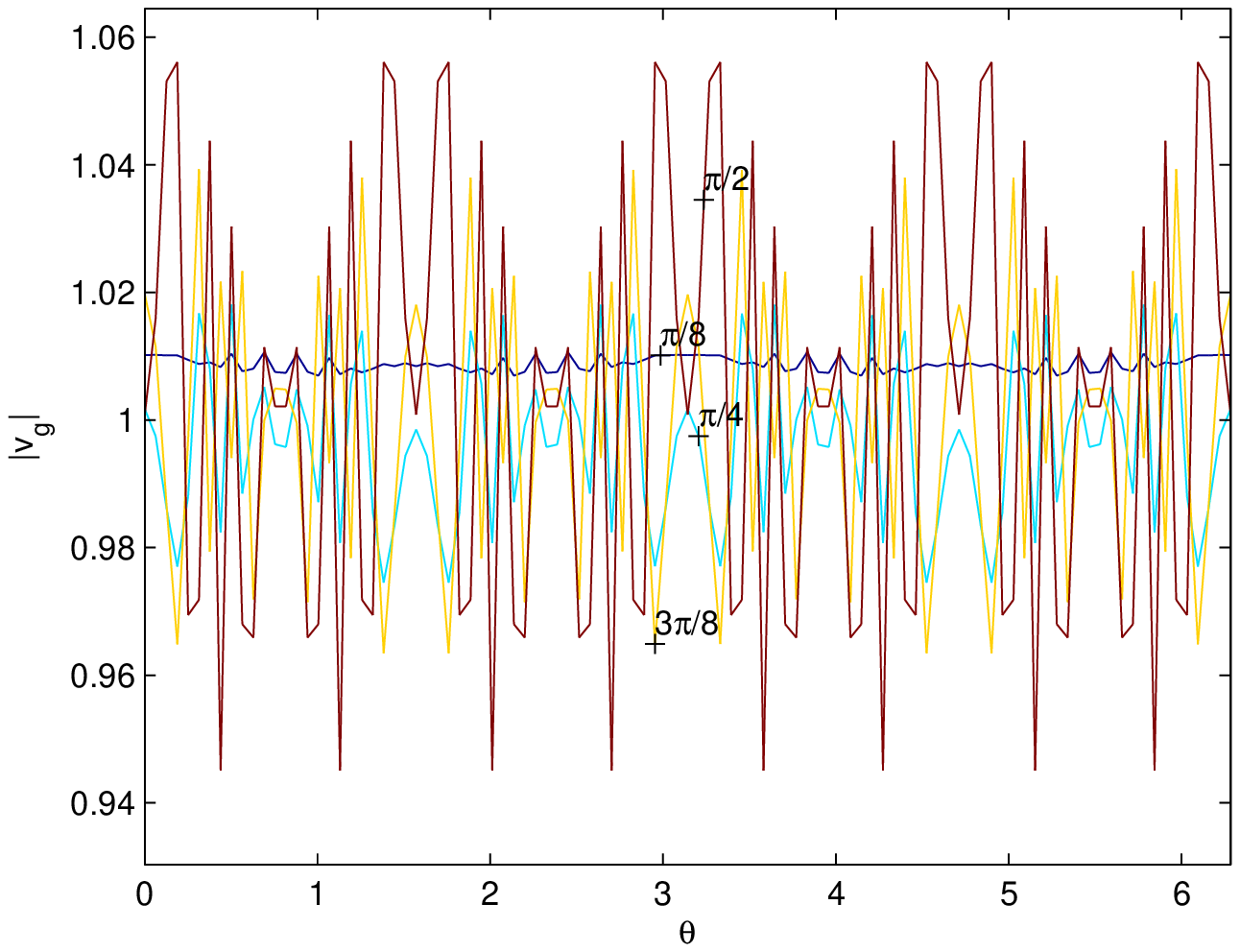}
{\footnotesize (b) $|{\bm \kappa}|=5\pi$}
\end{minipage}\\
\begin{minipage}[t]{60mm}
\includegraphics[width=60mm]{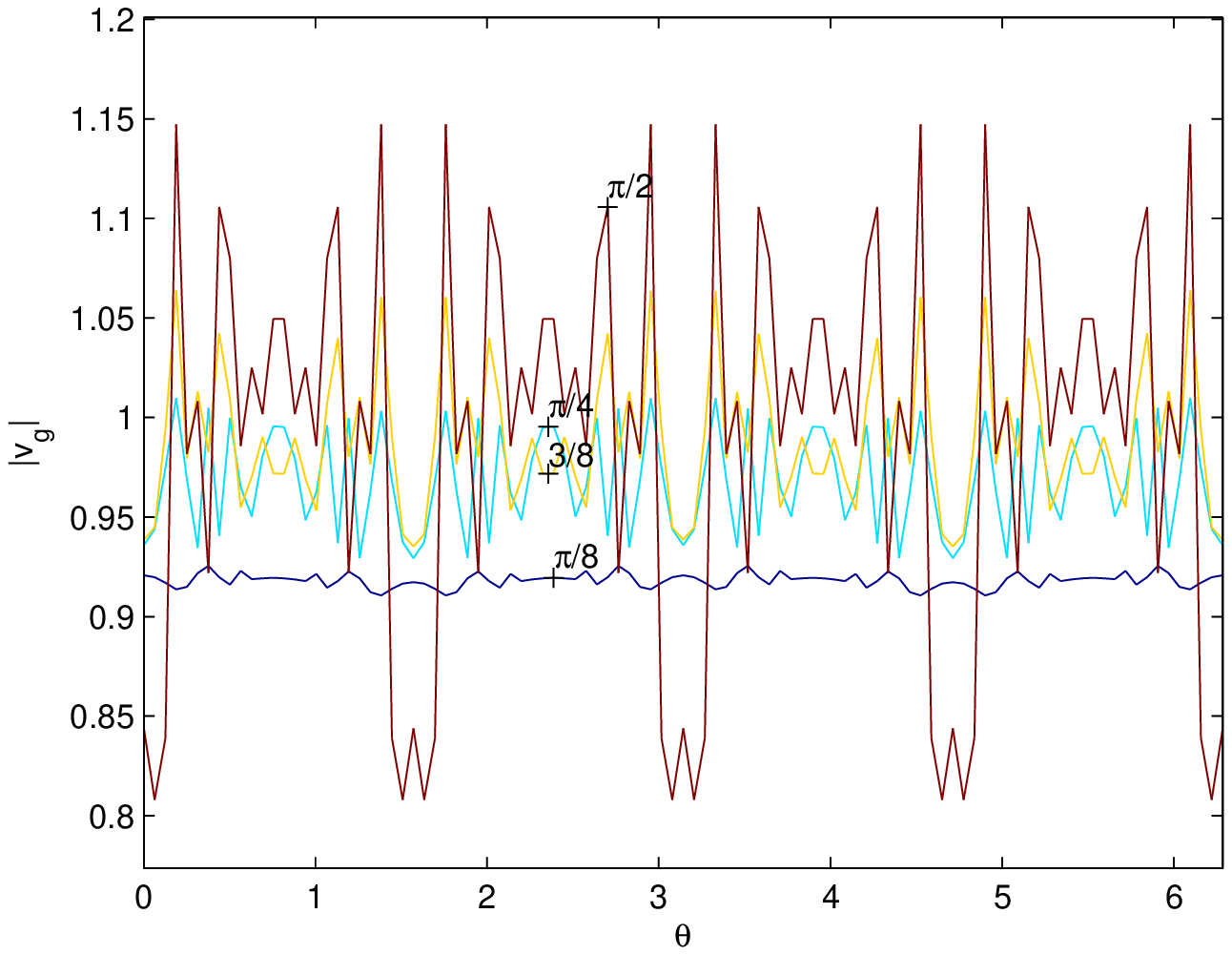}
{\footnotesize (c) $|{\bm \kappa}|=7.5\pi$}
\end{minipage}\\
\caption{Numerical group velocities at different $|{\bm \kappa}|$ and $\phi$ with $\tau=0.01,\ h=0.1$ and $N=150$ for the Maxwell's equations.}\label{fig613}
\end{figure}

\begin{figure}[H]
\centering\begin{minipage}[t]{60mm}
\includegraphics[width=60mm]{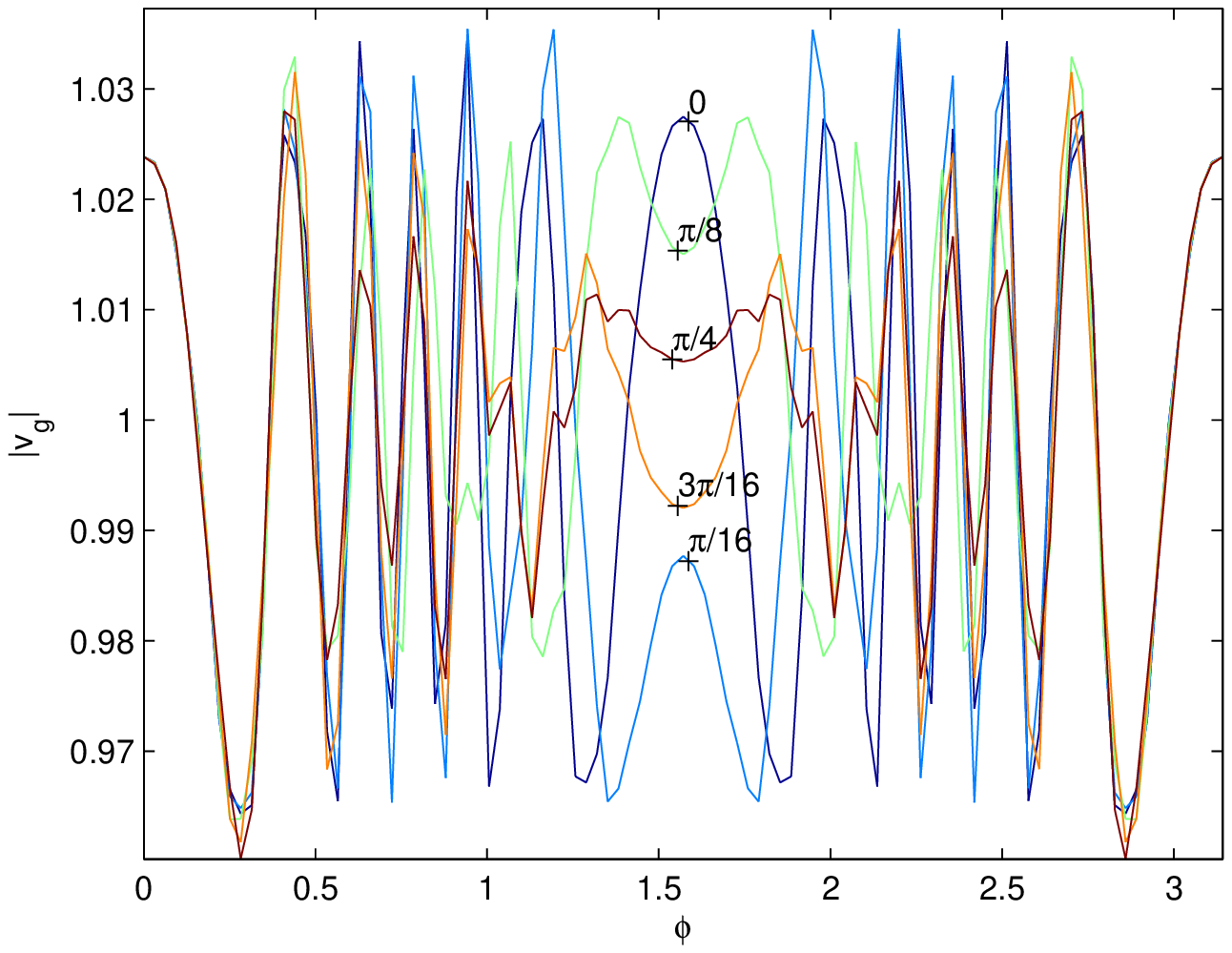}
{\footnotesize (a) $|{\bm \kappa}|=2.5\pi$}
\end{minipage}\
\begin{minipage}[t]{60mm}
\includegraphics[width=60mm]{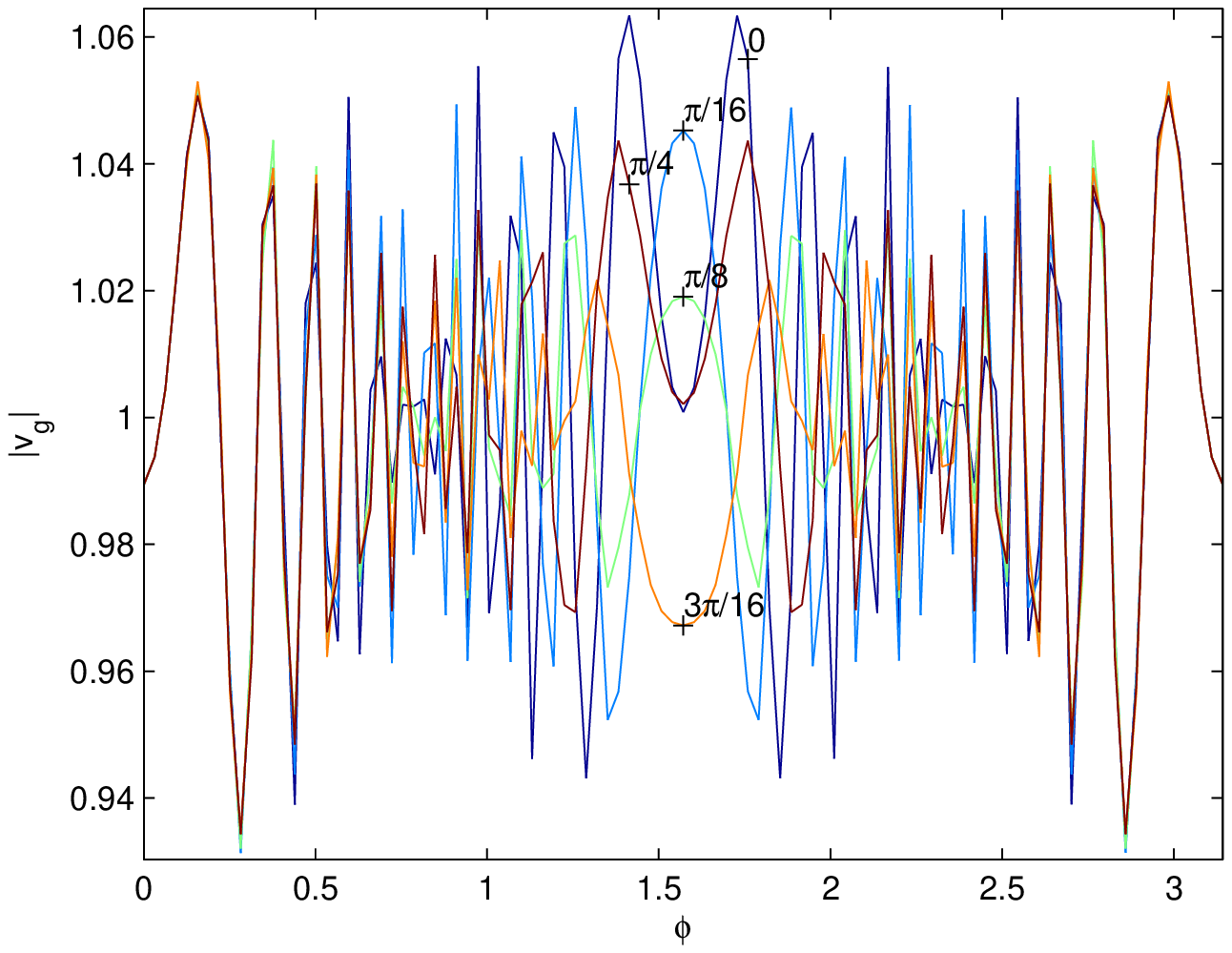}
{\footnotesize (b) $|{\bm \kappa}|=5\pi$}
\end{minipage}\\
\begin{minipage}[t]{60mm}
\includegraphics[width=60mm]{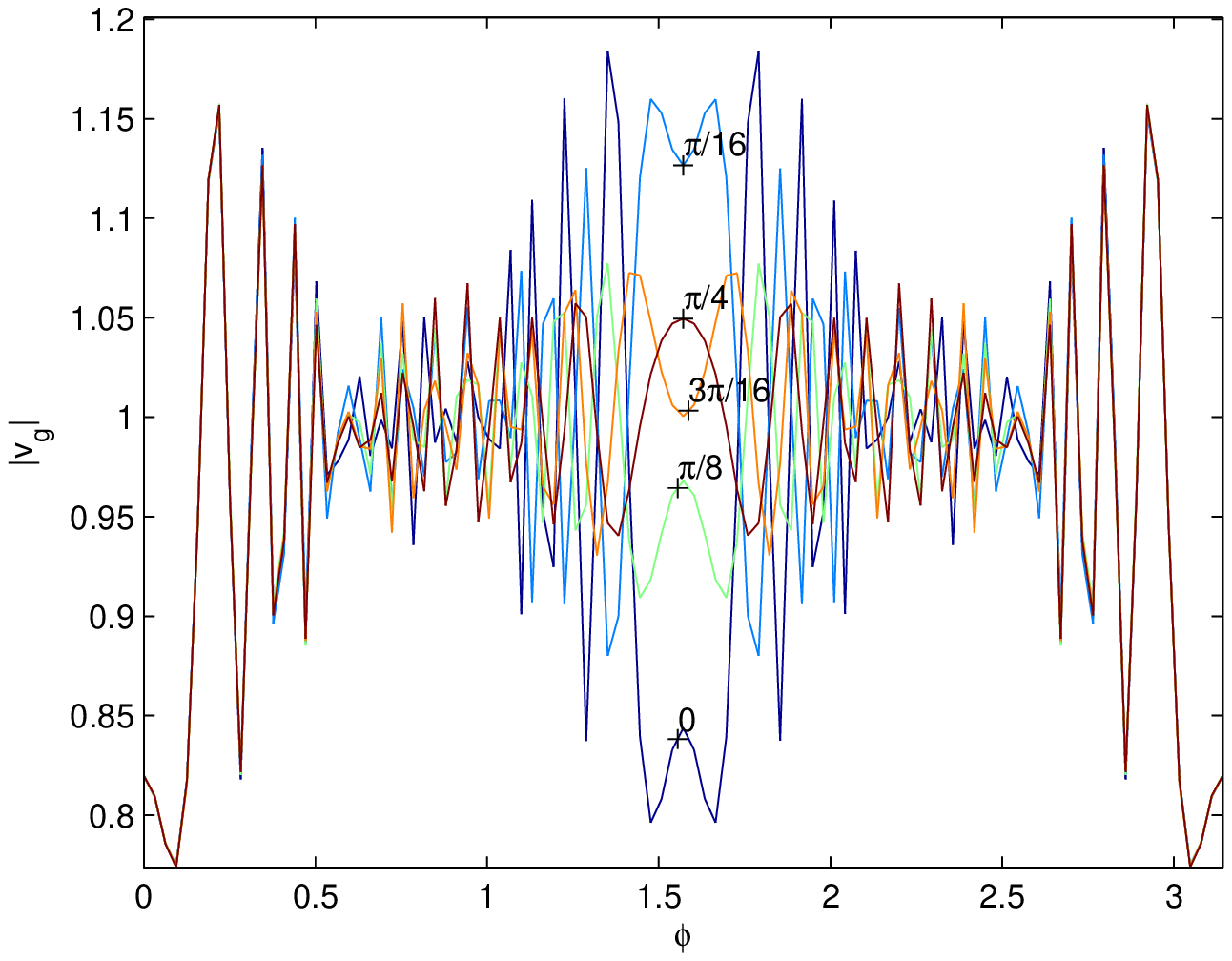}
{\footnotesize (c) $|{\bm \kappa}|=7.5\pi$}
\end{minipage}\\
\caption{Numerical group velocities at different $|{\bm \kappa}|$ and $\theta$ with $\tau=0.01,\ h=0.1$ and $N=150$ for the Maxwell's equations.}\label{fig614}
\end{figure}

Fig. \ref{fig613} shows the relation between the numerical group
velocities and the wave number angles $\theta$ at different $|{\bm \kappa}|$ and $\phi$.
The group velocity, which characterizes the speed of energy transport in wave packets \cite{WB99}, is fundamental to the understanding of linear waves \cite{WB99}.
The relation between the numerical group velocities and the wave number angles $\phi$ at different $|{\bm \kappa}|$ and $\theta$
is displayed in Fig. \ref{fig614}.

We can make some observations from the figures:
\begin{itemize}
\item[(1)] For all $|{\bm \kappa}|$, numerical group velocity occurs oscillation phenomena with exact value one.
But it's worth noting that, in Fig. \ref{fig613},
for fixed $\phi=\frac{\pi}{8}$, when $|{\bm \kappa}|=2.5\pi$ or $5\pi,\ |v_{g}|$ is greater
than the exact value one. However, for $|{\bm \kappa}|=7.5\pi$, $|v_{g}|$ is less than one.
\item[(2)]The maximum value of the numerical group velocity increases as $|{\bm \kappa}|$ increases, which implies the magnitude of $v_{g}$
depends on the vector wave number $|{\bm \kappa}|$. Therefore, the proposed scheme (\ref{SO-ECS1}) is dispersive. That is, a wave-packet
of this scheme with different $|{\bm \kappa}|$ will spread out.
\item[(3)] In Fig. \ref{fig613}, it is clear to see that the $|v_{g}|$ is symmetric with respect to $\theta=\pi$. If dividing the domain of $\theta$ into four equal
parts $[0,\frac{\pi}{2}]$, $[\frac{\pi}{2},\pi]$, $[\pi,\frac{3\pi}{2}]$, and $[\frac{3\pi}{2},2\pi]$, respectively, we can find that, whatever $|{\bm \kappa}|$ and $\phi$
chosen, $v_{g}$ is symmetric with respect to $\theta=\frac{\pi}{4},\frac{3\pi}{4},\frac{5\pi}{4},\frac{7\pi}{4}$ in each part of the domain.
Furthermore, from Fig. \ref{fig614}, we can observe that, whatever $|{\bm \kappa}|$ chooses, the $|v_{g}|$ is also symmetric with respect to $\phi=\frac{\pi}{2}$
which means the $(x,y)$-plane.
\end{itemize}

We then calculate the propagation angles of the group velocity $\alpha$ and $\beta$.
Note that the exact dispersion relation is as follows
\begin{align}
&v_{g}=|v_{g}|\left(\sin\alpha\cos\beta,\sin\alpha\sin\beta,\cos\alpha\right)=\left(\frac{\partial \omega}{\partial \kappa_{x}},
\frac{\partial \omega}{\partial \kappa_{y}},\frac{\partial \omega}{\partial \kappa_{z}}\right)\nonumber\\
&~~~=\left(\frac{\kappa_{x}}{\omega},
\frac{\kappa_{y}}{\omega},\frac{\kappa_{z}}{\omega}\right)
=|\kappa|\left(\frac{\sin\phi\cos\theta}{\omega},\frac{\sin\phi\sin\theta}{\omega},\frac{\cos\phi}{\omega}\right),
\end{align}
we can obtain $\alpha=\phi$ and $\beta=\theta$. Similarly, by virtue of the numerical dispersion relation, we can derive
\begin{align}
&\alpha=\arctan\Bigg(\sqrt{\frac{(\frac{\partial \omega}{\partial\kappa_{x}})^{2}
+(\frac{\partial \omega}{\partial\kappa_{y}})^{2}}{(\frac{\partial\omega}{\partial \kappa_{z}})^{2}}}\Bigg),\\
&\beta=\arctan\Bigg(\frac{\frac{\partial\omega}{\partial\kappa_{y}}}{\frac{\partial\omega}{\partial\kappa_{x}}}\Bigg),
\end{align}
where $\frac{\partial\omega}{\partial\kappa_{x}}, \frac{\partial\omega}{\partial\kappa_{y}}$ and $\frac{\partial\omega}{\partial\kappa_{x}}$ are defined
as above.

\begin{figure}[H]
\centering\begin{minipage}[t]{60mm}
\includegraphics[width=60mm]{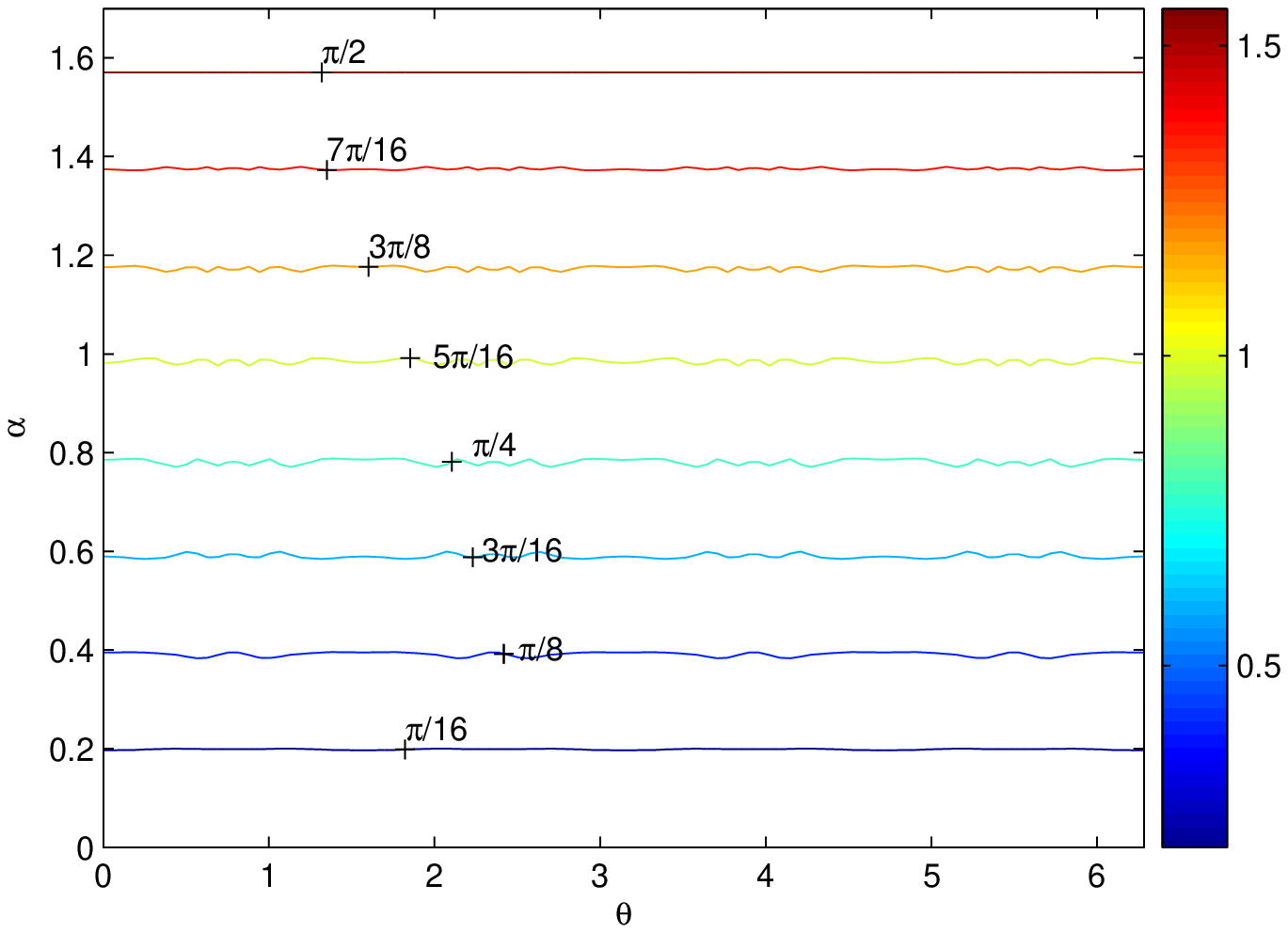}
{\footnotesize (a) $|{\bm \kappa}|=2.5\pi$}
\end{minipage}\begin{minipage}[t]{60mm}
\includegraphics[width=60mm]{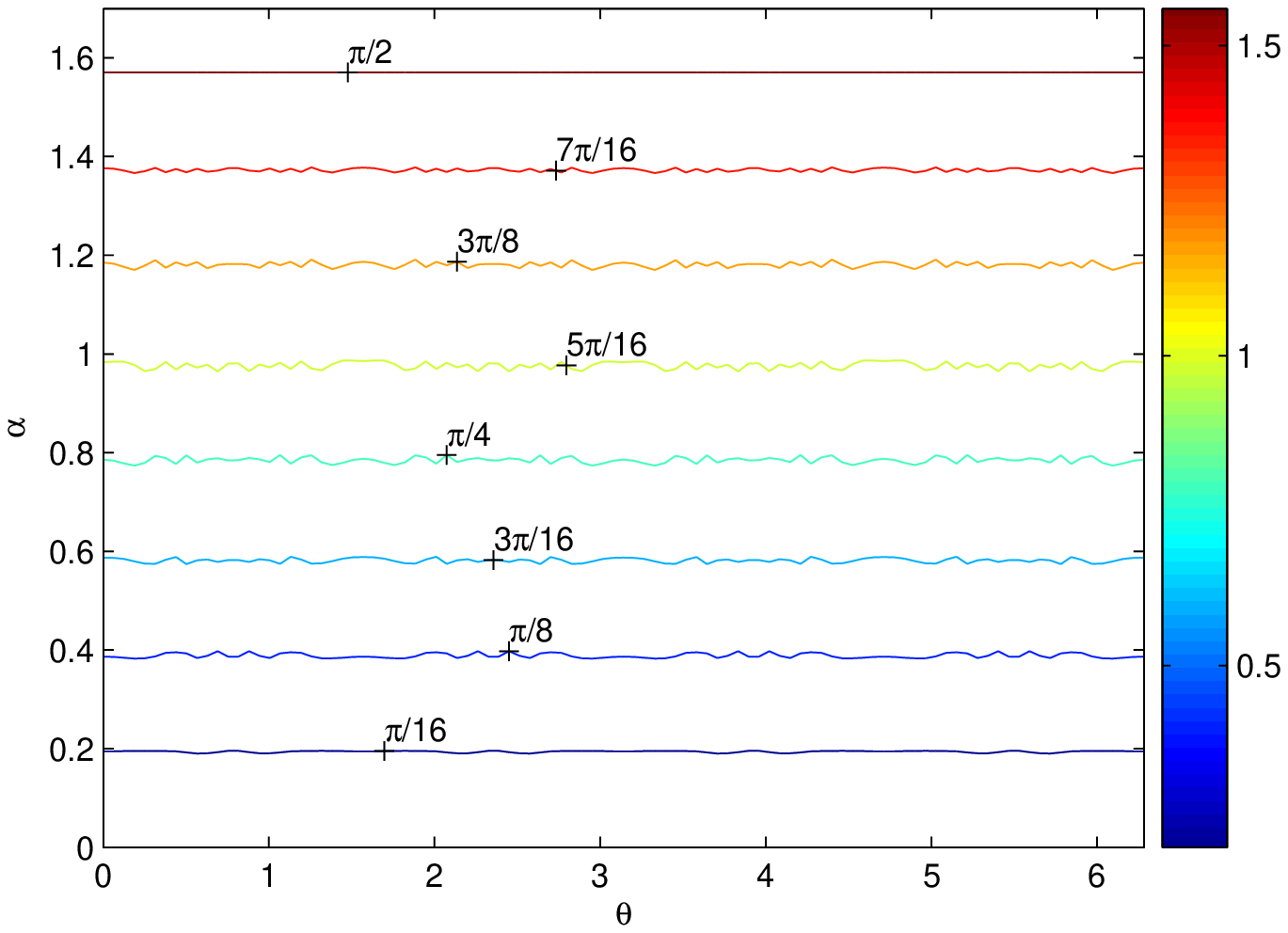}
{\footnotesize (b) $|{\bm \kappa}|=5\pi$}
\end{minipage}\\
\begin{minipage}[t]{60mm}
\includegraphics[width=60mm]{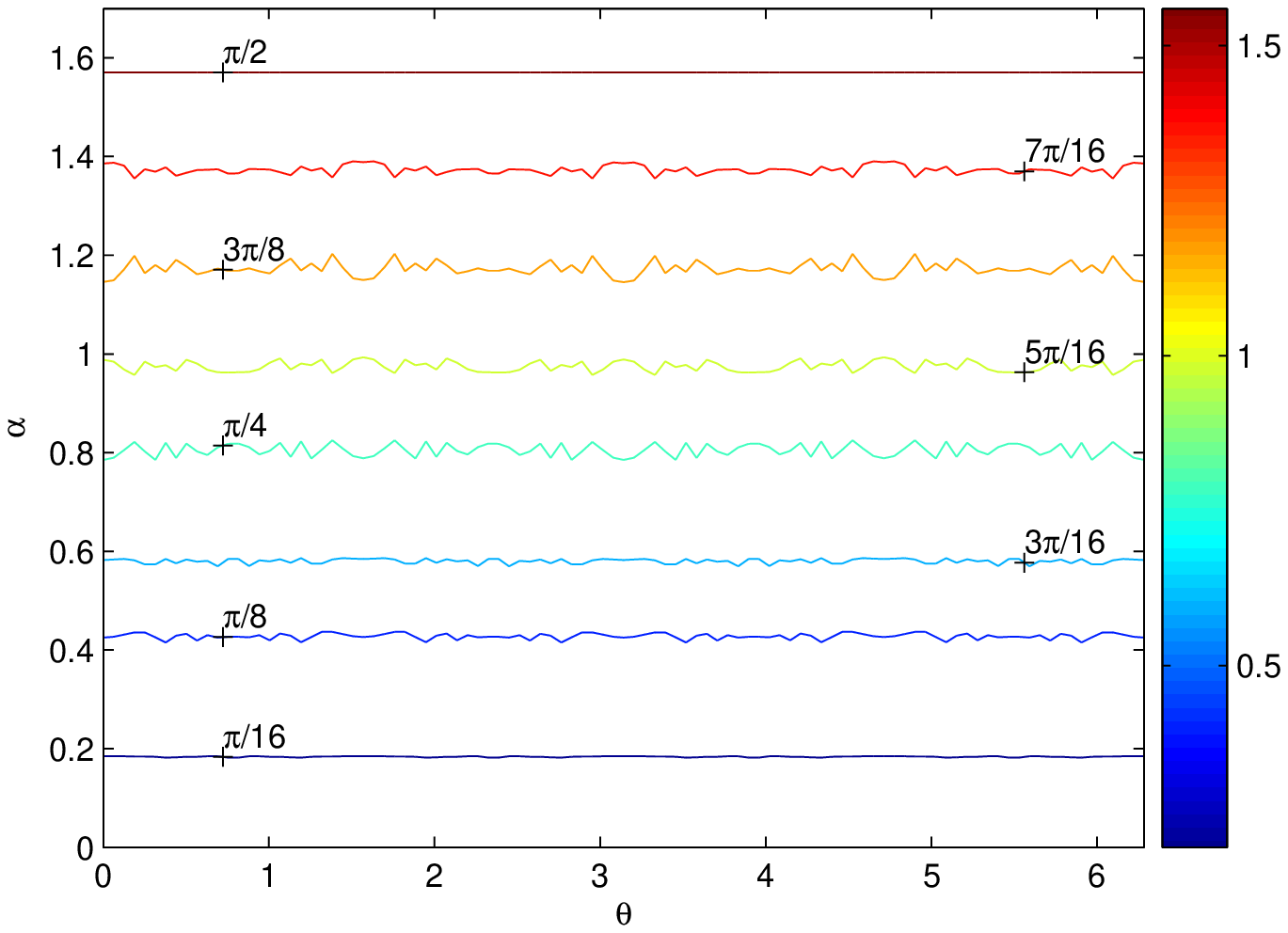}
{\footnotesize (c) $|{\bm \kappa}|=7.5\pi$}
\end{minipage}\\
\caption{Wave propagation angle $\alpha$ versus the wave number angle $\theta$ and $\phi$ at different $|{\bm \kappa}|$ with $\tau=0.01,\ h=0.1$ and $N=150$ for the Maxwell's equations.}\label{fig615}
\end{figure}

\begin{figure}[H]
\centering\begin{minipage}[t]{60mm}
\includegraphics[width=60mm]{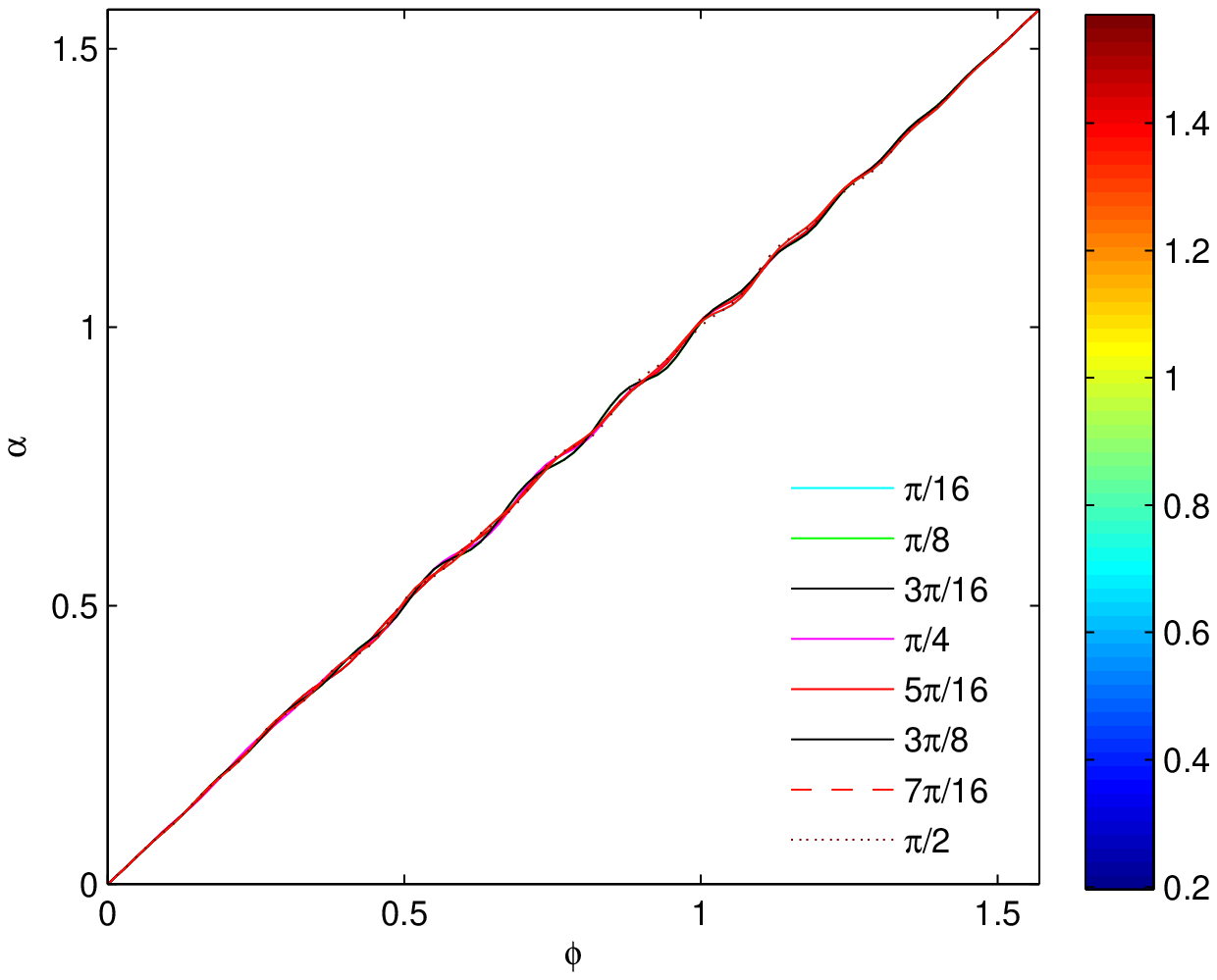}
{\footnotesize (a) $|{\bm \kappa}|=2.5\pi$}
\end{minipage}\
\begin{minipage}[t]{60mm}
\includegraphics[width=60mm]{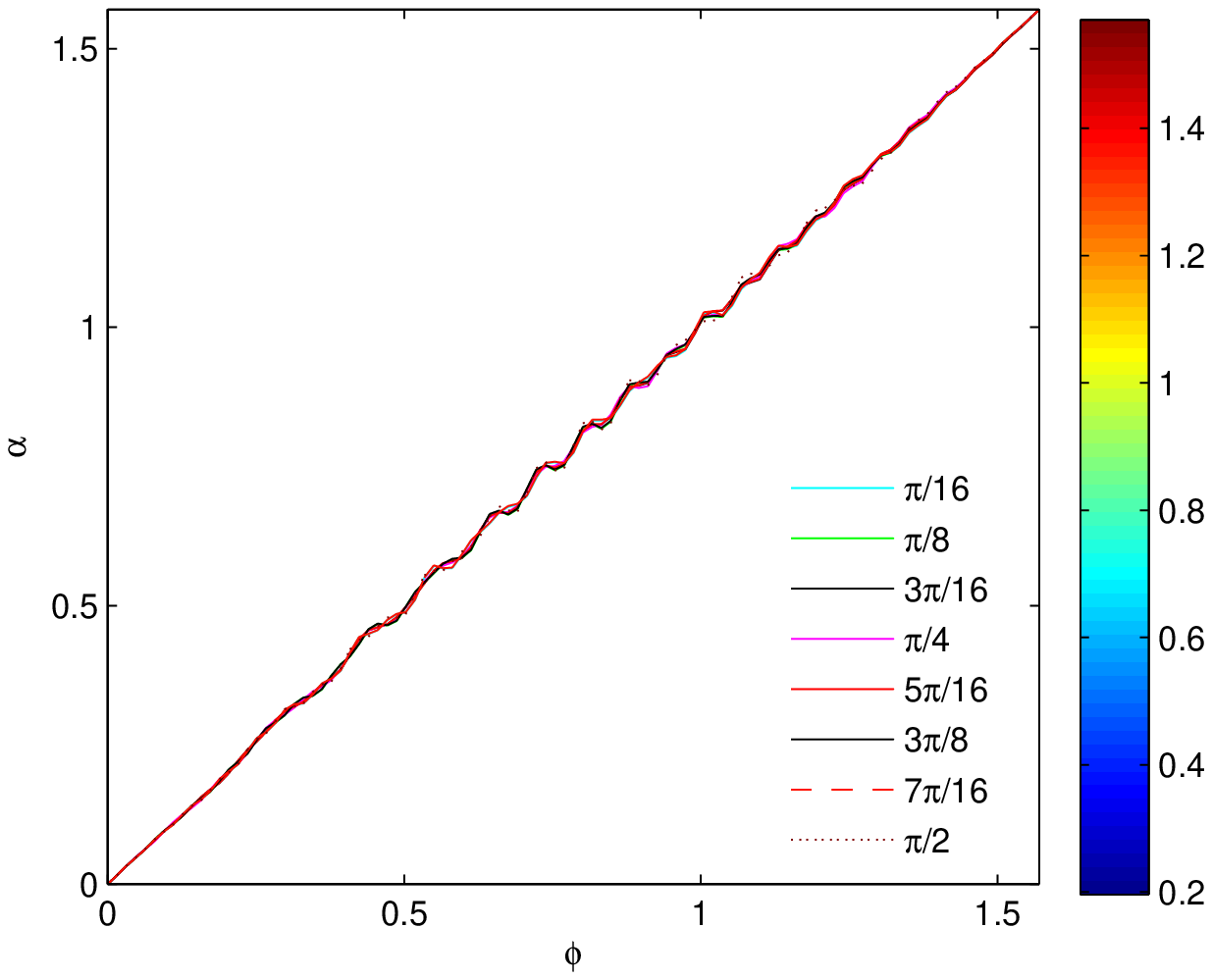}
{\footnotesize (b) $|{\bm \kappa}|=5\pi$}
\end{minipage}\\
\begin{minipage}[t]{60mm}
\includegraphics[width=60mm]{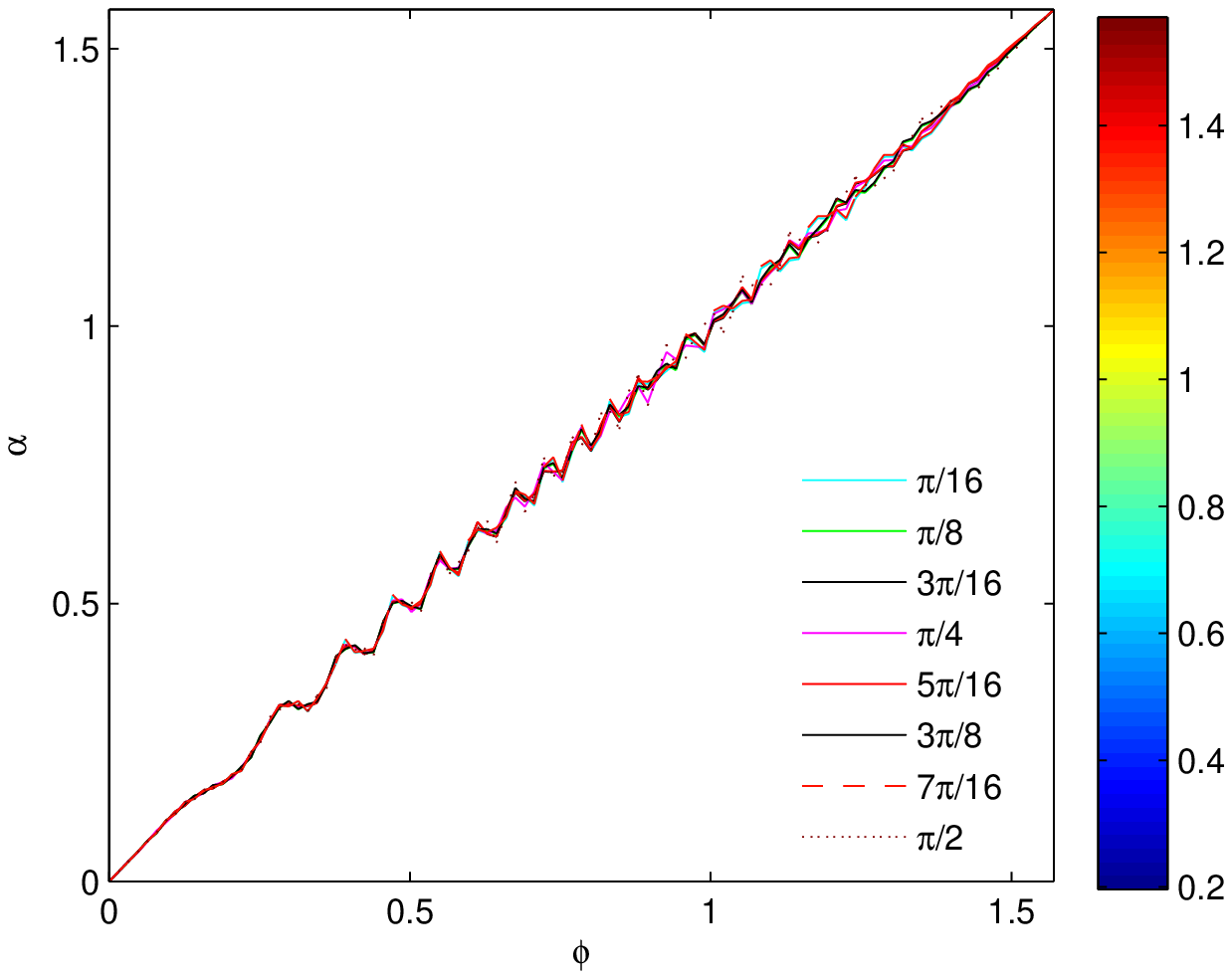}
{\footnotesize (c) $|{\bm \kappa}|=7.5\pi$}
\end{minipage}\\
\caption{Wave propagation angle $\alpha$ versus the wave number angle $\phi$ at different $|{\bm \kappa}|$ with $\tau=0.01,\ h=0.1$ and $N=150$ for the Maxwell's equations.}\label{fig616}
\end{figure}

\begin{figure}[H]
\centering\begin{minipage}[t]{60mm}
\includegraphics[width=60mm]{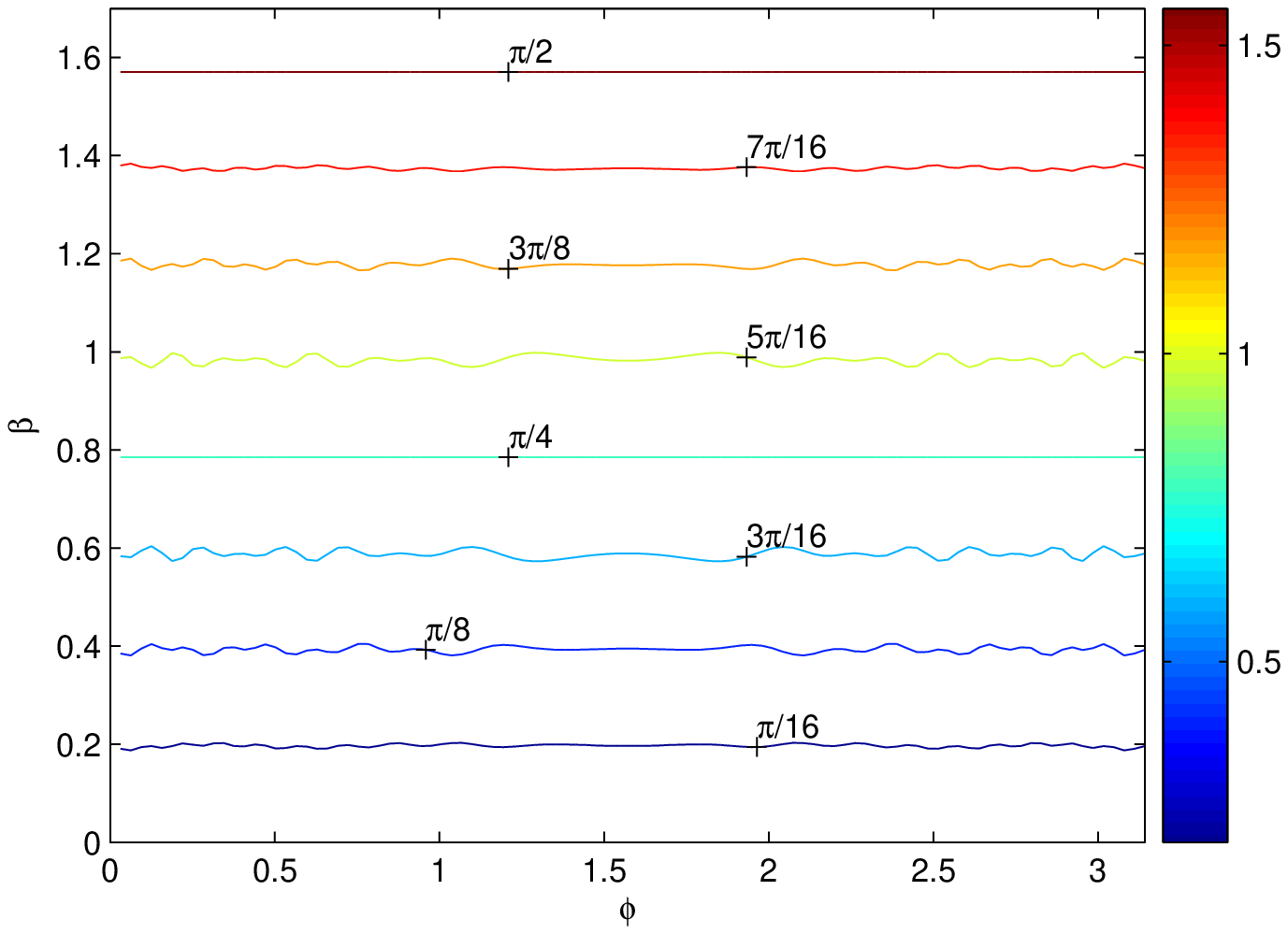}
{\footnotesize (a) $|{\bm \kappa}|=2.5\pi$ }
\end{minipage}\begin{minipage}[t]{55mm}
\includegraphics[width=60mm]{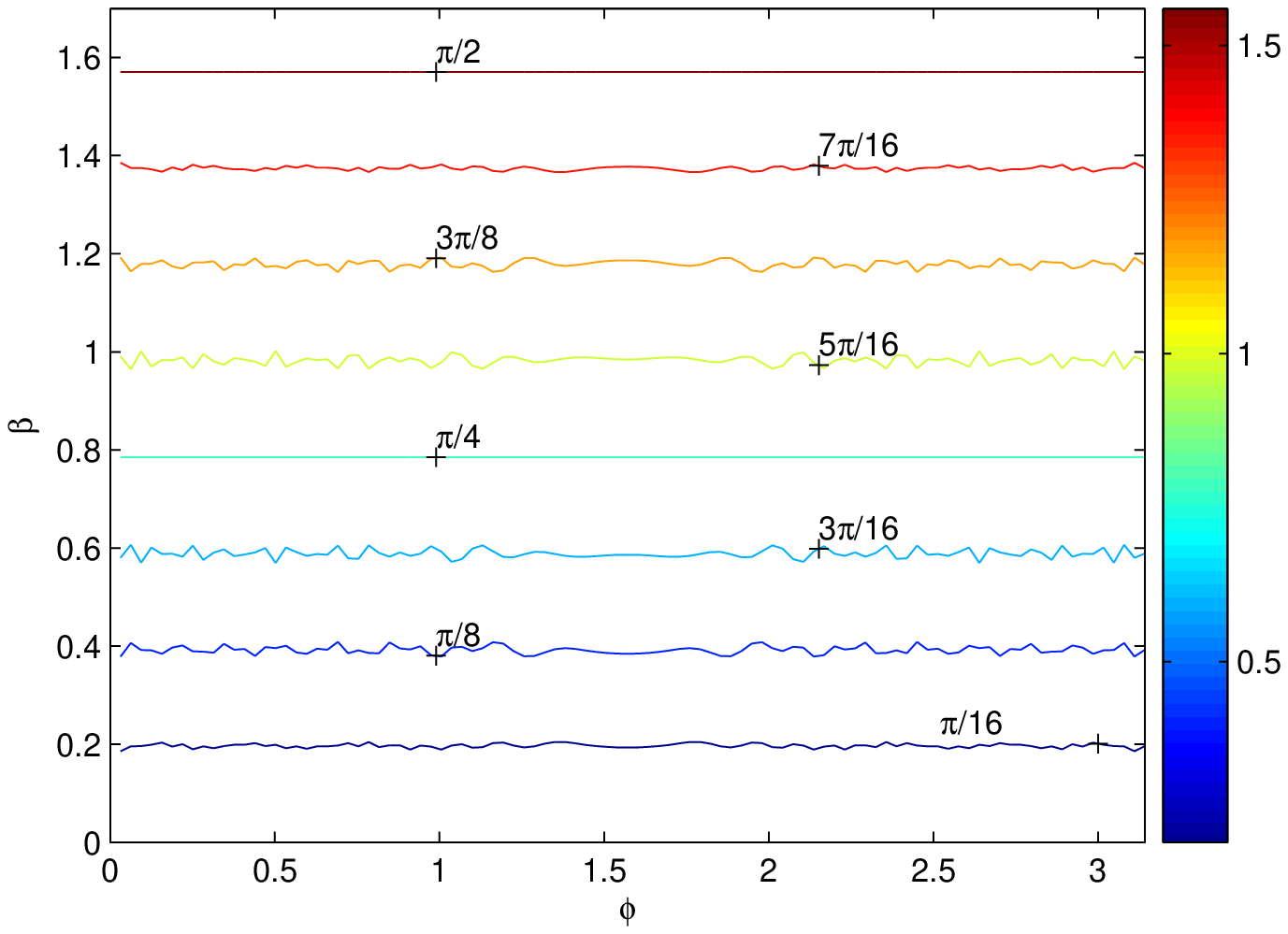}
{\footnotesize (b) $|{\bm \kappa}|=5\pi$}
\end{minipage}\\
\begin{minipage}[t]{60mm}
\includegraphics[width=60mm]{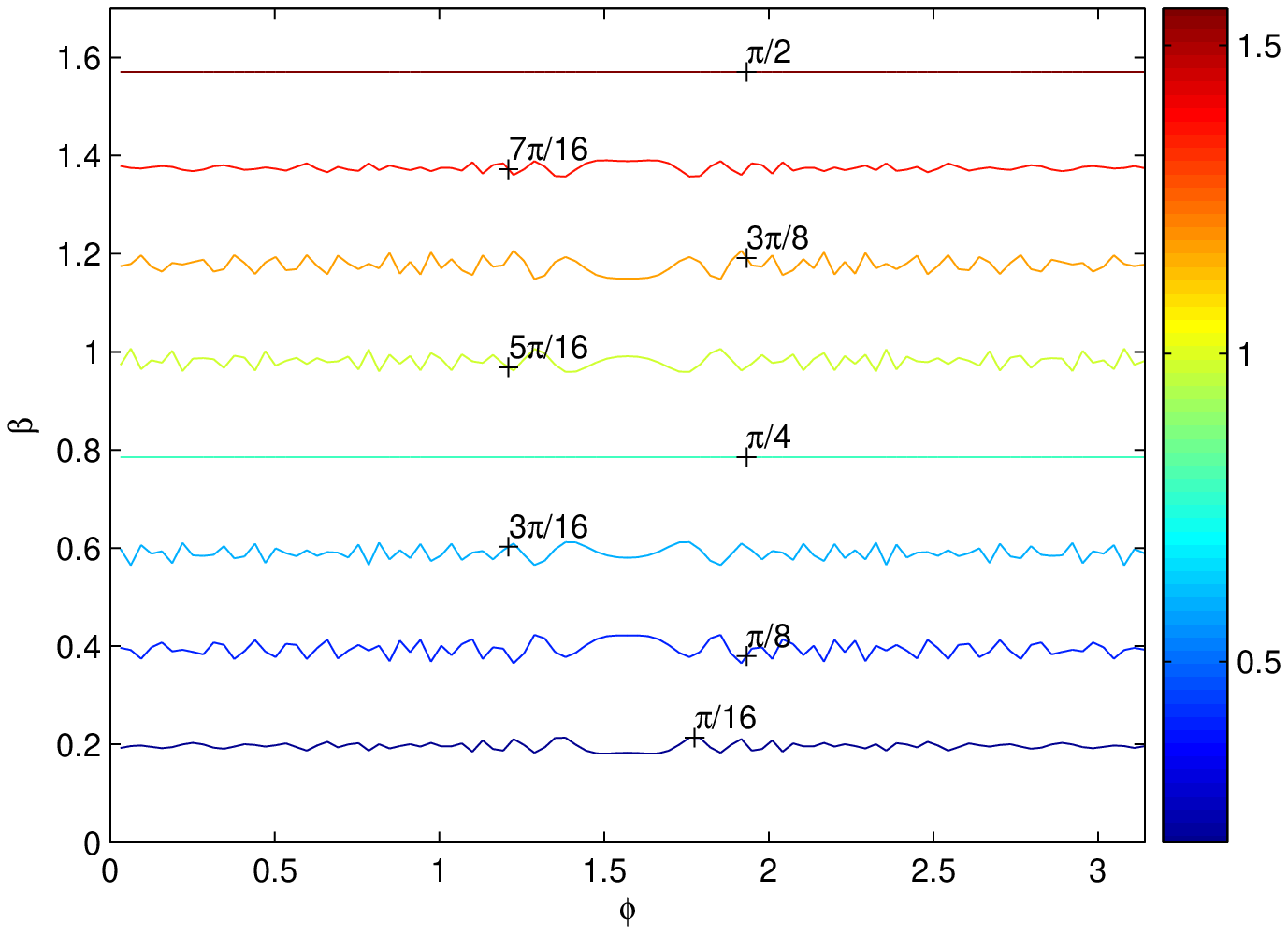}
{\footnotesize (c) $|{\bm \kappa}|=7.5\pi$}
\end{minipage}\\
\caption{Wave propagation angle $\beta$ versus the wave number angle $\phi$ and $\theta$ at different $|{\bm \kappa}|$ with $\tau=0.01,\ h=0.1$ and $N=150$ for the Maxwell's equations.}\label{fig617}
\end{figure}

\begin{figure}[H]
\centering\begin{minipage}[t]{60mm}
\includegraphics[width=60mm]{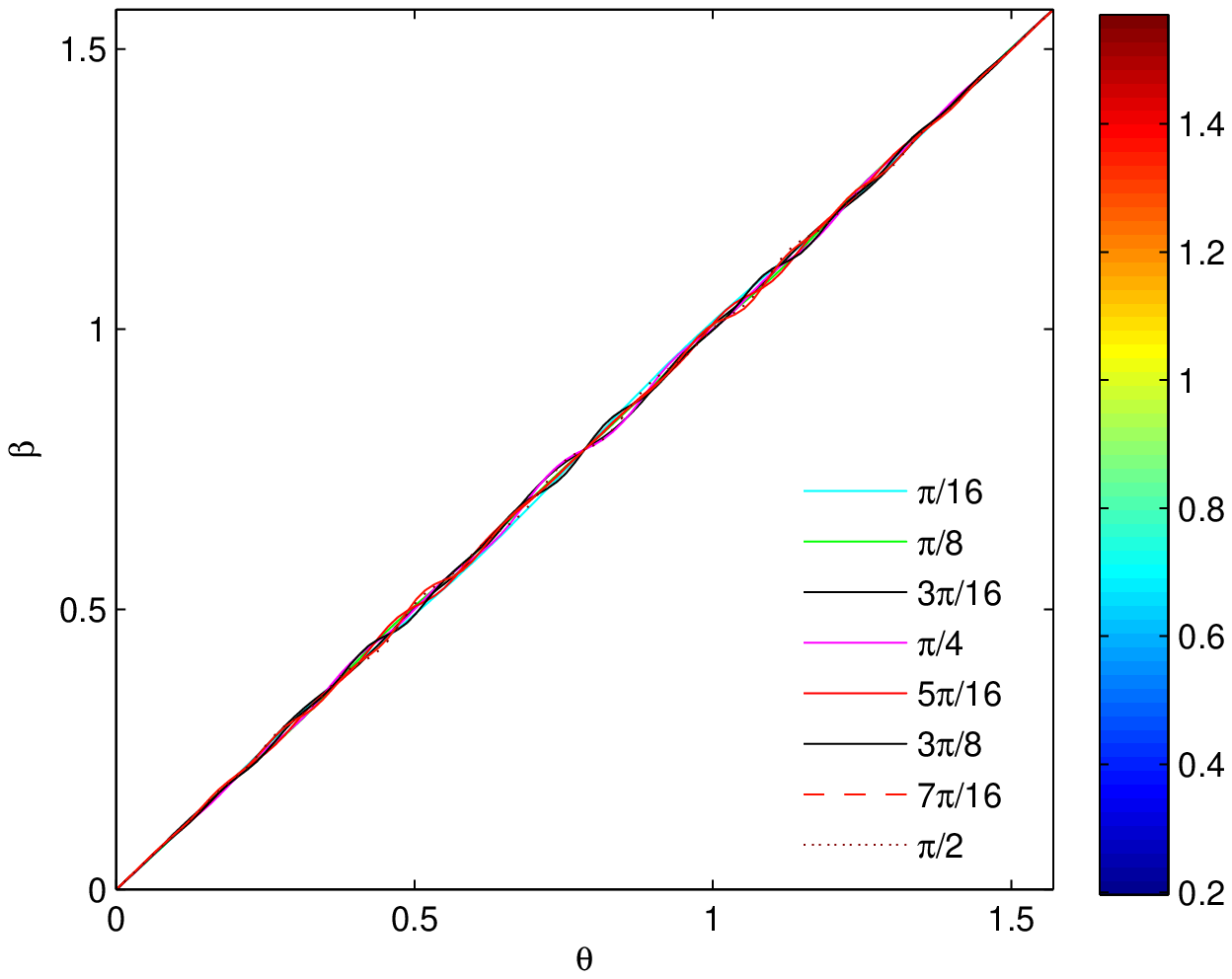}
{\footnotesize (a) $|{\bm \kappa}|=2.5\pi$ }
\end{minipage}\
\begin{minipage}[t]{60mm}
\includegraphics[width=60mm]{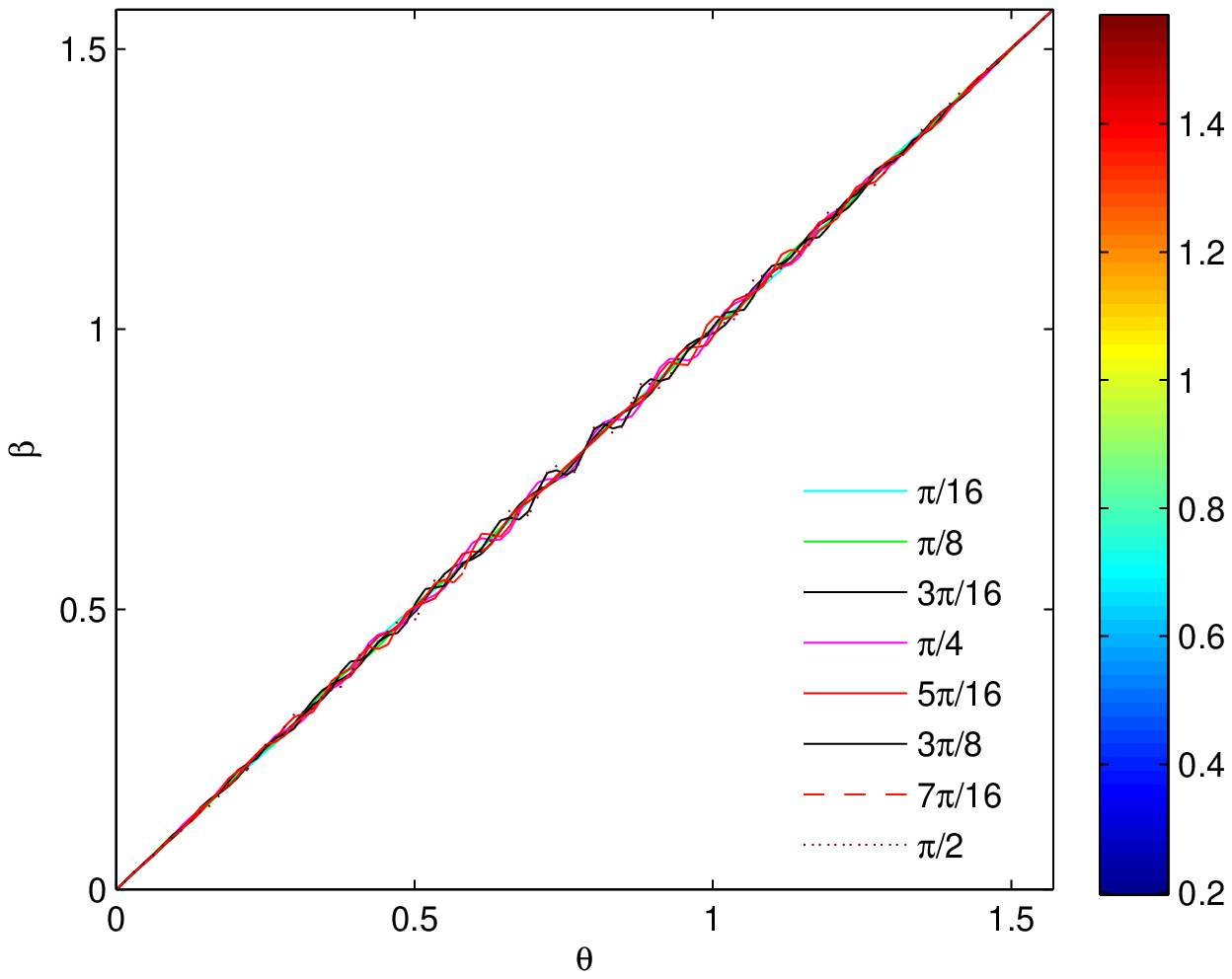}
{\footnotesize (b) $|{\bm \kappa}|=5\pi$}
\end{minipage}\\
\begin{minipage}[t]{60mm}
\includegraphics[width=60mm]{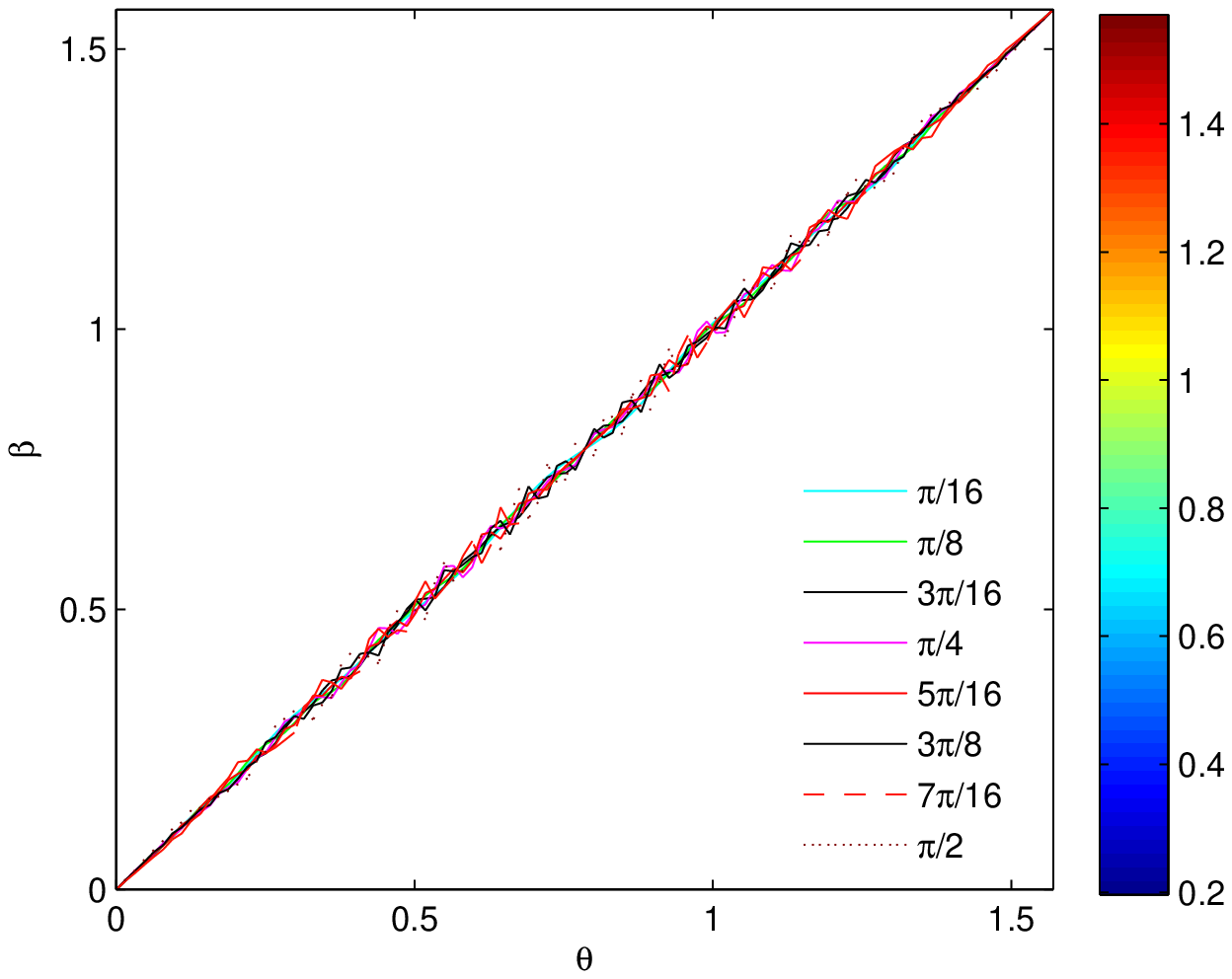}
{\footnotesize (c) $|{\bm \kappa}|=7.5\pi$}
\end{minipage}\\
\caption{Wave propagation angle $\beta$ versus the wave number angle $\theta$ at different $|{\bm \kappa}|$ with $\tau=0.01,\ h=0.1$ and $N=150$ for the Maxwell's equations.}\label{fig618}
\end{figure}

Due to the symmetry of the $\phi$ and the range value of the function arctan, we only need to consider $\phi\in[0,\frac{\pi}{2}]$, here.
The relations between the propagation angle $\alpha$ and the wave number angle $\theta\in[0,2\pi]$ with $|{\bm \kappa}|=2.5\pi,
|{\bm \kappa}|=5\pi$ and $|{\bm \kappa}|=7.5\pi$ at different $\phi\in[0,\frac{\pi}{2}]$
are plotted in Fig. \ref{fig615}.
We can find the similar results about symmetry with respect to $\theta$.
Thus, we just need
to discuss the relation between $\alpha$ and $\phi$ in the domain $[0,\frac{\pi}{2}]$ of $\theta$, which is represented by the contour plots in
Fig. \ref{fig616}.
Then, the relations between the propagation angle $\beta$ and the wave number angle $\phi\in[0,2\pi]$ with $|{\bm \kappa}|=2.5\pi,
|{\bm \kappa}|=5\pi$ and $|{\bm \kappa}|=7.5\pi$ at different $\theta\in[0,\frac{\pi}{2}]$
are plotted in Fig. \ref{fig617}, which shows that whatever the $\phi$ chosen, the $\beta$ is symmetric
with respect to $\phi=\frac{\pi}{2}$.
 This fact implies that we only need to discuss the
 relation between $\beta$ and $\theta$ in the domain $[0,\frac{\pi}{2}]$ of $\phi$,
which is displayed via contour plots in Fig. \ref{fig618}.
As illustrated in Fig. \ref{fig615}, we can observe that, for a fixed angle $\alpha$, as the wave number $|{\bm \kappa}|$ increases, $\alpha$ does vary.
In contrast, as the angle $\theta$ vary,
$\alpha$ has almost no change. This implies that the angle $\alpha$ only depends on the wave number $|{\bm \kappa}|$ and
is independent of the angle $\theta$.
Similarly, as illustrated in Fig. \ref{fig617}, we can see that the angle $\beta$ is independent of the angle $\phi$ as well.
We can conform the facts again in Fig. \ref{fig616} and Fig. \ref{fig618} respectively.
Now, we can conclude that the grid-anisotropy of the proposed method is direction-independent.

\section{Concluding remarks}

In this paper, a sixth order energy-conserved method is first developed for the 3D time-domain
Maxwell's equations. The proposed scheme can preserve all of the desired structures, including the symmetry, the five energy conservation
laws, the two divergence-free fields, the three momentum conservation laws as well as the symplectic conservation law.
An optimal error estimate for the proposed scheme is established with the best constant $O(T)$ in discrete $L^{2}$-norm.
Numerical results confirm our error estimate.
In addition, numerical dispersion analysis show that the non-physical
solution branches of the Fourier pseudo-spectral method only occur with the large wave numbers and its grid-anisotropy is direction-independent.
It is well-known that the numerical computation of Maxwell's equations in metamaterials plays a very important
role in seeking new designs and applications of metamaterials.
Therefore, the generalization of our energy-conserved method and numerical analysis for the Maxwell's equations
 in metamaterials will be the subject of our future research.

\section*{Acknowledgments}
This work is supported by the National Natural Science Foundation of China (Grant Nos.
11771213 and 41504078), the National Key Research and Development Project of China (Grant No. 2016YFC0600310),
 the Major Projects of Natural Sciences of University in Jiangsu Province of China (Grant No. 15KJA110002) and the Priority Academic Program Development of Jiangsu Higher Education Institutions.

\section*{Appendix:}
\appendix
\begin{appendices}
\section{A fast solver to the proposed scheme }
In Ref. \cite{CWG16}, an iterative method need be employed to solve the resulting linear equations.
This leads to high computational cost.
In this Appendix, a fast solver is presented to increasing computational efficiency. The idea is
based on the matrix diagonalization method (see Ref. \cite{ST06} and references therein) and the Fast Fourier Transform (FFT) algorithm,
which is different from those of the ADI-FDTD methods (e.g., see Refs. \cite{Na99,ZCZ00}) and the LOD-FDTD methods (e.g., see \cite{SMYN05}).
Below we will list the key points of this solver.

We rewrite (\ref{SO-ECS1}) as
\begin{align}\label{SX_A.2}
 &2\mu{ {\bm H}}^{n+\frac{1}{2}} +\tau{\bm D}{ {\bm E}}^{n+\frac{1}{2}}+\frac{c^{2}\tau^{3}}{12}{\bm D} ^{3}{{\bm E}}^{n+\frac{1}{2}}+\frac{c^{4}\tau^{5}}{120}{\bm D}^{5}{ {\bm E}}^{n+\frac{1}{2}}=2\mu{ {\bm H}}^{n}, \\\label{SX_A.3}
 &2\epsilon{ {\bm E}}^{n+\frac{1}{2}} -\tau{\bm D}{{\bm H}}^{n+\frac{1}{2}}-\frac{c^{2}\tau^{3}}{12}{\bm D} ^{3}{{\bm H}}^{n+\frac{1}{2}}-\frac{c^{4}\tau^{5}}{120}{\bm D}^{5}{{\bm H}}^{n+\frac{1}{2}}=2\epsilon{ {\bm E}}^{n}.
\end{align}
According to Lemma \ref{lem3.3} and $\mathcal{F}^{-1}_{N_{w}}\mathcal{F}_{N_{w}}={\bm I}_{N_{w}}$, ${\bm D}$ can be rewritten as
\begin{align}
{\bm D}={\bm F}^{-1}
{\bm \Lambda}{\bm F},
\end{align}
where
\begin{align*}
&{\bm \Lambda}=\left(\begin{array}{ccc}
             0&-\Lambda_{z}\otimes I_{N_{y}}\otimes I_{N_{x}}& I_{N_{z}}\otimes\Lambda_{y}\otimes I_{N_{x}}\\
              \Lambda_{z}\otimes I_{N_{y}}\otimes I_{N_{x}}&0&-I_{N_{z}}\otimes I_{N_{y}}\otimes \Lambda_{x} \\
             -I_{N_{z}}\otimes\Lambda_{y}\otimes I_{N_{x}}&I_{N_{z}}\otimes I_{N_{y}}\otimes \Lambda_{x}&0
             \end{array}
\right),\\
&{\bm F}=\left(\begin{array}{ccc}
              \mathcal{F}_{N_{z}}\otimes \mathcal{F}_{N_{y}}\otimes \mathcal{F}_{N_{x}}     &  \  &\   \\
              \    &  \mathcal{F}_{N_{z}}\otimes \mathcal{F}_{N_{y}}\otimes \mathcal{F}_{N_{x}} &\
               \\
              \ &\  & \mathcal{F}_{N_{z}}\otimes \mathcal{F}_{N_{y}}\otimes \mathcal{F}_{N_{x}}\\
             \end{array}
\right).
\end{align*}
Let
\begin{align}\label{fft}
 &{ \widetilde{\bm H}}_{w}^{n+\frac{1}{2}}=\mathcal{F}_{N_{z}}\otimes \mathcal{F}_{N_{y}}\otimes \mathcal{F}_{N_{x}}{\bm H}_{w}^{n+\frac{1}{2}},\
 { \widetilde{\bm H}}_{w}^{n}=\mathcal{F}_{N_{z}}\otimes \mathcal{F}_{N_{y}}\otimes \mathcal{F}_{N_{x}}{\bm H}_{w}^{n},\\
 &{ \widetilde{\bm E}}_{w}^{n+\frac{1}{2}}=\mathcal{F}_{N_{z}}\otimes \mathcal{F}_{N_{y}}\otimes \mathcal{F}_{N_{x}}{\bm E}_{w}^{n+\frac{1}{2}},\
 { \widetilde{\bm E}}_{w}^{n}=\mathcal{F}_{N_{z}}\otimes \mathcal{F}_{N_{y}}\otimes \mathcal{F}_{N_{x}}{\bm E}_{w}^{n}.
\end{align}
Left-multiplying \eqref{SX_A.2} and \eqref{SX_A.3} with ${\bm F}$, respectively, we have
\begin{align}\label{A5}
&2\mu\left(\begin{array}{ccc}
{ \widetilde{\bm H}}_{x}^{n+\frac{1}{2}}\\
{ \widetilde{\bm H}}_{y}^{n+\frac{1}{2}}\\
{ \widetilde{\bm H}}_{z}^{n+\frac{1}{2}}
\end{array}
\right)+\Big(\tau{\bm \Lambda}+\frac{c^{2}\tau^{3}}{12}{\bm \Lambda}^{3}+
\frac{c^{4}\tau^{5}}{120}{\bm \Lambda}^{5}\Big)
\left(\begin{array}{ccc}
{ \widetilde{\bm E}}_{x}^{n+\frac{1}{2}}\\
{\widetilde{\bm E}}_{y}^{n+\frac{1}{2}}\\
{\widetilde{\bm E}}_{z}^{n+\frac{1}{2}}
\end{array}
\right)=2\mu
\left(\begin{array}{ccc}
{ \widetilde{\bm H}}_{x}^{n}\\
{\widetilde{\bm H}}_{y}^{n}\\
{\widetilde{\bm  H}}_{z}^{n}
\end{array}
\right),\\\label{A6}
&2\epsilon\left(\begin{array}{ccc}
{ \widetilde{\bm E}}_{x}^{n+\frac{1}{2}}\\
{ \widetilde{\bm E}}_{y}^{n+\frac{1}{2}}\\
{ \widetilde{\bm E}}_{z}^{n+\frac{1}{2}}
\end{array}
\right)-\Big(\tau{\bm \Lambda}+\frac{c^{2}\tau^{3}}{12}{\bm \Lambda}^{3}+
\frac{c^{4}\tau^{5}}{120}{\bm \Lambda}^{5}\Big)
\left(\begin{array}{ccc}
{ \widetilde{\bm H}}_{x}^{n+\frac{1}{2}}\\
{ \widetilde{\bm H}}_{y}^{n+\frac{1}{2}}\\
{ \widetilde{\bm H}}_{z}^{n+\frac{1}{2}}
\end{array}
\right)=2\epsilon
\left(\begin{array}{ccc}
{ \widetilde{\bm E}}_{x}^{n}\\
{ \widetilde{\bm E}}_{y}^{n}\\
{\widetilde{\bm E}}_{z}^{n}
\end{array}
\right).
\end{align}
Denoting
\begin{align}\label{cfs}
&
{\bm a}_{j,k,m}=\Bigg(\begin{array}{cccc}
 0&-\Lambda_{z_{mm}}& \Lambda_{y_{kk}}\\
              \Lambda_{z_{mm}}&0&-\Lambda_{x_{jj}} \\
             -\Lambda_{y_{kk}}&\Lambda_{x_{jj}}&0
             \end{array}
\Bigg),\
\end{align}
Eqs. \eqref{A5}-\eqref{A6} can be rewritten into the following subsystems
\begin{align}\label{SY1}
{\bm A}_{j,k,m}\left(\begin{array}{c}
{\widetilde{\bm H}}_{x_{j,k,m}}^{n+\frac{1}{2}}\\
{\widetilde{\bm H}}_{y_{j,k,m}}^{n+\frac{1}{2}}\\
{\widetilde{\bm H}}_{z_{j,k,m}}^{n+\frac{1}{2}}\\
{\widetilde{\bm E}}_{x_{j,k,m}}^{n+\frac{1}{2}}\\
{\widetilde{\bm E}}_{y_{j,k,m}}^{n+\frac{1}{2}}\\
{\widetilde{\bm E}}_{z_{j,k,m}}^{n+\frac{1}{2}}
\end{array}
\right)={\bm B}_{j,k,m}\left(\begin{array}{c}
{\widetilde{\bm H}}_{x_{j,k,m}}^{n}\\
{\widetilde{\bm H}}_{y_{j,k,m}}^{n}\\
{\widetilde{\bm H}}_{z_{j,k,m}}^{n}\\
{\widetilde{\bm E}}_{x_{j,k,m}}^{n}\\
{\widetilde{\bm E}}_{y_{j,k,m}}^{n}\\
{\widetilde{\bm E}}_{z_{j,k,m}}^{n}
\end{array}
\right),\left. \begin{array}{c}
j=1,\cdots, N_{x},\\
k=1,\cdots,
N_{y},\\m=1,\cdots, N_{z},
\end{array}
\right.
\end{align}
where
\begin{align*}
{\bm A}_{j,k,m}=
\left(\begin{array}{cccc}
2\mu{\bm I}_{3}&\bar{\bm a}_{j,k,m}\\
 -\bar{\bm a}_{j,k,m}&2\epsilon{\bm I}_{3}
\end{array}
\right),~{\bm B}_{j,k,m}=
\left(\begin{array}{cccc}
2\mu{\bm I}_{3}&\ \\
 \ &2\epsilon {\bm I}_{3}
\end{array}
\right),
\end{align*}
and
\begin{align*}
 \bar{\bm a}_{j,k,m}=\tau{\bm a}_{j,k,m}+\frac{c^{2}\tau^{3}}{12}{\bm a}_{j,k,m}^{3}+\frac{c^{4}\tau^{5}}{120}{\bm a}_{j,k,m}^{5}.
 \end{align*}
Since ${\bm A}_{j,k,m}$ is a constant nondegenerate matrix (see Remark \ref{A1}), for a fixed $(j,k,m)$, we have
\begin{align}\label{SY2}
\ \left(\begin{array}{c}
{\widetilde{\bm H}}_{x_{j,k,m}}^{n+\frac{1}{2}}\\
{\widetilde{\bm H}}_{y_{j,k,m}}^{n+\frac{1}{2}}\\
{\widetilde{\bm H}}_{z_{j,k,m}}^{n+\frac{1}{2}}\\
{\widetilde{\bm E}}_{x_{j,k,m}}^{n+\frac{1}{2}}\\
{\widetilde{\bm E}}_{y_{j,k,m}}^{n+\frac{1}{2}}\\
{\widetilde{\bm E}}_{z_{j,k,m}}^{n+\frac{1}{2}}
\end{array}
\right)={\bm A}_{j,k,m}^{-1}{\bm B}_{j,k,m}\left(\begin{array}{c}
{\widetilde{\bm H}}_{x_{j,k,m}}^{n}\\
{\widetilde{\bm H}}_{y_{j,k,m}}^{n}\\
{\widetilde{\bm H}}_{z_{j,k,m}}^{n}\\
{\widetilde{\bm E}}_{x_{j,k,m}}^{n}\\
{\widetilde{\bm E}}_{y_{j,k,m}}^{n}\\
{\widetilde{\bm E}}_{z_{j,k,m}}^{n}
\end{array}
\right):={\mathcal A}_{j,k,m}\left(\begin{array}{c}
{\widetilde{\bm H}}_{x_{j,k,m}}^{n}\\
{\widetilde{\bm H}}_{y_{j,k,m}}^{n}\\
{\widetilde{\bm H}}_{z_{j,k,m}}^{n}\\
{\widetilde{\bm E}}_{x_{j,k,m}}^{n}\\
{\widetilde{\bm E}}_{y_{j,k,m}}^{n}\\
{\widetilde{\bm E}}_{z_{j,k,m}}^{n}
\end{array}
\right),
\end{align}
Solving the above equations gives ${\widetilde{\bm H}}_w^{n+\frac{1}{2}}$
and ${\widetilde{\bm E}}_w^{n+\frac{1}{2}}$. Using the relations ${ {\bm H}}_{w}^{n+1}=2\mathcal{F}_{N_{z}}^{-1}\otimes \mathcal{F}_{N_{y}}^{-1}\otimes \mathcal{F}_{N_{x}}^{-1}{ \widetilde{\bm H}}^{n+\frac{1}{2}}_{w}-{ {\bm H}}^{n}_{w}$ and ${ {\bm E}}_{w}^{n+1}=2\mathcal{F}_{N_{z}}^{-1}\otimes \mathcal{F}_{N_{y}}^{-1}\otimes \mathcal{F}_{N_{x}}^{-1}{ \widetilde{\bm E}}^{n+\frac{1}{2}}_{w}-{ {\bm E}}^{n}_{w}$ yields the solutions
${ {\bm H}}_{w}^{n+1}$ and ${ {\bm E}}_{w}^{n+1}$.
Note that the matrix $\mathcal{A}$ in Eq. \eqref{SY2} need only be computed once  and the Fast Fourier Transform (FFT) algorithm can be applied to the above process.
\begin{rmk}\label{A1} Note that there exists an orthogonal matrix ${\bm O}$ such that
\begin{align*}
{\bm a}_{12}={\bm O}^{T}{ \widetilde{\bm\Lambda}}{\bm O},\  \widetilde{\bm \Lambda}=\text{\rm diag}\Big(\lambda_1,\lambda_2,\lambda_3\Big),\ \lambda_1=0,\ \lambda_{2,3}=\pm\sqrt{-\Lambda_{x_{jj}}^2-\Lambda_{y_{kk}}^2-\Lambda_{z_{mm}}^2},
\end{align*}
which further implies that
\begin{align*}
\bar{{\bm a}}_{12}={\bm O}^{T}\widetilde{\bm \Lambda}_{\tau}{\bm O},\ \widetilde{\bm \Lambda}_{\tau}=\text{\rm diag}\Big(\widetilde{\lambda}_1,\widetilde{\lambda}_2,\widetilde{\lambda}_3\Big),
\widetilde{\lambda}_1=0,\ \widetilde{\lambda}_{2,3}=\tau\lambda_{2,3}+\frac{c^2\tau^3}{12}\lambda_{2,3}^3+\frac{c^4\tau^5}{120}\lambda_{2,3}^5.
\end{align*}
Thus, it is clear to see that
\begin{align*}
|{\bm A}|&=\Big|\left(\begin{array}{cccc}
{\bm O}^{T}&\ \\
 \ &{\bm O}^{T}
\end{array}
\right)
\left(\begin{array}{cccc}
2\mu{\bm I}_{3\times 3}&\widetilde{\bm \Lambda}_{\tau}\\
 -\widetilde{\bm \Lambda}_{\tau}&2\epsilon{\bm I}_{3\times 3}
\end{array}
\right)\left(\begin{array}{cccc}
{\bm O}&\ \\
 \ &{\bm O}
\end{array}
\right)\Big|\nonumber\\
&=\Big|
\left(\begin{array}{cccc}
2\mu{\bm I}_{3\times 3}&\widetilde{\bm \Lambda}_{\tau}\\
 \ &2\epsilon{\bm I}_{3\times 3}+\frac{1}{2\mu}\widetilde{\bm \Lambda}_{\tau}^2
\end{array}
\right)\Big|\not=0.
\end{align*}
\end{rmk}
\end{appendices}


\end{document}